\documentclass[11pt,leqno]{amsart}

\usepackage{amssymb}
\usepackage{mathrsfs}
\usepackage{color}
\usepackage[latin1]{inputenc}
\usepackage[T1]{fontenc}
\usepackage{hyperref}
\usepackage{graphicx}
\usepackage{tikz}
\usetikzlibrary{decorations.pathmorphing}

\newcommand\E{{\mathbb E}}
\newcommand\N{{\mathbb N}}
\newcommand\R{{\mathbb R}}

\newcommand\Sp{{\mathbb S}}

\def\BB{{\mathcal B}}

\def\DD{{\mathcal D}}
\def\EE{{\mathcal E}}
\def\FF{{\mathcal F}}
\def\GG{{\mathcal G}}
\def\HH{{\mathcal H}}
\def\II{{\mathcal I}}

\def\LL{{\mathcal L}}

\def\NN{{\mathcal N}}
\def\OO{{\mathcal O}}
\def\PP{{\mathcal P}}

\def\SS{{\mathcal S}}
\def\TT{{\mathcal T}}
\def\UU{{\mathcal U}}
\def\VV{{\mathcal V}}

\def\VV{{\mathcal V}}

\def\ZZ{{\mathcal Z}}

\def\PPS{{\mathcal{P}}_{\mbox{{\tiny sym}}}}

\def\eps{{\varepsilon}}

\newcommand{\wto}{\rightharpoonup}

\def\SN{\mathfrak{S}_N}
\def\SSS{\mathfrak{S}}

\newtheorem{theo}{Theorem}

\newtheorem{lem}[theo]{Lemma}

\newtheorem{rem}[theo]{Remark}

\newtheorem{defin}[theo]{Definition}
\newtheorem{ex}[theo]{Example}

\newcommand{\beqn}{\begin{equation}}
\newcommand{\eeqn}{\end{equation}}
\newcommand{\bear}{\begin{eqnarray}}
\newcommand{\eear}{\end{eqnarray}}
\newcommand{\bean}{\begin{eqnarray*}}
\newcommand{\eean}{\end{eqnarray*}}

\newcommand{\Black}{\color{black}}

\def\signsm{\bigskip \begin{center} {\sc St\'ephane Mischler\par\vspace{3mm}
Universit\'e Paris-Dauphine\par
CEREMADE, UMR CNRS 7534\par
Place du Mar\'echal de Lattre de Tassigny
75775 Paris Cedex 16\par
FRANCE\par\vspace{3mm}
e-mail:} \tt{mischler@ceremade.dauphine.fr} \end{center}}

\def\signcm{\bigskip \begin{center} {\sc 
Cl\'ement Mouhot\par\vspace{3mm}
University of Cambridge\par
DPMMS, Centre for Mathematical Sciences\par
Wilberforce Road, 
Cambridge CB3 0WA, 
UK\par\vspace{3mm}
e-mail:} \tt{C.Mouhot@dpmms.cam.ac.uk} \end{center}}

\def\signbw{\bigskip \begin{center} {\sc 
Bernt Wennberg\par\vspace{3mm}
Department of Mathematical Sciences\par
Chalmers University of Technology \par
and\par
Department of Mathematical Sciences\par
University of Gothenburg\par
41296 G\"oteborg\par
SWEDEN\par\vspace{3mm}
e-mail:} \tt{wennberg@chalmers.se} \end{center}}

\begin{document}

\title[New approach to quantitative propagation of chaos]{A new
  approach to quantitative propagation of chaos for drift,
  diffusion and jump processes}

\author{S. Mischler}
\author{C. Mouhot}
\author{B. Wennberg}

\maketitle


\begin{abstract} This paper is devoted the the study of the mean field
  limit for many-particle systems undergoing jump, drift or diffusion
  processes, as well as combinations of them. The main results are
  quantitative estimates on the decay of fluctuations around the
  deterministic limit and of correlations between particles, as the
  number of particles goes to infinity. To this end we introduce a
  general functional framework which reduces this question to the one
  of proving a purely functional estimate on some abstract generator
  operators ({\em consistency estimate}) together with fine stability
  estimates on the flow of the limiting nonlinear equation ({\em
    stability estimates}). Then we apply this method to a Boltzmann
  collision jump process (for Maxwell molecules), to a McKean-Vlasov
  drift-diffusion process and to an inelastic Boltzmann collision jump
  process with (stochastic) thermal bath.  To our knowledge, our
  approach yields the first such quantitative results for a
  combination of jump and diffusion processes.
  \end{abstract}

\vspace{0.3cm}

 \smallskip


\medskip

\textbf{Acknowledgments}: B.W. would like the CEREMADE at University
Paris-Dauphine for the invitation in june and october 2006 where this
work was initiated. S.M. and C.M. would like to thank the mathematics
departement of Chalmers University for the invitation in november
2008. The authors also thank F. Bolley, J. A. Ca\~nizo, N. Fournier,
A. Guillin, J. Rousseau and C. Villani for fruitful
discussions. The authors  also wish to mention the inspirative
courses of P.-L. Lions at Coll\`ege de France on ``Mean Field Games''
in 2007-2008 and 2008-2009. The authors thank the anonymous referee
for many helpful critics and remarks. 

\newpage

\tableofcontents



\section{Introduction}
\label{sec:intro}
\setcounter{equation}{0}
\setcounter{theo}{0}



\medskip

Fundamental in Boltzmann's deduction of the equation bearing his name
is the ``\emph{stosszahlansatz''}, or chaos assumption. Vaguely
expressed this assumption means that when two particles collide, they
are statistically uncorrelated just before the collision. After the
collision they are not, of course, because, for example, the knowledge
of the position and velocity of one particle that just collided gives
some information on the position of the collision partner. And while
the correlations created by collisions decrease with time, they never
vanish in a system of finitely many particles, and hence the Boltzmann
assumption could only be true in the limit of infinitely many
particles.

A mathematical framework for studying this limit is the so-called
BBGKY hierarchy (see Grad \cite{MR0135535}, Cercignani
\cite{MR0449375} and the book \cite{MR1307620} by Cercignani et al),
which consists of a family of Liouville equations, each describing the
evolution of an $N$-particle system (deterministic, in this case), and
whose solutions are densities in the space of $N$-particle
configurations in phase space. The BBGKY hierarchy describes in a
systematic way the evolution of marginal distributions. Formally, and
under appropriate assumptions, most notably the chaos assumption, the
one-particle marginal of solutions to the $N$-particle Liouville
equation, converges to solutions of the Boltzmann equation.

It would take until 1974 before a mathematically rigorous proof of
this statement was given by Lanford \cite{Lanford}. While this is
a remarkable result, it only proves that the Boltzmann equation is a
limit of the $N$-particle systems for a fraction of the mean
free time between collisions, and this is essentially where the
problem stands today (see however \cite{MR849204} for a large
time limit in a near to the vacuum framework; it is worth
emphasizing that in such a framework no more collisions occur than
in Lanford's framework).

In order to avoid some of the difficulties related to the
deterministic evolution of a real particle system, Kac \cite{Kac1956}
invented a Markov process for a particular $N$-particle system, and
gave a mathematically rigorous definition of {\em propagation of
  chaos}. He then proved that this holds for his Markov system, and
thus obtained a mathematically rigorous derivation of a simplified
(spatially homogeneous) Boltzmann equation in this case, usually
called the \emph{Kac equation}.

Kac's work provides the framework of this paper, and we will now
describe our main results. We let $E$ be the state space of one
particle (usually $\R^d$, but metric, separable and locally compact is
fine). A sequence of probability measures ${(f^N)}_{N=0}^\infty$,
where each $f^N\in\PP(E^N)$ is \emph{symmetric} in the sense that it
is invariant under permutation of the coordinates, is said to be {\em
  $f$-chaotic} for some probability measure $f\in\PP(E)$ if for each
$k\ge 1$ and functions $\phi_j\in C_b(E), j=1, \dots,k$, (continuous
bounded),
\begin{eqnarray}
\label{eq:0010}
  \lim_{N\rightarrow\infty} \int_{E^N} \prod_{j=1}^{k} \phi_j(z_j) \, 
  f^N({\rm d}z_1,\dots,{\rm d}z_N) &=&  \prod_{j=1}^{k}
  \int_{E} \phi_j(z) \, f({\rm d}z).
\end{eqnarray}
We next consider a family of time dependent probability measures
${(f_t^N)}_{N=0}^\infty$, being the distributions of the states of
Markov processes in $E^N$. The Markov process is said to {\em
  propagate chaos} if given an initial family of $N$-particle
distributions ${(f_{\mathit{in}}^N)}_{N =0}^\infty$, that is
$f_{\mathit{in}}$-chaotic, there is a time dependent distribution
$f_t$ such that ${(f_t^N)}_{N=1}^{\infty}$ is $f_t$-chaotic. In this
paper we are interested in specific equations that govern the
evolution of $f_t$; but it is important to bear in mind that they are
in general nonlinear, and it may be difficult to prove well-posedness
in function spaces relevant for proving the propagation of chaos.
\medskip

The main results of this paper are abstract. We consider:
\begin{itemize}
\item the family of $N$-particle systems represented by a Markov
  processes $(\ZZ^N_t)_{t\ge 0}$ in some product space $E^N$, with
  $\ZZ^N_t=(\ZZ_{1,t},\dots,\ZZ_{N,t})$, and the corresponding
  probability distributions ${(f_t^N)}_{N=1}^{\infty}$, solving
  Kolmogorov's backward equations\footnote{As we will see in the
    following section, the formalism also works if $\ZZ^N_t$ satisfies
    a deterministic evolution, in which case Liouville's equations
    replace Kolmogorov's backward equation}, and where the
  distributions $f_t^N$ belong to a suitable subspace of $\PP(E^N)$,
\item a (nonlinear) equation defined on a subspace of $\PP(E)$,
  which is the formal limit of the equations governing one-particle
  marginals of $f^N_t$:
\begin{eqnarray}
\label{eq:0020a}
  \frac{\partial}{\partial t} f_t &=& Q(f_t), \qquad   f_{\mbox{\tiny in}} \in \PP(E).
\end{eqnarray}
\end{itemize}
Then we:
\begin{itemize}
\item provide conditions on the processes and related function spaces
  that guarantee that ${(f_t^N)}_{N=1}^{\infty}$ is $f_t$-chaotic for
  $t \ge0$,
\item give explicit estimates of the rate of convergence
  in~(\ref{eq:0010}): more precisely, for any $T>0$,  $\ell \in \N^*$ and
  $\phi_j\in\FF\subset C_b(E)$ ($j=1,\dots,\ell$), there is a constant
  $\epsilon(N)$ converging to zero as $N\rightarrow\infty$ such that 
  \begin{equation}
    \label{eq:0015}
    \sup_{t \in [0,T]}
    \left| \int_{E^N} \prod_{j=1}^{\ell} \phi_j(z_j) \,
      f^N_t({\rm d}z_1,\dots,{\rm d}z_N)-  \prod_{j=1}^{\ell}
      \int_{E}\phi_j(z) \, f_t({\rm d}z)
    \right| \le \epsilon(N),
  \end{equation}
 which holds for $N\ge2\ell$
 and a suitably chosen space $\FF$; if $\FF$ is dense in $C_b(E)$,
 this implies in particular the propagation of chaos.
\end{itemize}
\medskip

To this end, our starting point is a technique that goes back at least to
Gr¸nbaum~\cite{Grunbaum}, which consists in representing an
$N$-particle configuration $\ZZ^N_t$ as a sum of Dirac measures,
\begin{eqnarray*}
  \ZZ^N_t=(\ZZ_{1,t},\dots,\ZZ_{N,t}) &\quad \longleftrightarrow \quad
  &
  \mu^N_{\ZZ^N_t}=\frac{1}{N}\sum_{j=1}^N \delta_{\ZZ_{j,t}} \in \PP(E)
\end{eqnarray*}
and proving that, in a weak sense, $\mu^N_{\ZZ^N_t}$ converges to
$f^N_t$. In fact, because $\ZZ^N_t$ is random, $\mu^N_{\ZZ^N_t}$ is a
random measure in $\PP(E)$, which has a probability distribution
$\Psi^N_t\in \PP(\PP(E))$. Proving the propagation of chaos is here
equivalent to proving that $\Psi^N_t\rightarrow \delta_{f_t}$ in
$\PP(\PP(E))$ when $N\rightarrow\infty$. The error $\epsilon(N)$ is
dominated by, on the one hand, how well the initial measure
$f_{\mbox{\tiny in}}$ can be approximated by a sum of $N$ Dirac measures,
and, on the other hand, estimates comparing the equations for $f^N_t$
and $f_t$.  These estimates depend on rather technical assumptions,
and although the abstract main theorem is stated in
Section~\ref{sec:abstractchaos}, the assumptions are stated in full
detail only in Section~\ref{sec:lemmas}.

\medskip

With the main theorem of this paper in hand, proving the propagation
of chaos for a particular $N$-particle system is reduced to proving:
\begin{itemize}
\item[(i)] a purely functional estimate on the dual generator $G^N$ of
  the $N$-particle dynamics which establishes and quantifies that, at 
  first order, $G^N$ is linked to the mean field limit generator $Q$ 
  (\emph{consistency estimate});
\item[(ii)] some fine stability estimates on the flow of the mean
  field limit equation involving the differential of the semigroup
  with respect to the initial data (\emph{stability estimates}).
\end{itemize}

\medskip

Point (i) of our method is largely inspired from the ``duality
viewpoint" of Gr\"unbaum's paper \cite{Grunbaum} where he considered
the propagation of the chaos issue for the Boltzmann equation
associated to hard-spheres (unbounded) kernel. As he confessed himself
the proof in \cite{Grunbaum} was incomplete due to the lake of
suitable stability estimates, i.e. precisely the point (ii) of our
method.

It is worth emphasizing that after we had finished writing our paper,
we were told about the recent book \cite{Kolokoltsov} by Kolokoltsov and
his series of papers on nonlinear Markov processes and kinetic
equations.  These interesting works focus on fluctuation estimates of
LLN and CLT types in the general framework of nonlinear Markov
processes, and in some sense they generalise to several other kinetic
models the Gr\"unbaum's duality viewpoint (although Kolokoltsov seems
to not be aware of that earlier work).  However 
we were not able to extract from these works a full proof in the cases
when the generator is an unbounded operator and weak distances have to
be used.  While the comparison of generators for the many-particle and
the limit semigroup present in both \cite{Grunbaum} and
\cite{Kolokoltsov} is reminiscent of our work, we believe that the main
novelty of the present paper is to achieve, for the first time, both
the fine stabilities estimates in point (ii) and the consistency
estimate in point (i) in appropriate spaces (with weak topologies), in
such a way that they may be combined and they lead to the already
mentioned propagation of chaos result with quantitative estimates.

\medskip

We illustrate the method by proving the propagation of chaos for three
different well known examples:

\begin{itemize}
\item [(a)] We first consider \emph{the Boltzmann equation for Maxellian molecules with
  angular cutoff}. For such a bounded kernel case the result is
  well-known since the pioneering works of Kac~\cite{Kac1956,Kac1957}
  and McKean \cite{McKean1967} (who prove the propagation of chaos
  without any rate) and from the works by Graham and M\'el\'eard
  \cite{GrahamM1,GrahamM2,GrahamM3,Meleard1996} (where the authors
  establish the propagation of chaos with optimal rate $\OO(1/N)$). In
  these papers, the cornerstone of the proof is a combinatorial
  argument applied to the equation on the law (Wild sum expansion) or
  to the stochastic flow (stochastic tree). These approaches are
  restricted to a constant (or at least bounded) collision rate.
\item[(b)] The second example is the \emph{McKean-Vlasov model}. For
  such a model again, propagation of chaos is well-known and has been
  extensively studied. One of the most popular and efficient approach
  to deal with this model is the so-called ``coupling method''
  introduced in the 1970's, which yields the optimal convergence rate
  $\OO(1/\sqrt{N})$ (note that the difference between these two
  optimal rates in (a) and (b) comes from the fact that they are not
  measured with the same distance). We refer to the lecture notes
  \cite{S6,Meleard1996} as well as to the references therein for a
  detailed discussion of that method.  We also refer to
  \cite{BolleyGM} and the references therein for recent developments
  on the subject.
\item[(c)] The third example is a \emph{mixed collision-diffusion
    equation} which arises from granular gas modeling. For such a
  model, it seems that both the ``combinatorics method'' and the
  ``coupling method'' fail while our present method is robust enough
  to apply and yield quantitative chaos estimates. Let us also
  emphasize that the BBGKY method and the nonlinear martingale method
  (see again \cite{S6,Meleard1996} or \cite{ArkCapIan,SMeddimo}) may
  also apply but would give a propagation of chaos without any rate
  since they are based on compactness arguments.
\end{itemize}

\medskip

Let us emphasize that it is not difficult to write a \emph{uniform in
  time} version of Theorem 2.1: in short, if the assumptions {\bf
  (A1)} to {\bf (A5)} are satisfied with $T=+\infty$, then the
conclusion of the main abstract Theorem~\ref{theo:abstract} holds with
$T=+\infty$ and the proof is unchanged. But such an abstract theorem
does not readily apply to the examples (a), (b) and (c) discussed
above. More precisely, it is indeed possible to prove quantitative
uniform in time propagation of chaos by our method (for the elastic
Boltzmann model for instance), but the price to pay is a significant
modification to the set of assumptions {\bf (A1)} to {\bf (A5)}. This
issue is addressed in our companion paper \cite{MM**} where the
abstract method is developed in a more general framework in order to
(1) apply it to Boltzmann collision models associated to unbounded
collisions rates, (2) develop a theory of uniform in time propagation
of chaos estimates. We shall consider the question of uniform in time
chaoticity estimates for the McKean-Vlasov equation in future works. 

\smallskip 
These three examples illustrate the generality of the method that we
study: the same abstract framework can be used to prove propagation of
chaos for $N$-particle systems that have not yet been analysed as well
as for models that have been studied before but with conceptually
different methods. But we chose to emphasize generality over optimality of the result. 
By optimizing the method of proof for a specific problem one can certainly obtain sharper
results, for example in terms of the rate of convergence as a function
of $N$ or in the choice of topologies for which convergence can be
proven. We did not pursue this goal in this paper. 

Also,  in applying the  abstract theorem  to a concrete model,
one is faced with the challenge of finding functional spaces that
satisfy the conditions of our abstract convergence theorem and are adapted to the
model. In many cases, like for the three examples presented here,
existing theory for the $N$-particle systems and for the limiting
equations may give a strong hint on what choices to make, but in other
cases this could present serious difficulties. Another guiding principle is
the consistency estimates between the generators of the $N$-particle system 
and the limit equation which constrain the norms or metrics that can be used 
and hint at the losses on the norms or metrics at the basis of the scale of spaces
used in the stability estimates. 
 
In spite of a cost of technicality, the
original approach proposed by Grunbaum, even if originally incomplete, seems to us 
intuitively very attractive, and with
Theorem~\ref{theo:abstract},  Theorem~\ref{theo:BddBoltz}, Theorem~\ref{theo:Vlasov} and Theorem~\ref{theo:Mix} we make this approach into a mathematically rigorous theory.

\smallskip The plan of the paper is as follows.  In
Section~\ref{sec:abstractchaos}, we present the method in an abstract
framework, by first setting up a functional framework that is
appropriate for comparing the $N$-particle dynamics with the limiting
dynamics, and we establish the abstract quantitative propagation of
chaos (Theorem~\ref{theo:abstract}). The main steps of the proof are
given as well, but these rely on  some technical assumptions and
lemmas which are postponed to Section~\ref{sec:lemmas}.  The
functional framework is developed with the necessary details in
Section~\ref{sec:topology}, where we also develop a differential
calculus for functions on $\PP(E)$, as needed for studying the
nonlinear semigroup. In Section~\ref{sec:BddBoltzmann}, we apply the
method to the Boltzmann equation associated to the Maxwell molecules
collision kernel with Grad's cut-off. In Section~\ref{sec:McK}, we
apply the method to the McKean-Vlasov equation, and finally, in
Section~\ref{sec:thermostat}, it is applied to some mixed jump and
diffusion equations motivated by granular gases.

\section{Propagation of chaos for abstract $N$-particle systems}
\label{sec:abstractchaos}
\setcounter{equation}{0}
\setcounter{theo}{0}

In this section we introduce the mathematical notation used in the
paper, make precise statements of the results, and describe the main
steps of the proof, leaving the details of the proofs to the next
sections. 


\medskip
\begin{figure}[t!]
\centering

\begin{tikzpicture}
[space/.style={minimum width=8em,align=center}]
\matrix[column sep={0.15\textwidth}, row sep={0.08\textwidth}]
{
\node(a)  [space] {$E^N/\SN$} ; & \node(b)  [space] {
$\PP_{\mbox{{\scriptsize sym}}}(E^N)$}; & \node(c)  [space] {$C(E^N)$}; \\
\node(d) [space]{ $\PP_N(E)\subset \PP(E)$}; & \node(e)
[space]{$\PP(\PP(E))$};  <&
\node(f) [space]{$C(\PP(E))$};\\
};
\begin{scope}
\draw [decorate,decoration=snake,->] (a) -- node[above,midway]
{master eq.} (b);
\draw [dashed,->] (b) -- node[above,midway]  {duality} (c); 
\draw [->] (a) -- node[left,midway] {$\mu_Z^N$} (d);
\draw[decorate,decoration=snake,->] (d) -- (e);
\draw[->] (b) -- node[left,midway] {$\bar \mu_{f^N}^N$} (e);
\draw [dashed,->] (e) -- node[above,midway]  {duality} (f);
\draw[->] (f) to [bend left] node [left, midway] {$\pi^N$} (c);
\draw[->] (c) to [bend left] node [right, midway] {$R$} (f);
\draw[->] (a) to [loop above ] node (a1) [above] {$\ZZ_t^N$} (a);
\draw[->] (b) to [loop above ] node (b1) [above] {$S_t^N$} (b);
\draw[->] (c) to [loop above ] node (c1) [above] {$T_t^N\big| G^N$} (c);
\draw[->] (d) to [loop below ] node (d1) [below] {$S^{N\!L}_t$} (d);
\draw[->] (f) to [loop below ] node (f1) [below]
{$T_t^{\infty}\big| G^{\infty}$} (f);
\draw[->] (a1) to (b1);
\draw[->] (b1) to (c1);
\draw[dashed,->] (d1) to node [below, midway] {pullback} (f1);
\end{scope}
\end{tikzpicture}
\caption{A summary of spaces and their relations. Semigroups are in
  most cases given together with their generators, as in $S^N_t \big
| A$. } 
\label{fig:diagram}
\end{figure}
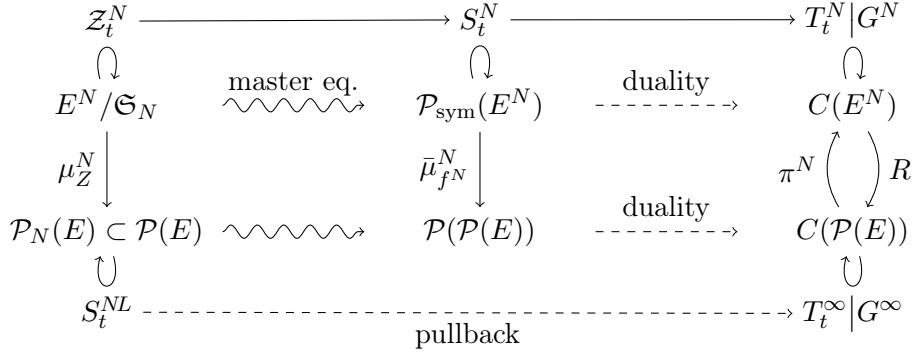

\subsection{The $N$-particle system} 
\label{sec:nparticlesystem}

The phase space of the $N$-particle system is\footnote{The phase space
  of a realistic $N$-particle system may be a subspace of the $N$-fold
  product, determined e.g. by energy constraints or by the fact that
  particles of finite size may not overlap. However, in the limit
  $N\rightarrow\infty$, all of $E$ should be accessible for any given
  particle.} $E^N/\SN$. Here $E$ is assumed to be a {\em locally
  compact, separable metric space}, and $\SN$ denotes the symmetric
group of order $N$. This means that we identify all points in the
$N$-particle phase space that can be obtained by permutation of the
particles, so that if $Z=(z_1,\dots,z_N) \in E^N/\SN$ we have
$(z_1,\dots,z_N)\sim (z_{\sigma_1},\dots,z_{\sigma_N})$, where
$(\sigma_1,\dots,\sigma_N)$ is any permutation of $\{1, \dots, N\}$.
The evolution in phase space may be a stochastic Markov process or the
solution to an Hamiltonian system of equations. In both cases, we
denote by ${(\ZZ^N_t)}_{t\ge0}$ the flow of the process.

Figure~\ref{fig:diagram} illustrates the relation between the
different objects that we consider. The $N$-particle system is
represented in the upper left corner of Figure~\ref{fig:diagram}. The
different mathematical objects in this diagram are explained along the
following subsections. 

\subsection{Master equations, Liouville's equations and their duals}
\label{sec:masterequations}

Let $\PPS(E^N)$ denote the proability measures on $E^N$ that are
invariant under permutation of the indices in $Z=(z_1,\dots,z_N)\in
E^N$.  The flow $\ZZ^N_t$ induces a semigroup of operators $S^N_t$ on
$\PPS(E^N)$ defined through the formula 
\begin{eqnarray}\label{eq:0030} 
&&
\forall \, f^N_{\mbox{{\tiny in}}} \in \PPS(E^N), \ \varphi \in C_b(E^N), \\ \nonumber
&&\qquad \left \langle S^N _t (f^N_{\mbox{{\tiny in}}} ), \varphi \right \rangle = \E
\left(\varphi\left(\ZZ^N _t\right)\right) := \int_{E^N} \E_{Z_0}
\left(\varphi\left(\ZZ^N _t\right)\right) \, f^N_{\mbox{{\tiny in}}} ({\rm d}Z_{\mbox{{\tiny in}}}),
\end{eqnarray}
where the bracket denotes the duality bracket between $\PP(E^N)$ and
$C_b(E^N)$: 
\begin{eqnarray*}
    \langle f, \phi \rangle &=& \int_{E^N} \phi(Z) \, f({\rm d}Z),
  \end{eqnarray*} 
and $\E_{Z_{\mbox{{\tiny in}}} }$ denotes the conditional expectation
with respect to the initial condition $\ZZ^N_{\mbox{{\tiny in}}}  = Z_{\mbox{{\tiny in}}} $. This semigroup
is the solution to Kolmogorov's forward equation in the case where
$(\ZZ^N_t)_{t\ge0}$ is a random process (this equation is often called
the {\em master equation}), and of the Liouville equation in the
Hamiltonian case. We always assume that $S^N_t$ preserves the
symmetry under permutation, and therefore restricts to an evolution
semigroup on $\PPS(E^N)$.
There is a dual semigroup of $S^N_t$, that acts on $C_b(E^N)$, the set
of bounded continuous functions on $E^N$.  We write this semigroup
$T^N_t$, and denote its generator by $G^N$.  The two semigroups are
related by
\begin{equation}\label{eq:0031}
  \left\langle f^N, T^N_t\phi \right\rangle = \left\langle S^N_t f^N, \phi \right\rangle.
  \end{equation}
  
  Markov processes such that (1) $T^N_t$ is a contraction in
  $C_0(E^N)$, the set of continuous functions that vanish at infinity,
  and (2) $t\mapsto T^N_t\phi$ is continuous for any $\phi \in
  C_0(E^N)$, are known as \emph{Feller processes}.

  Hence the upper part of the diagram represents a (random) process,
  and its distributions. The lower part of the diagram essentially
  shows the same thing as induced by the map $\mu^N_Z$, as we shall
  now see.

\subsection{The limiting dynamics}  
\label{sec:limitingdynamcis}

The components $(z_1,\dots,z_N)$ of $Z \in E^N / \SN$ represent the
positions (in generalized sense, i.e. in the phase space $E$) of the
$N$ particles. These $N$ particles can also be uniquely\footnote{The
  equivalence of particle configurations under permutation is used
  here.} represented as an \emph{empirical measure}, that is a sum of
Dirac measures:
\begin{equation}\label{eq:EmpiricalMeasure} 
Z=(z_1,\dots,z_N) \ \mapsto \ \mu^N_Z
= \frac{1}{N}\sum_{j=1}^N \delta_{z_j}.
\end{equation} 
The resulting measure is normalized so as to give a probability
measure, which is obviously independent of any permutation of the
indices. The set of such empirical measures is denoted by
$\PP_N(E)$. Probability measures on $E$ are denoted by $\PP(E)$, so
that $\PP_N(E) \subset \PP(E)$.

When the number of particles go to infinity, we may have $\mu^N_Z
\rightarrow f\in \PP(E)$, where now $f$ is a distribution of particles
in $E$.  We call a ``limiting equation for the $N$-particle systems''
the (usually nonlinear) equation of the form~\eqref{eq:0020a} that is
satisfied by the probability distribution $f_t$ obtained as the limit
of $\mu^N_Z$, and we write its solution in the form of a nonlinear
semigroup $S^{N\!L}_t$: 
  \begin{equation*}
    f_t = S^{N\!L}_t (f_{\mbox{\tiny in}})
  \end{equation*}
is the solution of 
\begin{equation*}
    \partial_t f_t = Q(f_t)\,, \qquad f_0 = f_{\mbox{\tiny in}}\,.
\end{equation*}

The main result of this paper can be seen as a perturbation result:
given a solution to the limiting equation, we consider a sequence of
measures $\mu_{\ZZ_{\mbox{\tiny in}}}^N\in \PP_N(E)$ such that
\begin{eqnarray*}
  \mu_{\ZZ_{\mbox{\tiny in}}}^N &\rightarrow& f_{\mbox{\tiny in}}
\end{eqnarray*}
and  prove that for all $t$ in some interval $0\le t\le T$,
\begin{eqnarray*}
  \mu_{\ZZ_t}^N &\rightarrow& S^{N\!L}_t(f_{\mbox{\tiny in}}) \,.
\end{eqnarray*}
The convergence is established in the weak topology for the law of the
random empirical measures, as will be explained next.

\subsection{ $N$-particle dynamics of random measures and weak
  solutions of the limiting equation}
\label{sec:randommeasures}
A random point $\ZZ \in E^N$ with law $f^N \in \PPS(E^N)$ can be
identified with a random measure, denoted by $\mu^N_\ZZ \in \PP_N(E)$,
whose law is induced from $f^N$. We denote this law $\pi^N_P f^N \in
\PP(\PP(E))$. Note that since $E$ is a separable metric space, then so
is $\PP(E)$ by Prokhorov's Theorem, and we may define the space
$\PP(\PP(E))$ of probability measures on $\PP(E)$, as well as the set
of continuous bounded functions, denoted $C_b(\PP(E))$, which is the
dual of $\PP(\PP(E))$. In Section~\ref{sec:topology} we will discuss
how the choice of topology on $\PP(E)$ influences $C_b(\PP(E))$.

The nonlinear dynamics given by the semigroup $S^{N\!L}_t$ is
deterministic and defines a semigroup of operators on $C_b(\PP(E))$,
in a way that is reminiscent to the
equations~(\ref{eq:0030})-\eqref{eq:0031} but now using the duality
structure between $C_b(\PP(E))$ and $\PP(\PP(E))$. We define, for any
$f_{\mbox{\tiny in}}\in\PP(E)$ and $\Phi\in C_b(\PP(E))$,
   \begin{equation*}
     T^{\infty}_t \Phi( f_{\mbox{\tiny in}} ) = \Phi\left( S^{N\!L}_t(f_{\mbox{\tiny in}})\right).
  \end{equation*}
  The semigroup $T^{\infty}_t $ is called the {\em pullback semigroup}
  of $S^{N\!L}_t$. As we shall see, it plays a similar role for the
  deterministic limiting flow $S^{N\!L}_t$ associated with the
  equation~\eqref{eq:0020a} (lower half of the diagram), as the one
  the semigroup $T^{N}_t$ plays for the flow $\ZZ^N_t$ at the level of
  the $N$-particle system (upper half of the diagram): in both cases
  these are the \emph{dual statistical flows}.  Note that for making
  sense of the pullback semigroup $T^\infty_t$ on $C_b(\PP(E))$, one
  needs the map $f_{\mbox{\tiny in}} \mapsto S^{N\!L}_t(f_{\mbox{\tiny in}})$
  to be continuous, and this will be a major issue when all arguments
  are made precise.

  To connect the $N$-particle dynamics with the limiting dynamics, we
  also need mappings between $C_b(E^N)$ and
  $C_b(\PP(E))$. On the one hand 
  \begin{equation*}
    \pi^N: \left\{ 
      \begin{array}{cl} \displaystyle
        C_b(\PP(E)) & \rightarrow \ C_b(E^N) \vspace{0.2cm} \\
        \Phi & \mapsto \ \phi
      \end{array}
    \right.
  \end{equation*}
  is dual to the map $\PP(E)\ni f^N \mapsto \pi^N_P f^N
  \in\PP(\PP(E))$, and is defined by
  \begin{eqnarray*}
    \forall \, Z \in E^N,\quad \phi(Z) =  \left(\pi^N\Phi\right)(Z) = \Phi\left(\mu^N_Z\right),
  \end{eqnarray*}
  where the empirical measure $\mu^N_Z$ is defined through
  \eqref{eq:EmpiricalMeasure}.
  On the other hand 
\begin{equation}
  R^N: \left\{ 
    \begin{array}{cl}
      C_b(E^N) & \rightarrow \ C_b(\PP(E)) \vspace{0.2cm} \\
      \phi & \mapsto \ \Phi
    \end{array}
    \right.
  \end{equation}
  is defined in the following way: for each $\phi\in C_b(E^N)$,
  $R^N[\phi]$ is evaluated at the point $f \in\PP(E)$ as
  \begin{equation}
\label{eq:0036}
f \mapsto R^N[\phi](f)=
R^N_{\phi}(f) = \int_{E^N} \phi(z_1,\dots,z_N) \, f({\rm
  d}z_1) \, \cdots \, f({\rm d}z_N),
  \end{equation}
  which can be interpreted, as we shall see later, as a real valued
  polynomial taking probability measures as
  arguments. 

\subsection{The abstract theorem}
\label{sec:abstracttheorem}

We are now ready to give a more precise version of the main abstract
theorem. The exact statement involves rather technical definitions,
and its validity depends on five assumptions, {\bf (A1)} to {\bf (A5)}
which are properly stated in Section~\ref{sec:lemmas}. The first
assumption is the requirement of symmetry under permutations that has
already been stated. The remaining conditons are (1) estimates on the
regularity and stability of the nonlinear semigroup, and (2)
consistency estimates that quantifies that the $N$-particle systems
and the limiting semigroup are compatible. 

\begin{theo}[Fluctuation estimate]
  \label{theo:abstract}
  Consider a process $(\ZZ^N_t)_{t\ge 0}$ in $E^N/ \mathfrak{S}_N$, and
  the related semigroups $S^N_t$ and $T^N_t$ as defined above.
  Let $f_{\mbox{{\tiny \emph{in}}}}\in \PP(E)$,
  and consider a hierarchy of $N$-particle solutions $f^N_t = S^N _t
  (f_{\mbox{{\tiny \emph{in}}}}^{\otimes N})$, and a solution $f_t=S^{N\!L} _t (f_{\mbox{{\tiny \emph{in}}}})$ to the
  limit equation. We assume that {\bf (A1)} to {\bf (A5)} hold. 

  Then there is an absolute constant $C>0$ and, for any $T \in
  (0,\infty)$, there are  constants $C_{T}, \tilde C_{T} >0$
  (depending on $T$) such that for any $N, \ell \in \N^*$, with $N \ge
  2 \ell$, and for any
  \begin{equation*}
    \varphi = \varphi_1 \otimes \dots \otimes \,
    \varphi_\ell \in \FF^{\otimes \ell},
    \quad 
    \varphi_j \in \FF, \ \|\varphi_j \|_\FF \le 1,
  \end{equation*}
  we have
  \begin{eqnarray}
  \label{eq:cvgabstract1}
  &&\quad \sup_{[0,T)}\left| \left \langle \left( S^N_t(f_{\mbox{{\tiny \emph{in}}}}^{\otimes N}) 
        - \left(S^{N\! L}_t(f_{\mbox{{\tiny \emph{in}}}}) \right)^{\otimes N} \right), \varphi \otimes {\bf 1}^{N-\ell}
    \right\rangle \right| 
  \\ \nonumber 
  &&\qquad \qquad \qquad 
  \le  C \,  \frac{   \ell^2}{N} + C_{T} \,  \ell^2 \, 
  \varepsilon(N)  + 
  \tilde C_{T} \, \ell \, \Omega_N ^{\GG_3} (f_{\mbox{{\tiny \emph{in}}}}) , 
  \end{eqnarray}
   with
\begin{equation}\label{def:OmegaG3N}
  \Omega_N ^{\GG_3} (f_{\mbox{{\tiny \emph{in}}}}) : = \int_{E^N} 
  \mbox{{\em dist}}_{\GG_3}\left(\mu^N_Z,f_{\mbox{{\tiny \emph{in}}}}\right) \, f^{\otimes N}_{\mbox{{\tiny \emph{in}}}} ({\rm d}Z).
\end{equation}
   The space $\FF\subset C_b(E)$ 
   and the 
   distance $\mbox{{\em dist}}_{\GG_3}$ are defined later in
   Section~\ref{sec:lemmas}. 
\end{theo}

We have used the notation $\varphi = \varphi_1 \otimes \cdots \otimes
\, \varphi_\ell$ to denote
\begin{equation*}
\varphi(z_1,\dots,z_\ell) = \varphi_1(z_1)
\varphi_2(z_2) \cdots \varphi_\ell(z_\ell),
\end{equation*}
and $\varphi \otimes {\bf 1}^{N-\ell}$ to  denote
\begin{equation*}
(\varphi \otimes {\bf
  1}^{N-\ell})(z_1,\dots,z_N) = \varphi_1(z_1) \varphi_2(z_2) \cdots
\varphi_\ell(z_\ell).
\end{equation*}

We first note that the left hand side of (\ref{eq:cvgabstract1}) is
the same as the left hand side of (\ref{eq:0015}), the only difference
being that the inital data to the $N$-particle system are assumed to
factorize: $f^N_{\mbox{\tiny in}}= f_{\mbox{\tiny in}}^{\otimes N}$, where
$f_{\mbox{\tiny in}}$ is also the inital data to the limiting nonlinear
equation.  This is a stronger hypothesis than merely requiring the
initial data to be chaotic (where the initial data to the $N$-particle
system may factorize only in the limit of infinitely many
particles). This restriction simplifies the proof, but it can easily
be relaxed at the cost of some additional error terms.

The restriction to test functions of the form $\varphi_1 \otimes
\dots \otimes \, \varphi_\ell\otimes 1^{\otimes (N-\ell)}$, i.e.
functions depending only on the first $\ell$ variables, corresponds to
analysing $\ell$-particle marginals. Hence the theorem implies the
propagation of chaos as soon as $\FF$ is dense in $C_b(E)$ in the
topology of uniform convergence on compact sets. This condition is
satisfied in all examples given below.

\subsection{Main steps of the proof}
\label{sec:sketch-proof}

The proof begins by splitting the  quantity we want to estimate, 
\begin{eqnarray*}
  \left \langle \left( S^N_t(f_{\mbox{{\tiny in}}}^{\otimes N}) 
        - \left(S^{N\! L}_t(f_{\mbox{{\tiny in}}}) \right)^{\otimes N} \right), \varphi \otimes {\bf 1}^{N-\ell}
    \right\rangle\,,
\end{eqnarray*}
in three parts, each one corresponding to one of the error terms in
the right hand side of the equation~(\ref{eq:cvgabstract1}):
\begin{eqnarray}
  &&\left| \left \langle S^N_t(f_{\mbox{{\tiny in}}}^{\otimes N}) - \left( S^{N\!L}_t
        (f_{\mbox{{\tiny in}}})  \right)^{\otimes N} , 
      \varphi
      \otimes 1^{\otimes N-\ell} \right\rangle \right| \le
  \nonumber   \\
  &&\le \left| \left\langle S^N_t(f_{\mbox{{\tiny in}}}^{\otimes N}), 
      \varphi  \otimes
      1^{\otimes N-\ell} \right\rangle -
    \left \langle S^N_t(f_{\mbox{{\tiny in}}}^{\otimes N}), 
      \nonumber     R^\ell[\varphi] \circ \mu^N_Z \right\rangle  \right|
  \\
  &&+ \left| \left\langle f_{\mbox{{\tiny in}}}^{\otimes N}, T^N_t ( R^\ell[\varphi] \circ
      \mu^N_Z) \right\rangle
    - \left\langle f_{\mbox{{\tiny in}}}^{\otimes N},       (T_t ^\infty R^\ell[\varphi] ) \circ \mu^N_Z) \right\rangle  \right|
  \nonumber
  \\  
  &&+ \left| \left\langle f_{\mbox{{\tiny in}}}^{\otimes N}, (T_t ^\infty R^\ell[\varphi] ) \circ
      \mu^N_Z) \right\rangle - \left\langle (S^{N\!L}_t (f_{\mbox{{\tiny in}}}))^{\otimes \ell} ,
      \varphi \right\rangle \right| =: \TT_1 + \TT_2 + \TT_ 3.
\label{eq:cvtabst10}
\end{eqnarray} 


Then each of these terms is estimated separately: 

\begin{enumerate} 
\item The first term, $\TT_1$, is bounded by $C\ell^2/N$, as proven in
  Lemma~\ref{lem:symmetrization}. From the definitions of $R^{\ell}$
  and $\mu^N_Z$, it follows that $R^{\ell}[\varphi](\mu^N_Z)$ is a sum
  of terms of the form $\varphi(z_{j_1},\dots,z_{j_{\ell}})$, where
  each index $j_1,\dots,j_{\ell}$ describes the set
  $\{1,\dots,N\}$. The error is then due to the fraction of terms for
  which two or more of the indices $j_1,\dots,j_{\ell}$ are the
  same. Hence the estimate of $\TT_1$ is of purely combinatorial
  nature, and only depends on the symmetry under permutation,
  i.e. the assumption {\bf (A1)}.

\item The estimate of the second term, $\TT_2$, relies on the
  convergence of the $N$-particle semigroups $T^N_t$ to the limiting
  semigroup $T^{\infty}_t$. This is where the $N$-particle dynamics
  and the limiting dynamics are compared. The estimate can be found in
  Lemma~\ref{lem:steptwo}, which depends on the consistency assumption
  {\bf (A3)} on the generators and the stability assumption {\bf (A4)}
  on the limiting dynamics. 
  While the generator $G^N$ of $T^N_t$ can be defined in a
  straightforward manner, defining and estimating the generator
  $G^{\infty}$ of $T^{\infty}_t$ requires a more detailed analysis. It
  indeed involves \emph{derivatives} of functions acting on $\PP(E)$,
  and the required \emph{differential structure} depends on the
  topology and metric structure chosen on $\PP(E)$. The generator
  $G^{\infty}$ is characterized in Lemma~\ref{lem:nonlgenerator}. The
  assumption {\bf (A4)} is proved by establishing refined stability
  estimates on the limiting semigroup, showing their
  \emph{differentiability} according to the initial data in a metric
  compatible with the previous steps.

\item For the last term $\TT_3$, we note that
\begin{eqnarray*}
  T^{\infty}_t
  R^{\ell}[\varphi](\mu^N_Z)&=&R^{\ell}[\varphi]\left(S^{N\!L}_t(\mu^N_Z)\right)\,=\,
  \left\langle (S^{N\!L}_t(\mu^N_Z))^{\otimes\ell},\varphi\right\rangle,
\end{eqnarray*}
which means that the nonlinear limiting equation is solved taking a
sum of Dirac masses as initial data, and an $\ell$-fold product of the
solution is integrated against $\phi$. The resulting function of
$Z=(z_1,\dots,z_N)$ is then integrated against
$f_{\mbox{\tiny in}}^{\otimes N}$, which amounts to taking an average over
all intial data such that the position of the $N$ particles are
independently taken at random from the law $f_{\mathit in}$. When $N$
is large, the random empirical measures $\mu^N_\ZZ$ are close to
$f_{\mbox{\tiny in}}$, i.e. $\Pi^N _P(f_{\mbox{{\tiny in}}} ^{\otimes N}) \rightharpoonup
\delta_{f_{\mbox{{\tiny in}}}}$ in $\PP(\PP(E))$. This implies that the term $\TT_3$
vanish in the limit when $N$ goes to infinity. However the rate of
this convergence sensitively depends on the regularity of the test
function $\varphi$ and on the continuity properties of
$S^{N\!L}_t$. In all cases considered here, the error $\TT_3$
dominates the other error terms, and effectively determines the rate
of convergence in the propagation of chaos. The precise result,
together with the required assumptions, is given in
Lemma~\ref{lem:stepthree}.
\end{enumerate}

\section{Metrics on $\PP(E)$ and differentiability of functions
  on $\PP(E)$}
\label{sec:topology}
\setcounter{equation}{0}
\setcounter{theo}{0}


This section contains the techical details concerning the space of
probability measures on $\PP(E)$ and its dual, that is the space of
continuous functions acting on $\PP(E)$.

\subsection{The metric issue}
\label{sec:funct-set}

$\PP(E)$ is our fundamental ``state space'', where we compare
the marginals of the $N$-particle density $f^N_t$ and the chaotic
infinite-particle dynamics $f_t $ through their observables, i.e. the
evolution of continuous bounded functions on $E^N$ and $\PP(E)$
respectively under the dual dynamics $T^N _t$ and $T^\infty_t$.

There are two canonical choices of topology on the space of
probabilities, which determine two different sets $C_b(\PP(E))$.

On the one hand, for a given locally compact and separable metric
space $\EE$, the space $M^1(\EE)$ of finite Borel measures on $\EE$ is
a Banach space when endowed with the total variation norm:
\begin{eqnarray*}
  \forall \, f \in M^1(\EE), \quad \| f \|_{TV} &:=& f^+(\EE) + f^-(\EE) \\
  &=& \sup_{\phi \in C_b(Z), \, \| \phi \|_\infty \le 1} \langle
  f,\phi \rangle = \sup_{\phi \in C_0(\EE), \, \| \phi \|_\infty \le
    1} \langle f,\phi \rangle, 
\end{eqnarray*} 
where $f = f^+ - f^-$ stands for the Hahn decomposition and the
equality between the two last terms comes from the fact that $\EE$ is
locally compact and separable. 

We recall that $f_k \xrightarrow[]{TV} f$ (strong topology) when
$(f_k)$ and $f$ belongs to $M^1(\EE)$ and $\| f_k - f \|_{TV} \to 0$
when $k\to\infty$, and that $f_k \wto f$ (weak topology), if
$$
\forall \, \varphi \in C_b(Z) \qquad \langle f, \varphi \rangle =
\lim_{k\to\infty} \langle f_k , \varphi \rangle\,.
$$
The associated topology is denoted by
$\sigma(M^1(\EE),C_b(\EE))$. However, the weak convergence can be
associated with different, non-equivalent metrics, and the choice of
metric plays an important role as soon as one wants to perform {\em
  differential calculus} on $\PP(\EE)$.

In the sequel, we will denote by $C_b(\PP(E),w)$ the space of
continuous and bounded functions on $\PP(E)$ {\em endowed with the
  weak topology}, and $C_b(\PP(E),TV)$ the space of continuous and
bounded functions on $\PP(E)$ {\em endowed with the total variation
  norm}.  It is clear that $C_b(\PP(E),w) \subset C_b(\PP(E),TV)$
since $f_k \xrightarrow[]{TV} f$ implies $f_k \wto f$.

However, the supremum norm $\| \Phi \|_{L^\infty(\PP(E))}$ {\em does
  not} depend on the choice of topology on $\PP(E)$, and endows the
two previous sets with a Banach space topology. The transformations
$\pi^N$ and $R^N$ satisfy:
\begin{equation}\label{eq:compat:infty}
  \left\| \pi^N \Phi \right\|_{L^\infty(E^N)} \le \| \Phi
  \|_{L^\infty(\PP(E))} \ \mbox{ and } \ \| R^N[\phi] \|_{L^\infty(\PP(E))} \le
  \| \phi \|_{L^\infty(E^N)}.
\end{equation}

The transformation $\pi^N$ is well defined from $C_b(\PP(E),w)$ to
$C_b(E^N)$, but  it does not map $C_b(\PP(E),TV)$ into
$C_b(E^N)$. 

In the other way round, the transformation $R^N$ is well defined from
$C_b(E^N)$ to $C_b(\PP(E),w)$, and therefore also from $C_b(E^N)$ to
$C_b(\PP(E),TV)$: for any $\phi \in C_b(E^N)$ and for any sequence
$f_k$ so that the weak convergence $f_k \wto f$ holds, we have
$f_k^{\otimes N} \wto f^{\otimes N}$, and then $R^N[\phi](f_k)
\to R^N[\phi](f)$.

\smallskip 

The different metric structures associated with the weak topology are
not seen at the level of $C_b(\PP(E),w)$. However any norm (or
semi-norm) ``more regular'' than the uniform norm on $C_b(\PP(E))$ (in
the sense of controlling some modulus of continuity or some
differential) strongly depends on this choice, as is illustrated by
the abstract Lipschitz spaces defined below.

\begin{defin}
\label{defGG} 
Let $m_\GG : E \to \R_+$ be given. Then we define the following
weighted subspace of probability measures
  $$
  \PP_{\GG}(E)  := \{ f \in \PP(E); \,\, \langle f, m_\GG \rangle < \infty \}, 
  $$
  together with a  corresponding space of ``increments'',
  $$
  \mathcal{I} \PP_\GG(E) := \left\{f_1 - f_2 \ ; \ f_1, \, f_2 \in \PP_{\GG}(E)  \right\}.
  $$
  If moreover there is a vector space $\GG$ (with norm denoted by $\|
  \cdot \|_\GG$) which contains $\mathcal{I} \PP_\GG(E)$, we then define
  the following distance on $\PP_{\GG}(E)$
  $$
  \forall \, f_1, \, f_2 \in \PP_{\GG}(E), \quad
  \mbox{\emph{dist}}_\GG(f_1,f_2) := \| f_1 - f_2 \|_\GG.
  $$
 \end{defin}

 \begin{rem}
   Note carefully that the space of increments $\mathcal{I} \PP_\GG(E)$ is
   \emph{not} a vector space in general. 
 \end{rem}

Now we can define a precised notion of equivalence of metrics: 
\begin{defin}
  We say that $\PP_\GG(E)$ has a \emph{bounded diameter} if there exits
  $K_\GG >0$ such that
  \[ 
  \forall \, f \in \PP_\GG(E), \quad \mbox{\emph{dist}}_\GG(f,g) \le K_\GG
  \] 
  for some given fixed $g \in \PP_\GG(E)$. 

  Two metrics $d_0$ and $d_1$ on $\PP_\GG(E)$ are said to be \emph{H\"older
    uniformly equivalent on bounded sets} if there exists $\kappa \in
  (0,\infty)$ and for any $a \in (0,\infty)$ there exists $C_a \in
  (0,\infty)$ such that
  $$
  \forall \, f_1, \, f_2 \in \BB P_{\GG,a}, \quad
  \frac{1}{C_a} \, [d_2(f_1,f_2)]^\kappa \le d_1(f_1, f_2) \le C_a \, [d_2(f_1,f_2)]^\kappa
  $$
  where 
  \[
  \BB P_{\GG,a} := \left\{ f \in \PP_\GG(E) \ ; \ \langle f, m_\GG
    \rangle \le a \right\}.
  \]

  Finally, we say that two normed spaces $\GG_0$ and $\GG_1$ are
  \emph{H\"older uniformly equivalent (on bounded sets)} if this
  is the case for the corresponding metrics. 
\end{defin}

We also define the vector space $UC( \PP_\GG(E);\R)$ of
  uniformly continuous and bounded function $\Psi : \PP_\GG(E) \to \R$,
  where the continuity is related the metric topology on $\PP_\GG(E)$
  defined by $\mbox{dist}_{\GG}$ above. Observe that this is a Banach
  space when endowed with the supremum norm.

\begin{ex}\label{expleTV} 
  With the choice $m_\GG := 1$, $\| \cdot \|_\GG := \|\cdot \|_{TV}$
  we obtain $\PP_\GG(E)(E) = (\PP(E),TV)$ endowed with the total variation
  norm.
  \end{ex}



 \subsection{Examples of distances on measures when $E = \R^d$}
 \label{subsec:ExpleMetrics} 

 There are many ways to define distances on $\PP(E)$ which are
 topologically equivalent to the weak topology of measures, see for
 instance \cite{BookRachev,coursCT}.

 We list below some well-known distances on $\PP(\R^d)$ or on its
 subsets
 $$
\PP_q(\R^d) := \{ f \in \PP(\R^d); \,\, M_q(f)  < \infty \}, 
\quad q \ge 0, 
$$
where the moment $M_q(f)$ of order $q$ of a probability measure is defined as
$$
M_q(f) :=  \left\langle f \,,\, \langle v \rangle^q \right\rangle, \quad  \langle v \rangle^2 = 1 + |v|^2.
$$
These distances are all H\"older uniformly equivalent to the weak
topology $\sigma (\PP(E),C_b(E))$ on the bounded subsets
\[
\BB P_{q,a}(E) := \left\{f \in \PP_q(\R^d), \,\, M_q(f)  \le a \right\}
\]
for any $a \in (0,\infty)$ and for $q$ large enough. For more
informations we refer to~\cite{Dudley2002}.

\begin{ex}[Dual-H\"older, or Zolotarev's, distances]\label{expleZolotarev}
  Denote by $\mbox{{\em dist}}_E$ a distance on $E$ and fix $z_0 \in
  E$ (e.g. $z_0=0$ when $E = \R^d$ in the sequel). Denote by
  $\mbox{{\em Lip}}_0(E)$ the set of Lipschitz functions on $E$
  vanishing at one arbitrary point $z_0 \in E$ endowed with the norm
  $$
  [\varphi ]_{\mbox{{\scriptsize {\em Lip}}}} = [\varphi ]_1 :=
  \sup_{z,\tilde z \in E, \ z \not = \tilde z} {|\varphi(z) -
    \varphi(\tilde z)| \over \hbox{{\em
        dist}}_E(z,\tilde z)} .
  $$
  We then define the dual norm: take $m_\GG := 1$ and endow $\PP_\GG(E)$
  with
 \begin{equation}\label{def;[]*s} 
   \forall \, f,g \in \PP_\GG(E),
  \quad [g-f]^*_1 := \sup_{\varphi \in Lip_0(E)} {\langle
    g - f, \varphi \rangle \over [\varphi]_1 }. 
  \end{equation}
\end{ex}

\begin{ex}[Monge-Kantorovich-Wasserstein distances]\label{expleWp} 
  For $q \in [1,\infty)$, define
  $$
  \PP_\GG(E) (E) = \PP_q(E):= \left\{ f \in \PP(E); \,\, \langle f, m_\GG \rangle :=\left \langle f,
    \mbox{{\em dist}}(\cdot,v_0)^q \right \rangle < \infty \right\}
  $$
  and the Monge-Kantorovich-Wasserstein (MKW) distance $W_q$ by 
  \begin{equation}\label{eq:wasserstein}
    \forall \, f,g \in \PP_q(E), \quad  W^q_q(f,g) 
    :=  \inf_{p \in \Pi(f,g)}  \int_{E\times E}  \mbox{{\em
        dist}}_E(z,\tilde z)^q \, p({\rm d}z,{\rm d}\tilde z),
\end{equation}
where $\Pi(f,g)$ denote the set of probability measures $p \in
\PP(E \times E)$ with marginals $f$ and $g$ ($p(A\times E) = f
(A)$, $p(E\times A) = g (A)$ for any Borel set $A \subset E$).  Note
that for $Z,\tilde Z \in E^N$ and any $q \in [1,\infty)$, one has
\begin{equation}\label{Wqellq} 
  W_q\left(\mu^N_Z,\mu^N_{\tilde Z}\right) =
  d_{\ell^q(E^N/\SN)}
  (Z, \tilde Z) := \min_{\sigma \in \SN} \left( {1 \over N} \sum_{i=1}^N
    \mbox{{\em dist}}_E(z_i,\tilde z_{\sigma(i)})^q \right)^{1/q}, 
\end{equation} 
and that
\begin{equation}\label{W1KR} 
  \forall \, f, \, g \in
  P_1(E) , \quad W_1 (f,g) 
  = [f-g]^*_1 = \sup_{\phi \in \mbox{\tiny {\em Lip}}_0(E)}\, \left \langle f-g,
  \phi \right \rangle
\end{equation} 
as well as 
\begin{equation}\label{estim:W1Wq}
 \forall \, q  \in [1,\infty), \,\,\, 
  \forall \, f, \, g \in \PP_q(R^d),  \quad 
  W_1(f,g) \le W_q(f,g) .
\end{equation}
We refer to \cite{VillaniTOT} and the references therein for more
details on the Monge-Kantorovich-Wasserstein distances and for a proof
of these claims.
\end{ex}

\begin{ex}[Fourier-based norms] \label{expleFourier} For
  $E=\R^d$, $m_{\GG} := 1$, let
  \[ 
  \forall \, f \in \mathcal{T} \PP_\GG(E), \quad \| f \|_{\GG} = |f|_s
  := \sup_{\xi \in \R^d} \frac{|\hat f(\xi)|}{ \langle \xi \rangle^s}, \quad s > 0.
  \]
  We denote by $\HH^{-s}$ (which includes $\II \PP_\GG(E)$ for $s$
  large enough) the Banach space associated to the norm $| \cdot
  |_{s}$. Such norms first appeared in connection with kinetic theory
  in~\cite{GabettaTW95}. 
\end{ex}

\begin{ex}[Negative Sobolev norms]\label{expleH-s} 
  For $E=\R^d$, $m_{\GG} := 1$, let
  \[ 
  \forall \, f \in \mathcal{T} \PP_\GG(E), \quad \| f \|_{\GG} = \| f
  \|_{ H^{-s} (\R^d)} := \left\| \frac{\hat f(\xi)}{\langle \xi
      \rangle^s} \right\|_{L^2(\R^d)}, \quad s >0.
  \]
  We denote by $H^{-s}$ (which includes $\II \PP_\GG(E)$ for $s$ large
  enough) the Hilbert space associated to the norm $\| \cdot
  \|_{H^{-s}}$.

  For $E=\R^d$, $m_{\GG} := 1$, and some integers $k,\ell \ge 0$, we
  also define
 \[ 
  \forall \, f \in \mathcal{T} \PP_\GG(E), \quad \| f \|_{\GG}
  = \| f \|_{  H^{-k}_{-\ell} (\R^d)} := \sup_{\varphi \in H^k_\ell} \langle f,\varphi \rangle, 
  \]
  where 
  \[ 
   \left\| \varphi \right\|_{H^k_{\ell}}^2 := \sum_{  |\alpha| \le k} \int_{\R^d} |\partial^\alpha \varphi (z)|^2 \, 
   \langle z \rangle^{2\ell} \, {\rm d}z.
  \]
  We denote by $H^{-k}_{-\ell}(\R^d)$ (which includes $\II \PP_\GG(E)$
  for $k$ large enough) the Hilbert space associated to the norm $\|
  \cdot \|_{H^{-k}_{-\ell}(\R^d)}$.
  \end{ex}


\subsection{Differential calculus for functions of probability measures }
\label{sec:diff-measures}

We start with a definition of Lipschitz regularity, for which a mere
metric structure is sufficient.

\begin{defin}\label{def:Holdercalculus}
  For metric spaces $\tilde \GG_1$ and $\tilde\GG_2$ we denote by
  $C^{0,1}(\tilde \GG_1,\tilde\GG_2)$ the space of functions from
  $\tilde \GG_1$ to $\tilde\GG_2$ with Lipschitz regularity, i.e. the
  set of functions $\Psi : \tilde \GG_1 \to \tilde\GG_2$
  such that there exists a constant $C >0$ so that
  \begin{equation}
    \label{eq:devphi}
    \forall \, f,g \in \tilde \GG_1, \quad 
    \mbox{{\em dist}}_{\tilde\GG_2} \left(\Psi(g) , \Psi(f) \right)
    \le C \, \mbox{{\em dist}}_{\tilde\GG_1} (g,f).
  \end{equation}
  We then define the semi-norm $[ \cdot ]_ {C^{0,1}(\tilde
    \GG_1,\tilde\GG_2)}$ on $C^{0,1}(\tilde \GG_1,\tilde\GG_2)$ as the
  infimum of the constants $C > 0$ such that \eqref{eq:devphi} holds.
 \end{defin}

 The next step consists in defining a higher order differential
 calculus; this is where the assumption that metrics are inherited
 from a normed vector space structure plays a role. 

\begin{defin}\label{def:diffcalculus}
  Let $\GG_1$ and $\GG_2$ be normed spaces, and let $\tilde \GG_1$ and
  $\tilde \GG_2$ be two metric spaces such that $\tilde \GG_i - \tilde
  \GG_i \subset \GG_i$. For $k \in \N$, we define $C^{k,1}(\tilde
  \GG_1; \tilde \GG_2)$ to be the set of bounded continuous functions
  $\Psi : \tilde \GG_1 \to \tilde \GG_2$ such that there exists $D^j
  \Psi : \tilde \GG_1 \to \BB^j(\GG_1,\GG_2)$ continuous and bounded,
  where $\BB^j(\GG_1,\GG_2)$ is the space of bounded $j$-multilinear
  applications from $\GG_1$ to $\GG_2$ (endowed with its canonical
  norm) for $j = 1, \dots, k$, and some constants $C_j >0$, $j=0,
  \dots , k$, so that for any $j=0, \dots, k$
  \begin{equation}
    \label{eq:devdist}
    \forall \, f,g \in \tilde \GG_1, \quad
    \left\| \Psi(g) -  
    \ \sum_{i=0}^j \left \langle D^i \Psi(f) ,  
    (g-f)^{\otimes i} \right \rangle \right\|_{\GG_2}  
  \le C_j \, \| g-f \|_{\GG_1}^{j+1}
\end{equation}
(with the convention $D^0 \Psi = \Psi$).  


We also define the following seminorms on $C^{k,1}(\tilde \GG_1,\tilde
\GG_2)$
\[ [\Psi]_{j,0} := \sup_{f \in \tilde \GG_1} \left\| D^j \Psi(f)
\right\|_{\BB^j(\GG_1,\GG_2)}, \quad j=1, \dots, k,
\]
with 
\[
\left\| L \right\|_{\BB^j(\GG_1,\GG_2)} := \sup_{h_i, \, \| h_i \|_{\GG_1} \le 1, \, 1 \le i \le j}
\left\|L\left(h_1,\dots, h_j\right) \right\|_{\GG_2},
\]
and 
\[
[\Psi]_{j,1} :=   \sup_{f,g \in \tilde\GG_1} \frac{ \Big\| \Psi(g) -
  \sum_{i=0}^j \langle D^i \Psi(f) , (g-f)^{\otimes i}
  \rangle \Big\|_{\GG_2}}{\| g-f \|_{\GG_1}^{j+1} }.
\]

Finally we combine these semi-norms into the norm 
\[
\| \Psi \|_{C^{k,1}(\tilde\GG_1,\tilde\GG_2)} = 
\sum_{j=1} ^k  \, [\Psi ]_{j,0} + [\Psi]_{k,1}. 
\]
\end{defin}
%


\begin{rem}
  Observe that for any $j \ge 1$ the LHS of \eqref{eq:devdist} makes
  sense since 
  \begin{multline*}
    \Psi(g) - \sum_{i=0}^j \left \langle D^i \Psi(f) ,  
      (g-f)^{\otimes i} \right \rangle \\ = \Big[\Psi(g) - \Psi(f) \Big]
    - \sum_{i=1} ^j \left \langle D^i \Psi(f), (g-f)^{\otimes i}
    \right \rangle \in \GG_2
  \end{multline*}
  since $\Psi(g) - \Psi(f) \in \tilde \GG_2 - \tilde \GG_2 \subset
  \GG_2$ and $D^i \Psi(f) \in \BB^j(\GG_1,\GG_2)$. Then note that our
  definition is very close to the usual Fr\'echet definition of
  differentiability in Banach spaces for the function $h \mapsto
  \Psi(f+h)$ with $h = g-f \in \GG_1$, except that the domain and
  range are restricted to subsets that have no vectorial structures
  and are not open within $\GG_1$ and $\GG_2$. We also only consider
  Lipschitz differentiability.
\end{rem}

The following lemma confirms that this differential calculus is
well-behaved for composition, which seems to be a minimal requirement
for further applications.
\begin{lem}\label{lem:DL} 
  Consider $\UU \in C^{k,1} (\tilde \GG_1, \tilde \GG_2)$ and $\VV \in
  C^{k,1} (\tilde \GG_2, \tilde \GG_3)$. Then the composition $\Psi :=
  \VV \circ \UU$ belongs to $C^{k,1} (\tilde \GG_1, \tilde
  \GG_3)$. Moreover the following chain rule holds at first order
  $k=1$
\begin{equation}\label{ChainRule1} 
\forall \, f \in \tilde \GG_1, \quad  D \Psi[f] = D \VV   [\UU(f)] \circ D \UU[f], 
\end{equation}
with the estimates
\begin{equation*}
\left\{ 
\begin{array}{l} \displaystyle
  [\Psi ]_{0,1} \le [\VV ]_{0,1} \, [\UU]_{0,1}, \vspace{0.2cm} \\ \displaystyle
  [\Psi ]_{1,0} \le [\VV ]_{1,0} \, [\UU]_{1,0}, \vspace{0.2cm} \\ \displaystyle
 [\Psi ]_{1,1} \le [ \VV ]_{1,0} \, [ \UU ]_{1,1} + [ \VV ]_{1,1} \, [
 \UU ]_{0,1}^2.
\end{array}
\right.
\end{equation*}

At second order $k=2$ one also has the chain rule
\begin{equation}\label{ChainRule2} 
\forall \, f \in \tilde \GG_1, \quad  D^2 \Psi[f] = D^2 \VV [\UU(f)]
\circ (D \UU[f] \otimes D \UU[f]) + D \VV [\UU(f)] \circ D^2
\UU[f].  
\end{equation}
\end{lem}

\noindent{\sl Proof of Lemma~\ref{lem:DL}. } It is straightforward by
writing and compounding the expansions of $\UU$ and $\VV$ provided by
Definition~\ref{def:diffcalculus}. \qed


\subsection{The subalgebra of polynomials in  $C_b(\PP(E))$}
\label{sec:compatibility}

In Section~\ref{sec:abstractchaos} we defined a map $R^{\ell}:
C(E^{\ell})\rightarrow C_b(\PP(E))$, which may be used to define a
subalgebra of polyomials in $ C_b(\PP(E))$. We first define the
monomials:  

\begin{defin}\label{def:duality}
  A \emph{monomial} in $C_b(\PP(E))$ of degree $\ell$ is a function
  $R^\ell[\varphi]$ with $\varphi = \varphi_1 \otimes \cdots \otimes
  \varphi_\ell$, $\varphi_i \in C_b(E)$ and $\ell \in \N$.
Explicitly 
\begin{eqnarray*}
  R^{\ell}[\varphi](f) &=&  \int_{E^\ell} \varphi(z_1,\dots,z_\ell) \,
  {\rm d}f^{\otimes \ell}(z_1,\dots,z_\ell)\\
&=& \prod_{j=1}^{\ell} \int_E \varphi_j(z) \, {\rm d}f(z)\,,
\end{eqnarray*}
which is well-defined for all $f\in \PP(E)$.
\end{defin}
The product of two monomials is defined in a natural way  by 
$R^{\ell_1}[\phi]R^{\ell_2}[\psi] =
R^{\ell_1+\ell_2}[\phi\otimes\psi]$, and the polynomial functions are
 linear combinations of monomials. These form a  subalgebra of
$C_b(\PP(E))$  that
contains the constants and separates points in $\PP(E)$, and hence
the Stone-Weierstrass Theorem implies that this subalgebra is dense in
$C_b(\PP(E))$, where the meaning of ``dense'' depends on the topology
chosen on $\PP(E)$. 

While the polynomials in $\R$ always are differentiable, the smoothness
of the polynomials  depends on the metric structure.
 We need first some preliminary definitions.
\begin{defin}\label{def:DualitySpaces}

\begin{itemize}
\item \textbf{Duality of type 1}: We say that a pair $(\FF,\GG)$ of
  normed vector spaces such that $\FF \subset C_b(E)$ and
  $\PP(E)-\PP(E) \subset \GG$ \emph{satisfies a duality inequality}
  if
\begin{equation} \label{eq:dualiteFG} 
\forall \, f,g \in \PP(E),
  \,\, \forall \, \varphi \in \FF, \quad |\langle(f-g), \varphi \rangle
  |\le C \, \|f-g \|_\GG \, \| \varphi \|_\FF.
\end{equation} 

\item \textbf{Duality of type 2}: More generally we say that a pair
  $(\FF,\PP_\GG(E))$ of a normed vector space $\FF \subset C_b(E)$ endowed
  with the norm $\|\cdot \|_\FF$ and a probability space $\PP_\GG(E)
  \subset \PP(E)$ endowed with a metric $d_\GG$ \emph{satisfies a
    duality inequality} if
\begin{equation} \label{eq:dualiteFGbis} \forall \, f,g  \in \PP_\GG(E),
  \,\, \forall \, \varphi \in \FF, \quad |\langle g-f, \varphi \rangle
  |\le C \, \mbox{\emph{dist}}_\GG(f,g)\, \| \varphi \|_\FF.
\end{equation} 
\end{itemize}
\end{defin}

\begin{lem}\label{lem:polyCk}
  If $\varphi \in \FF^\ell$ and the pair $(\FF,\GG)$ satisfy a duality of type
  1, 
  the polynomial function $R^\ell[\varphi]$ is of class
  $C^{k,1}(\PP_\GG(E),\R)$ for any $k\ge0$. In the more general case where
  the pair $(\FF,\PP_\GG(E))$ satisfies a duality of type
  2, 
  the polynomial function $R^\ell[\varphi]$ is at least of class
  $C^{0,1}(\PP_\GG(E),\R)$.
\end{lem}

\begin{proof} It is clearly enough to prove the lemma for monomials,
  and the proof then mainly follows from the multilinearity of
  $R$. In the case of duality of type 2, then the conclusion follows from 
\begin{eqnarray*}
  R^\ell[\varphi](f_2) - R^\ell[\varphi](f_1) 
  &=& \sum_{i=1}^\ell \left(  \prod_{1 \le k < i } \langle \varphi_k,
    f_2 \rangle \right) 
  \langle \varphi_i, f_2 - f_1 \rangle 
    \left(  \prod_{i < k \le \ell } \langle \varphi_k, f_1 \rangle \right).
\end{eqnarray*}

In the case of duality of type 1, we define 
\[
\GG \to \R, \quad h \mapsto DR^\ell[\varphi] (f) (h) :=
\sum_{i=1}^\ell \left( \prod_{j \not= i } \langle \varphi_j, f \rangle \right)
\langle \varphi_i, h \rangle,
\]
and we write 
\begin{multline*}
  R^\ell[\varphi](f_2) - R^\ell[\varphi](f_1) - DR^\ell[\varphi] (f_1)
  (f_2 - f_1) =
  \\
  = \sum_{1\le j < i \le \ell} \left( \prod_{1 \le k < j } \langle
    \varphi_k, f_2 \rangle \right) \langle \varphi_j, f_2 - f_1
  \rangle \left( \prod_{j < k <i} \langle \varphi_k, f_1 \rangle
  \right) \times \\ 
  \times \langle \varphi_i, f_2 - f_1 \rangle \left( \prod_{i < k \le
      \ell } \langle \varphi_k, f_1 \rangle \right).
\end{multline*}
We deduce then 
\begin{eqnarray*}
&& \left| R^\ell[\varphi](f_2) - R^\ell[\varphi](f_1) \right| \le 
\| \varphi \|_{1,\FF \otimes (L^\infty) ^{\ell-1}} \, \|f_2 - f_1 \|_\GG , 
\\
&& \left| DR^\ell[\varphi](f_1)(h) \right| \le 
\| \varphi \|_{1,\FF \otimes (L^\infty) ^{\ell-1}} \, \|h \|_\GG, \\ 
&&  \left|R^\ell[\varphi](f_2) - R^\ell[\varphi](f_1) - 
DR^\ell[\varphi](f_1)(f_2 - f_1)\right|\le 
\| \varphi \|_{1,\FF^2 \otimes (L^\infty) ^{\ell-2}} \, \|f_2 - f_1 \|^2_\GG ,  
\end{eqnarray*}
where
\begin{eqnarray*}
\| \varphi \|_{1,\FF^k \otimes (L^\infty) ^{\ell-k}} 
&:=&
\sum_{  \{i_1, \dots,  i_k \} \subset \{1,...,\ell\}}
\| \varphi_{i_1} \|_{\FF} \, \cdots \| \varphi_{i_k}  \|_{\FF} \, \prod_{j \neq
      (i_1,\dots,i_k)} \| \varphi_j \|_{L^\infty(E)} \\
&\le&
\left\{
  \begin{array}{lcl}
     \ell \,  \| \varphi \|_{\infty,\FF \otimes (L^\infty) ^{\ell-1}}
     &  \qquad \hbox{for} \quad&  k=1, \vspace{0.2cm} \\
     \displaystyle \frac{\ell (\ell -1)} 2 \, \| \varphi \|_{\infty,\FF^2 \otimes
       (L^\infty) ^{\ell-2}}& \qquad \hbox{for} \quad&  k=2,
  \end{array}
\right.
\end{eqnarray*}
and we have defined 
\begin{eqnarray*}
\| \varphi \|_{\infty,\FF^k \otimes (L^\infty) ^{\ell-k}} 
&:=& \max_{i_1, \dots,  i_k \mbox{ {\tiny distinct in }} [|1,\ell|]} 
\| \varphi_{i_1} \|_{\FF} \, \cdots \| \varphi_{i_k} \|_{\FF} \! 
\prod_{j \neq (i_1,\dots,i_k)} \| \varphi_j \|_{L^\infty(E)} \\
&\le& \| \varphi \|_{\FF^{\otimes \ell}},
\end{eqnarray*}
since $\| \cdot \|_{L^\infty(E)} \le \| \cdot \|_\FF$.  This proves
that $R^\ell[\varphi] \in C^{1,1}(\PP_\GG(E),\R)$. The cases $k \ge
2$ are proved similarly. 
\end{proof}

\section{Assumptions and technical lemmas}
\label{sec:lemmas}
\setcounter{equation}{0}
\setcounter{theo}{0}

In this section we collect the lemmas used in the 
proof of Theorem~\ref{theo:abstract}, and the technical assumptions
that are needed.

The first assumption is simply the statement that the $N$-particle
dynamics is well defined and invariant under permutation of the particles. 

\medskip

\fbox{
\begin{minipage}{0.9\textwidth}
\begin{itemize}
\item[{\bf (A1)}] {\bf $N$-particle semigroup.}  The family $T^N_t$
  consists of strongly continuous Markov semigroups on $C_b(E^N) $
  that are invariant under permutations of indices.  We denote their
  generator by $G^N$ and we denote by $S^N_t$ the dual semigroup on
  the $N$-particle distributions.
     \end{itemize}
\end{minipage}
 }
\medskip

\subsection{The generator of the pullback semigroup}
\label{sec:calculus-gen}

While the definition of $T^N_t$ and its generator is rather standard,
it takes more care when defining the pullback of the nonlinear
semigroup and the corresponding generator. The second assumption
that we need to impose on the system is related to this, and relates
to our definition of a differential calculus of functions on $\PP(E)$.

\bigskip

\fbox{
\begin{minipage}{0.9\textwidth}
\begin{itemize} 
\item[{\bf (A2)}] {\bf Nonlinear semigroup.} Consider a probability
  space $\PP_{\GG_1}(E)$ (defined in Definition~\ref{defGG}) associated
  to a weight function $m_{\GG_1}$, endowed with the metric induced
  from a normed space $\GG_1$, and with bounded diameter.  Assume that
  for any $\tau > 0$ we have:
  \begin{itemize}
  \item[(i)] The equation \eqref{eq:0020a} generates a continuous
    semigroup $S^{N\!L}_t$ on $\PP_{\GG_1}(E)$ which is  uniformly
    Lipschitz continuous: there exists $C_\tau>0$ such that
    \[ 
    \forall \, f, g \in \PP_{\GG_1}(E), \quad \sup_{t \in [0,\tau]}
    \mbox{dist}_{\GG_1}\left(S^{N\!L}_t f, S^{N\!L}_t g\right) \le C_\tau \,
    \mbox{dist}_{\GG_1}(f,g).
  \]
  
\item[(ii)] There exists $\delta \in (0,1]$ such that the (possibly
  nonlinear) generator $Q$ (introduced in equation \eqref{eq:0020a})
  is  bounded and   $\delta$-H\"older continuous from
  $\PP_{\GG_1}(E)$ into $\GG_1$ in the following sense:  there exist $L, K >0$ so that for any $f, g \in 
  \PP_{\GG_1}(E)$ 
    $$
    \| Q(f) \|_{\GG_1} \le K, \quad \| Q(f) - Q(g) \|_{\GG_1} \le L \, \| f - g
    \|^\delta_{\GG_1}.
    $$
\end{itemize}
\end{itemize}
\end{minipage}
}
\bigskip

\Black

This assumption is sufficient for defining the generator of $T^{\infty}_t$:

\begin{lem}\label{lem:nonlgenerator} 
  Under assumption {\bf (A2)} the pullback semigroup $T^\infty_t$
  is a contraction semigroup on the Banach space $UC(P_{\GG_1}(E);\R)$
  and its generator $G^\infty$ is an unbounded linear operator on $UC(
  P_{\GG_1}(E);\R)$ with domain $\hbox{\emph{Dom}}(G^\infty)$ containing
  $C^{1,1}( P_{\GG_1}(E); \R)$. It is defined by 
  \begin{equation}
    \label{eq:formulaGinfty}
    \forall \, \Phi \in C^{1,1}( P_{\GG_1}(E); \R), \ 
    \forall \, f \in P_{\GG_1}(E), \quad \left( G^\infty \Phi \right) (f) =
    \left \langle D\Phi[f], Q(f)\right\rangle.
  \end{equation}
\end{lem}

\begin{proof}
The proof is split into several steps.

\smallskip \noindent {\sl Step 1.}  We claim that for any $f_{\mbox{{\tiny in}}} \in
\PP_{\GG_1}(E)$ and $\tau > 0$ the map 
\[
\SS(f_{\mbox{{\tiny in}}}) : [0,\tau) \to \PP_{\GG_1}(E), \quad t \mapsto S^{N \! L}_t(f_{\mbox{{\tiny in}}})
\]
is right-differentiable at $t=0^+$ with $\SS(f_{\mbox{{\tiny in}}})'(0) =
Q(f_{\mbox{{\tiny in}}})$. 

Denote $f_t := S^{N\!L}_t f_{\mbox{{\tiny in}}}$.  First, since
$Q(f_t)$ is bounded in $\GG_1$ uniformly on $t \in [0,\tau]$ from {\bf
  (A2)-(ii)}, we have, 
%
uniformly on $f_{\mbox{{\tiny in}}}\in \PP_{\GG_1}(E)$,
\begin{equation}\label{eq:ft-f0A2} \| f_t - f_{\mbox{{\tiny in}}}
  \|_{\GG_1}= \left\|\int_0^t Q(f_s) \, {\rm d}s \right\|_{\GG_1} \le K \,
  t, \end{equation} 
and then using {\bf (A2)-(ii)} and the inequality
\eqref{eq:ft-f0A2} we obtain
\begin{eqnarray*} \| f_t - f_{\mbox{{\tiny in}}} - t \, Q(f_{\mbox{{\tiny in}}})\|_{\GG_1} &=&
  \left\|\int_0^t \left(Q(f_s) - Q(f_{\mbox{{\tiny in}}})\right) \, {\rm d}s
  \right\|_{\GG_1}
  \\
  &=& L \, \int_0^t \left\| f_s - f_{\mbox{{\tiny in}}} \right\|_{\GG_1}^\delta \,
  {\rm d}s
  \\
  &\le& L \, \int_0^t (K \, s)^\delta \, ds = L \, K^\delta \,
  {t^{1+\delta} \over 1+\delta}, \end{eqnarray*} 
which implies the claim.

\smallskip \noindent {\sl Step 2.} We claim that $(T^\infty_t)$ is a
$C_0$-semigroup of linear and bounded (in fact contraction) operators
on $UC(\PP_{\GG_1}(E);\R)$.  Indeed, first for any $\Phi \in UC(
\PP_{\GG_1}(E);\R)$ and denoting by $\omega_\Phi$ the modulus of
continuity of $\Phi$, we have \begin{eqnarray*} \left|(T^\infty_t
    \Phi)(g) - (T^\infty_t \Phi)(f)\right|
  &=& \left| \Phi(S^{N\!L}_t(g)) - \Phi(S^{N\!L}_t(g)) \right| \\
  &\le& \omega_\Phi \left(\mbox{dist}_{\GG_1}
    \left(S^{N\!L}_t(g),S^{N\!L}_t(f)\right)\right)
  \\
  &\le& \omega_\Phi \left(C_\tau \, \mbox{dist}_{\GG_1} (f,g)\right)
\end{eqnarray*} so that $T^\infty_t \Phi \in UC( \PP_{\GG_1}(E);\R)$ for any $t
\ge 0$.  Next, we have
$$
\| T^\infty_t \| = \sup_{\|\Phi \| \le 1} \| T^\infty_t \Phi \| =
\sup_{\|\Phi \| \le 1} \sup_{f \in \PP_{\GG_1}(E) } \left|
\Phi(S^{N\!L}_t(f))\right| \le 1, \quad \|\Phi \| = \sup_{h \in
  \PP_{\GG_1}(E) } |\Phi(h)|.
$$
Finally, from \eqref{eq:ft-f0A2}, for any $ \Phi \in
UC(\PP_{\GG_1}(E);\R) $, we have
$$
\left\|T^\infty_t \Phi - \Phi \right\| = \sup_{f \in \PP_{\GG_1}(E) }
\left|\Phi(S^{N\!L}_t(f)) - \Phi(f)\right| \le \omega_\Phi (K \, t)
\to 0 \quad \mbox{ as } t \to 0^+.
$$
As a consequence, Hille-Yosida Theorem (see for instance
\cite[Theorem~3.1]{MR710486}) implies that $(T^\infty_t)$ is
associated to a closed generator $G^\infty$ with dense domain
$\hbox{dom}(G^\infty) \subset UC(\PP_{\GG_1}(E);\R)$.

\smallskip \noindent {\sl Step 3.} A candidate for
this generator is defined as follows.  Let
$\tilde G^\infty$ be defined by
$$ 
\forall \, \Phi \in C^{1,1}(\PP_{\GG_1}(E); \R), \ \forall \, f \in \PP_{\GG_1}(E),
\quad ( \tilde G^\infty \Phi ) (f) := \left \langle D\Phi[f],
  Q(f)\right\rangle.
$$
The RHS is well defined since $D\Phi(f) \in \BB(\GG_1,\R) = \GG'_1$
and $Q(f) \in \GG_1$ by assumption. Moreover, since both $f \mapsto
D\Phi[f]$ and $f \mapsto Q(f)$ are uniformly continuous 
so is the
map $f \mapsto (\tilde G^\infty \Phi) (f)$. It yields $\tilde G^\infty
\Phi \in UC(\PP_{\GG_1}(E); \R)$. 

\smallskip \noindent {\sl Step 4.} Finally, by composition, 
\[
\forall\, f \in \PP_{\GG_1}(E)(E), \quad t \mapsto T^\infty _t \Phi (f) = \Phi
\circ S^{N\!L}_t (f)
\]
is right-differentiable at $t=0^+$ and
\begin{eqnarray*} 
  {{\rm d} \over {\rm d}t} (T^\infty_t \Phi) (f) |_{t=0}
  &:=&  {{\rm d} \over {\rm d}t} (\Phi \circ \SS(f)(t)) |_{t=0} \\
  &=& \left\langle D\Phi (\SS(f)(0)), {{\rm d} \over {\rm d}t} \SS(f)(0) \right\rangle \\
  &=& \left\langle D\Phi[f], Q(f) \right\rangle = \left( \tilde G^\infty \Phi
  \right) (f),
\end{eqnarray*}
which implies $\Phi \in \hbox{Dom}(G^\infty)$ and
\eqref{eq:formulaGinfty}.
\end{proof}

\subsection{Estimates in the proof of Theorem~\ref{theo:abstract}}

The proof of Theorem~\ref{theo:abstract} relies on the three following
lemmas, together with their assumptions. 

\subsubsection{Estimate of $\TT_1$}
\label{sec:estimate-tt_1}

The first error term is estimated by the following combinatorics
argument. 
\begin{lem}
\label{lem:symmetrization}
Let $S^N_t$, $\mu^N_Z$ and $R^{\ell}[\varphi]$ as defined in
Section~\ref{sec:abstractchaos} (see Figure~\ref{fig:diagram}), and let
$\varphi\in C^b(E^{\ell})$. For any $t\ge 0$ and any $N \ge 2 \ell$
\begin{equation}\label{estim:T1}
  \TT_1 := \left| \left\langle S^N_t(f_{\mbox{{\tiny \emph{in}}}}^{\otimes N}), 
      \varphi \otimes
      1^{\otimes N-\ell} \right\rangle - \left \langle S^N_t(f_{\mbox{{\tiny \emph{in}}}}^{\otimes N}),
      R^\ell[\varphi] \circ \mu^N_Z \right\rangle \right| \le \frac{2 \,
    \ell^2 \, \| \varphi \|_{L^\infty(E^\ell)}}{N}.
\end{equation} 
\end{lem}

\begin{proof}
Since $S^N_t(f_{\mbox{{\tiny in}}}^{\otimes N})$ is a symmetric probability measure, estimate
(\ref{estim:T1}) is a direct consequence of the following estimate:
For any $\varphi \in C_b(E^\ell)$ and any  $  \, N \ge 2 \ell$ we have
\begin{equation}\label{estim:symmetrization1} 
  \left| \left( \varphi 
      \otimes {\bf 1}^{\otimes N-\ell} \right)_{\mbox{{\tiny sym}}} -
    \pi_N R^\ell[\varphi] \right| 
  \le \frac{2 \, \ell^2 \, \| \varphi \|_{L^\infty(E^\ell)}}{N}\,.
\end{equation}
Here the symmetrized version of a function  $\phi \in C_b(E^N)$, is
defined as
\begin{equation}\label{def:sym}
  \phi_{\mbox{{\tiny sym}}} = 
  \frac{1}{|\SSS_N|} \, \sum_{\sigma \in \SSS_N} \phi_\sigma.
\end{equation}
As a consequence for any symmetric measure $f^N \in \PPS(E^N)$ we have 
\begin{equation} \label{eq:mNRphi} \left| \langle f^N,
    R^\ell[\varphi](\mu^N _Z) \rangle - \langle f^N, \varphi \rangle
  \right| \le {2 \, \ell^2 \, \| \varphi \|_{L^\infty(E^\ell)} \over
    N}.
\end{equation}
To establish the inequality~(\ref{estim:symmetrization1}), we let 
$\ell \le N/2$ and  introduce
\[ 
A_{N,\ell} := \left\{ (i_1, \dots, i_\ell) \in \{1,\dots,N\}^\ell \, : \
  \forall \, k \not= k', \ i_k \not= i_{k'} \ \right\} \quad
\mbox{and} \quad B_{N,\ell} := A_{N,\ell}^c.
\]
Since there are $N (N-1) \dots (N-\ell+1)$ ways of choosing $\ell$
distinct indices among $\{1,\dots,N\}$ we get
\begin{eqnarray*}
  {\left|B_{N,\ell}\right|  \over N^\ell} &=& 
  1 - \left(1 - {1 \over N}\right) \, \dots \, 
  \left(1 - {\ell-1 \over N}\right) 
  = 1 - \exp \left( \sum_{i = 0}^{\ell-1} 
    \ln \left(1 - \frac{i}N \right) \right) \\
  &\le& 1 - \exp \left( - 2 \sum_{i = 0}^{\ell-1} \frac{i}N \right) 
  \le {\ell^2 \over  N},
\end{eqnarray*}
where we have used 
\[
\forall \, x \in [0,1/2], \quad \ln ( 1 - x) \ge - 2 \, x \qquad \mbox{and}
 \qquad \forall \, x \in \R, \quad e^{-x} \ge 1 - x.
\]
Then we compute
\begin{eqnarray*} &&R^\ell[\varphi](\mu^N_Z) =
  {1 \over N^\ell} \sum_{i_1, \dots, i_\ell = 1}^N 
 \varphi (z_{i_1}, \dots, z_{i_\ell}) \\
  &&= {1 \over N^\ell} \sum_{(i_1, \dots, i_\ell) \in A_{N,\ell}} 
  \varphi (z_{i_1}, \dots, z_{i_\ell})
  +  {1 \over N^\ell} \sum_{(i_1, \dots, i_\ell) \in B_{N,\ell}}  
   \varphi (z_{i_1}, \dots, z_{i_\ell}) \\
  &&= {1 \over N^\ell} \, {1 \over (N-\ell)!} \sum_{\sigma \in \SSS_N}
  \varphi (z_{\sigma(1)}, \dots, z_{\sigma(\ell)})
  +  {\mathcal O} \left(  {\ell^2 \over N} \, \|\varphi \|_{L^\infty} \right) \\
  &&= {1 \over N!} \sum_{\sigma \in \SSS_N} 
  \varphi (z_{\sigma(1)}, \dots , z_{\sigma(\ell)})
  +  {\mathcal O} \left(  {2 \, \ell^2 \over N} \, 
    \| \varphi \|_{L^\infty} \right) \\
\end{eqnarray*}
and the proof of (\ref{estim:symmetrization1}) is complete. Next for any
$f^N \in \PP(E^N)$ we have
$$
\left\langle f^N, \varphi \right\rangle = \left\langle f^N,
  \left(\varphi \otimes {\bf
      1}^{\otimes N-\ell} \right)_{\mbox{{\tiny sym}}} \right\rangle,
 $$
 and \eqref{eq:mNRphi} follows from (\ref{estim:symmetrization1}).
\end{proof}
 
\subsubsection{Estimate of $\TT_2$}
\label{sec:estimate-tt_2}

The second error term is estimated thanks to a consistency result for
the generators of the $N$ particle system and the limiting dynamics,
and stability estimates on the limiting dynamics. To proceed we need to
introduce the two corresponding assumptions.

\bigskip

\fbox{
\begin{minipage}{0.95\textwidth}
\begin{itemize}
\item[{\bf (A3)}] {\bf Convergence of the generators.}  Let
  $\PP_{\GG_1}(E)$ be the probability space endowed with a metric
  considered in {\bf (A2)}. For $k=1$ or $2$ there is a
  function $\varepsilon(N)$ with
  $\lim_{N\rightarrow\infty}\varepsilon(N)=0$ such that for all  $N
  \in \N$: and all $ \Phi \in
  C^{k,1}(\PP_{\GG_1}(E);\R), $ the  generators 
  $G^N$ satisfy 
  \begin{equation}
   \label{eq:estimcvg1}
   \left\|   G^N \, (\pi_N \, \Phi) - \pi_N \, G^\infty(\Phi)  \right\|_{L^\infty(E^N)} \le \eps(N) \, 
   \| \Phi \|_{C^{k,1}(\PP_{\GG_1}(E))}.
  \end{equation}
  where $Q$ is the nonlinear operator involved in equation
  \eqref{eq:0020a}, and $G^\infty (\Phi)= \langle Q ,D\Phi \rangle$ is the
  generator of the pullback semigroup $T^{\infty}$.  
\end{itemize}
\end{minipage}
}
\bigskip

\fbox{
\begin{minipage}{0.95\textwidth}
\begin{itemize}  
\item[{\bf (A4)}] {\bf Differential stability of the limiting
    semigroup.}  We assume that the flow $S^{N \! L}_t$ is
  $C^{k,1}(\PP_{\GG_1}(E),\PP_{\GG_2}(E))$ in the sense that there
  exists $C_{T} >0$ such that
  \begin{equation}
     \label{eq:stab-lim1}
    \int_0^T \left\| S^{N\! L}_t
    \right\|_{C^{k,1}(\PP_{\GG_1}(E),\PP_{\GG_2}(E))} \, {\rm d}t \le C_{T}
  \end{equation}
  where $\PP_{\GG_2}(E)$ is the same subset of probabilities as
  $\PP_{\GG_1}(E)$, endowed with the norm associated to a normed space
  $\GG_2 \supset \GG_1$ possibly larger than $\GG_1$ and where $k$ is
  the same integer as above.
\end{itemize}
\end{minipage}
}

\begin{lem}
\label{lem:steptwo}
Suppose that the assumptions {\bf (A1)} to {\bf (A4)} are satisfied,
and let $\TT_2$ be as defined in equation~(\ref{eq:cvtabst10}). Then

 \begin{eqnarray}\label{estim:T2}
\quad \TT_ 2 &:= & \left| \left\langle f_{\mbox{{\tiny in}}}^{\otimes N},
     T^N_t ( R^\ell[\varphi] \circ \mu^N_Z) \right\rangle -
   \left\langle f_{\mbox{{\tiny in}}}^{\otimes N}, \left( (T_t ^\infty R^\ell[\varphi] ) \circ
       \mu^N_Z \right) \right\rangle \right| 
\\ \nonumber &\le&
 C(k,\ell) \, C_{T} \, \eps(N) \, \| \varphi \|_{\FF_1 ^k \otimes
   (L^\infty)^{\ell-k}} 
\end{eqnarray}
for an  explicitly given constant $C(k,\ell)$ depending only on $k$
and $\ell$.

\end{lem}

\begin{proof}
We start from the following identity
\begin{multline*}
T^N_t \pi_N - \pi_N T^\infty _t = - \int_0 ^t \frac{\partial}{\partial
  s} \left( T^N _{t-s} \, \pi_N \, T^\infty _s \right) \, {\rm d}s \\ =
\int_0 ^t T^N _{t-s} \, \left[ G^N \pi_N - \pi_N G^\infty \right] \,
T^\infty _s \, {\rm d}s\,,
\end{multline*}
which we evaluate on $\Phi=R^{\ell}[\phi]\in C_b(\PP(E))$.
From assumption {\bf (A3)} we have for any $t \in [0,T]$
\begin{eqnarray}\nonumber
  &&
\left| \left\langle f_{\mbox{{\tiny in}}}^{\otimes N},
T^N_t \pi_N R^\ell[\varphi] - \pi_N T^\infty _t R^\ell[\varphi]  \right\rangle \right| 
     \\
  \nonumber
  &&\qquad = \left| \int_0 ^t  \left\langle S^N _{t-s}\left(f_{\mbox{{\tiny in}}} ^N\right), 
    \left[ G^N \pi_N - \pi_N G^\infty \right] \, 
    (T^\infty_s R^\ell[\varphi]) \right\rangle  \, {\rm d}s \right|
  \\ \nonumber
  &&\qquad \le \int_0 ^T \left\| \left[ G^N \pi_N - \pi_N G^\infty \right] \, 
    (T^\infty_s R^\ell[\varphi]) \right\|_{L^\infty(E^N)} \, {\rm d}s \\ \label{eq:trotter}
  &&\qquad \le \eps(N) \,  \int_0 ^T 
  \left\| T^\infty_s R^\ell[\varphi]\right\|_{C^{k,1}(\PP_{\GG_1}(E))} \, {\rm d}s.
\end{eqnarray}
Since $T^\infty_t(R^\ell[\varphi]) = R^\ell[\varphi] \circ S^{N\!L}_t$
with $S^{N\!L}_t \in C^{k,1}(\PP_{\GG_1}(E),\PP_{\GG_2}(E))$ thanks to assumption
{\bf (A4)} and $R^\ell[\varphi] \in C^{k,1}(\PP_{\GG_2}(E),\R)$ because
$\varphi \in \FF_2^{\otimes\ell}$ (see
Subsection~\ref{sec:compatibility}), we obtain with the help of
Lemma~\ref{lem:DL} that $T^\infty_t(R^\ell[\varphi]) \in
C^{k,1}(\PP_{\GG_1}(E))$ with uniform bound. We hence conclude that
\begin{equation}
  \label{estim:TtPhi}
  \int_0^T \left\| T^\infty_s(R^\ell[\varphi]) \right\|_{C^{k,1} (\PP_{\GG_1}(E))} \, {\rm d}s
  \le C(k,\ell) \, C_{T} \, \left\| R^\ell[\varphi] \right\|_{C^{k, 1}(P_{\GG_2})},
\end{equation}
where $C(k,\ell) \le \ell^2$ since $k=1$ or $k=2$. 

Going  back to the computation \eqref{eq:trotter}, and plugging
(\ref{estim:TtPhi}) we deduce (\ref{estim:T2}).
\end{proof}

\subsubsection{Estimate of $\TT_3$}
\label{sec:estimate-tt_3}

The error term $\TT_3$ from equation~(\ref{eq:cvtabst10}) depends on
estimating how well the initial data for the nonlinear equation
$f_{\mbox{\tiny in}}$ can be approximated by an empirical measure, and how
well this error is then propagated along the semigroup. To this
purpose we need to make a stability assumption for the limiting
semigroup. 

\bigskip

\fbox{
\begin{minipage}{0.95\textwidth}
\begin{itemize}   
\item[{\bf (A5)}] {\bf Weak stability of the limiting semigroup.}  There is  probability space $P_{\GG_3}(E)$
  (corresponding to to a weight function $m_{\GG_3}$ and  metric
  $\hbox{dist}_{\GG_3}$) such that  every $T > 0$ there exists a constant
  $\tilde C_{T} >0$ such that
 \begin{equation}
  \forall  \, f_1, f_2 \in P_{\GG_3}(E), \quad 
  \label{eq:stab-weak-S}
  \sup_{[0,T)} \mbox{dist}_{\GG_3} 
    \left( S^{N \! L}_t (f_1), S^{N \! L}_t(f_2) \right) \le \tilde C_{T} \,  \mbox{dist}_{\GG_3} (f_1,f_2).
  \end{equation}
\end{itemize}
\end{minipage}
}
\smallskip

\begin{rem}
  Observe that when $\PP_{\GG_1}(E) = \PP_{\GG_1}(E) = \PP_{\GG_3}(E)$
  with the same weights and distances, the assumption {\bf (A5)} is
  included in {\bf (A4)}. However it is crucial to have the flexibility
  to play with different metric structures in these assumptions.
\end{rem}

\begin{lem}
\label{lem:stepthree}
Assume that the limiting semigroup $S^{N\!L}_t$ satisfies assumption
{\bf (A5)} for some probability space $\PP_{\GG_3}(E)$. Let $\FF_3$
satisfies a duality inequality with $\PP_{\GG_3}(E)$ as defined in
Definition~\ref{def:DualitySpaces}. Let $\varphi\in C_b(E^{\ell})$.

Then for any $f_{\mbox{\tiny \emph{in}}}\in\PP(E)$, $t>0$ and $N\ge 2\ell$ we
have
\begin{multline}\label{estim:T4}
  \TT_3 := 
  \left| \left\langle f_{\mbox{{\tiny \emph{in}}}}^{\otimes N}, 
    \left(T_t ^\infty R^\ell[\varphi] \right) \circ \mu^N_Z \right\rangle
    - \left\langle \left(S^{N\!L}_t (f_{\mbox{{\tiny \emph{in}}}})\right)^{\otimes k} , 
                \varphi \right\rangle \right|  \\ \le 
  \ell \, \tilde C_{T} \,  \Omega_N ^{\GG_3} (f_{\mbox{{\tiny \emph{in}}}}) \, 
  \| \varphi \|_{\FF_3 \otimes (L^\infty)^{\ell-1}},
\end{multline}
where $\Omega_N ^{\GG_3} (f_{\mbox{{\tiny \emph{in}}}}) $ is defined in
\eqref{def:OmegaG3N} and $ \| \varphi \|_{\FF_3 \otimes
  (L^\infty)^{\ell-1}}$ is defined in the proof of
Lemma~\ref{lem:polyCk}.
\end{lem}

\begin{proof}
We split $\TT_3$ in two terms, the first one being
\begin{eqnarray*}
  \TT_{3,1} &:=& \left\langle f_{\mbox{{\tiny in}}}^{\otimes N}, \left(T_t ^\infty R^\ell[\varphi]
    \right) \circ \mu^N_Z \right\rangle \\
  &=& 
  \int_{E^N} R^\ell[\varphi] \left(S^{N\!L} _t \left(\mu^N _Z\right)\right)
  \,  f_{\mbox{{\tiny in}}}({\rm d}z_1) \, \dots \, f_{\mbox{{\tiny in}}}({\rm d}z_N) \\
  &=&  \int_{E^N}  \left( \prod_{i=1}^\ell a_i (Z) \right)  \, 
f_{\mbox{{\tiny in}}}({\rm d}z_1) \, \dots \, f_{\mbox{{\tiny in}}}({\rm d}z_N),
\end{eqnarray*}
with 
$$
\forall \, i=1, \dots, \ell, \quad a_i = a_i(Z) := \int_{E}
\varphi_i(w) \, S^{N\!L} _t(\mu^N_Z)({\rm d}w). 
$$

Similarly, we write for the second term
\begin{eqnarray*}
\TT_{3,2} 
&=& \left\langle \left( S^{N\!L} _t  (f_{\mbox{{\tiny in}}}) \right)^{\otimes \ell}, 
 \varphi  \right\rangle
=  \int_{E^N}  \left( \prod_{i=1}^\ell b_i  \right)\, 
f_{\mbox{{\tiny in}}}({\rm d}z_1) \, \dots \, f_{\mbox{{\tiny in}}}({\rm d}z_N),
\end{eqnarray*} 
with 
$$
\forall \, i=1, \dots, \ell, \quad b_i := \int_{E} \varphi_i(w) \,
S^{N\!L} _t(f_{\mbox{{\tiny in}}})({\rm d}w).
$$
Using the identity 
\[
\prod_{i=1} ^\ell a_i - \prod_{i=1} ^\ell b_i = \sum_{i=1} ^\ell a_1 \dots a_{i-1}
\, (a_i - b_i) \, b_{i+1} \dots b_\ell,
\]
we get 
\begin{equation}\label{estim:T4-3}
\TT_3  \le 
\sum_{i=1} ^\ell \left( 
\prod_{j \neq i} \| \varphi_j \|_{L^\infty(E)} \right) \, 
\int_{E^N} \left| a_i(Z) - b_i \right| \, f_{\mbox{{\tiny in}}}({\rm d}z_1) \, \dots \, f_{\mbox{{\tiny in}}}({\rm d}z_N).
\end{equation}

Then by using the duality bracket together with assumption {\bf (A5)}
we have 
\begin{eqnarray} \nonumber \left| a_i(Z) - b_i \right| &:=& \left| \int_{E}
    \varphi_i(w) \, \left(S^{N\!L} _t(f_{\mbox{{\tiny in}}})({\rm d}w) - S^{N\!L}_t(\mu^N
      _V)({\rm d}w) \right)\right| \\ \nonumber &\le& \| \varphi_i \|_{\FF_3} \,
  \mbox{dist}_{\GG_3} \left(S^{N\!L} _t(f_{\mbox{{\tiny in}}}), S^{N\!L} _t(\mu^N
    _Z)\right) \\ \label{estim:T4-4} 
    &\le& \tilde C_{T} \, \| \varphi_i
  \|_{\FF_3} \, \mbox{dist}_{\GG_3}\left(f_{\mbox{{\tiny in}}}, \mu^N _Z\right).
\end{eqnarray} 
Therefore combining (\ref{estim:T4-3}) and (\ref{estim:T4-4}) (for any
$1 \le i \le \ell$), we conclude that (\ref{estim:T4}) holds.
\end{proof}

In order to use Lemma~\ref{lem:stepthree} we need an estimate on the
term $\Omega^{\GG_3}_N(f_{\mbox{{\tiny in}}}) $.  This information is provided by the
following quantitative version of the law of large number for
empirical measures taken from \cite{BookRachev}. We refer to
\cite{MM**} for a more detailed discussion of this issue. 
\begin{lem}\label{lem:Rachev} 
  For any $f_{\mbox{{\tiny \emph{in}}}} \in \PP_{d+5}(\R^d)$ and any $N \ge
  2$ there exists a constant $C$ which only depends on $d$ and
  $M_{d+5}(f_{\mbox{{\tiny \emph{in}}}})$ so that
  $$
  \Omega^{W^2_2}_N(f_{\mbox{{\tiny \emph{in}}}}) = 
  \int_{\R^{dN}} W_2(\mu^N_Z,f_{\mbox{{\tiny \emph{in}}}})^2 \, f^{\otimes N}_{\mbox{{\tiny \emph{in}}}}({\rm d}Z) 
  \le C \, N^{ - \frac{2}{d+4}}.
  $$
\end{lem}

 
\subsection{A remark on assumption {\bf (A4)}}
\label{sec:remarksMcK}

In this section we briefly explain how our key estimate {\bf
  (A4)} can be obtained in the case of a nonlinear operator $Q$ which
splits into a linear part and a bilinear part:
\begin{eqnarray}
  \forall f\in \PP(E), \quad Q(f) = Q_1(f) + Q_2(f,f)
\end{eqnarray}
with $Q_1$ linear and $Q_2$ bilinear symmetric. 

For two initial data $f_{\mbox{{\tiny in}}}$ and $g_{\mbox{{\tiny in}}}$ in a space $P_\GG(E)$ of
probability measures,  and some initial data $h_{\mbox{{\tiny in}}} \in \GG$ we introduce
the following evolution equations,
  $$
  \left\{ 
      \begin{array}{l}
        \partial_t g = Q(g) = Q_1(g) + Q_2(g,g), \quad
        g_{|t=0} = g_{\mbox{{\tiny in}}} , \vspace{0.3cm} \\
        \partial_t f = Q(f) = Q_1(f) + Q_2(f,f) , \quad f_{|t=0}
        = f_{\mbox{{\tiny in}}}, \vspace{0.3cm} \\
        \partial_t h = DQ[f] (h) = Q_1(h) + 2 \, Q_2(f,h),
        \quad h_{|t=0} = h_{\mbox{{\tiny in}}} 
        \end{array}
        \right.
  $$

  In the third equation the solution $h_t$ depends linearly on
  $h_{\mbox{{\tiny in}}}$ (but also nonlinearly on $f_{\mbox{{\tiny in}}}$): it is formally the
  first-order variation of the semigroup, i.e. $D S^{N\!L}_t
  f_{\mbox{{\tiny in}}}(h_{\mbox{{\tiny in}}}) = h_t$.

  We now want to write a \emph{second}-order variation of the
  semigroup. To this purpose we conside $h_t$ and $\tilde h_t$ two
  solutions to the third equation above and write
  \begin{equation*}
    \partial_t r = DQ[f](r) + {1 \over 2} D^2Q(f) (h,\tilde h) 
        =
        Q_1(r) + 2 \, Q_2(f,r) + Q_2 (h,\tilde h), \quad
        r_{|t=0} = 0.
  \end{equation*}
  In this fourth equation $r_t$ depends bilinearly on $h_t$ and
  $\tilde h_t$ which are two solutions to the third equation, and
  therefore bilinearly on $h_{\mbox{{\tiny in}}}, \tilde
  h_{\mbox{{\tiny in}}}$ (it also depends nonlinearly on
  $f_{\mbox{{\tiny in}}}$ again): it is formally the second-order
  derivative of the semigroup $D^2 S^{N\!L}_t [f_{\mbox{{\tiny
        in}}}](h_{\mbox{{\tiny in}}}, \tilde h_{\mbox{{\tiny in}}}) =
  r_t$. Observe that the initial data $r_{\mbox{{\tiny in}}}$ are
  always zero for this second variation problem since the map
  $(f_t)_{t \ge 0} \mapsto f_{\mbox{{\tiny in}}}$ is linear.

  Consider now $h_{\mbox{{\tiny in}}} = \tilde h_{\mbox{{\tiny in}}} = g_{\mbox{{\tiny in}}} - f_{\mbox{{\tiny in}}}$ (which
  implies $h_t = \tilde h_t$). Let us define $\mathsf s := f+g$,
  $\mathsf d := g - f$, $\omega := g - f - h$, $\psi := g - f - h -
  r$, for which we get the following evolution equations
 $$
  \left\{ 
      \begin{array}{l}
        \partial_t \mathsf d = Q_1(\mathsf d) +  Q_2(\mathsf s,\mathsf d),
        \quad \mathsf d_{|t=0} = \mathsf d_{\mbox{{\tiny in}}} = g_{\mbox{{\tiny in}}} - f_{\mbox{{\tiny in}}}, \vspace{0.3cm} \\
        \partial_t \omega  =
        Q_1(\omega) + Q_2(\mathsf s,\omega) + Q_2 (h,\mathsf d), \quad
        \omega_{|t=0} = 0,\vspace{0.3cm} \\
        \partial_t \psi  =
        Q_1(\psi) + Q_2(\mathsf s,\psi) + Q_2 (h,\omega) +  Q_2 (r,\mathsf d), \quad
        \psi_{|t=0} = 0.
        \end{array}
        \right.
$$

Now we can translate the regularity estimates on $S^{N\!L} _t$ in terms
of estimates on these solutions on some given time interval $[0,T]$: 
\begin{equation*}
  \left\{ 
    \begin{array}{l}\displaystyle
      \sup_{t\in[0,T]} \|\mathsf d_t\|_{\mathcal{G}_2} \le
      C_T\|\mathsf d_{\mbox{{\tiny in}}}\|_{\mathcal{G}_1} \implies S^{N\!L}_t\in
      C^{0,1}(P_{\mathcal{G}_1}(E), P_{\mathcal{G}_2}(E)) \vspace{0.2cm} \\ \displaystyle
      \left\{ \begin{array}{l} 
          S^{N\!L}_t\in
      C^{0,1}(P_{\mathcal{G}_1}(E), P_{\mathcal{G}_2}(E))
      \vspace{0.2cm} \\ 
          \sup_{t\in[0,T]} \|h_t\|_{\GG_2} \le
      C_T \|h_{\mbox{{\tiny in}}}\|_{\mathcal{G}_1} \vspace{0.2cm} \\
      \sup_{t\in[0,T]} \|\omega_t\|_{\GG_2} \le
      C_T \|\mathsf d_{\mbox{{\tiny in}}}\|_{\mathcal{G}_1}^2 
    \end{array}
    \right\} \implies S^{N\!L}_t\in
      C^{1,1}(P_{\mathcal{G}_1}(E), P_{\GG_2}(E)) \vspace{0.2cm} \\ \displaystyle
      \left\{ 
        \begin{array}{l} 
          S^{N\!L}_t\in
      C^{1,1}(P_{\mathcal{G}_1}(E), P_{\GG_2}(E)) \vspace{0.2cm} \\
      \sup_{t\in[0,T]} \|r_t\|_{\GG_2} \le
      C_T \|h_{\mbox{{\tiny in}}}\|_{\mathcal{G}_1}  \|\tilde h_{\mbox{{\tiny in}}}\|_{\mathcal{G}_1}  \vspace{0.2cm} \\
          \sup_{t\in[0,T]} \|\psi_t\|_{\mathcal{G}_2} \le
      C_T \|\mathsf d_{\mbox{{\tiny in}}}\|_{\mathcal{G}_1}^3 
    \end{array}
    \right\} 
    \implies S^{N\!L}_t\in C^{2,1}(P_{\mathcal{G}_1}(E), P_{\mathcal{G}_2}(E)).
    \end{array}
  \right.
\end{equation*}
     
Such estimates are typically obtained by energy estimates for the
equations satisfied by $\mathsf d$, $r$, $\omega$ and $\psi$, for a
well chosen ``cascade'' of norms connecting $\| \cdot \|_{\GG_1}$ to
$\|\cdot \|_{\GG_2}$ (see later in the applications).



\section{Maxwell molecule collisions with cut-off}
\label{sec:BddBoltzmann}
\setcounter{equation}{0}
\setcounter{theo}{0}


\subsection{The model}

In this section we assume that $E = \R^d$, $d \ge 2$, and we consider
an $N$-particle system undergoing a space homogeneous random Boltzmann
collisions according to a collision kernel $b \in L^1([-1,1])$ only
depending on the deviation angle and locally integrable. This is
usually called \emph{Maxwellian molecules with Grad's angular
  cut-off}, as introduced in \cite{Kac1956,Kac1957,McKean1967}. We
make the normalization hypothesis 
\begin{equation*}
\|b \|_{L^1} = \int_{\Sp^{d-1}} b(\sigma_1) \, {\rm d} \sigma = 1.
\end{equation*}

Let us now describe the stochastic process. Since the phase space
$E^N$ corresponds to the velocities of the particles, we shall denote
$Z=V$ in this section. Given a pre-collisional $N$-system of velocity
particles $V = (v_1, \dots, v_N) \in E^N = (\R^d)^N$, the stochastic
runs as follows:
\begin{itemize}
\item[(i)] for any $i'\neq j'$, we draw randomly for the pair of
  particles $(v_{i'},v_{j'})$ a random time $T_{i',j'}$ of collision
  according to an exponential law of parameter $1$, and then choose
  the collision time $T_1$ and the colliding pair $(v_i,v_j)$ (which
  is a.s. well-defined) in such a way that
$$
T_1 = T_{i,j} := \min_{1 \le i' \neq j' \le N} T_{i',j'};
$$
\item[(ii)] we then choose $\sigma \in \Sp^{d-1}$ at random according to
  the law $b(\cos \theta_{ij})$ where we define the angular deviation
  $\theta_{ij}$ by $\cos \theta_{ij} = \sigma \cdot
  (v_j-v_i)/|v_j-v_i|$;
\item[(iii)] the new state after collision at time $T_1$ becomes
$$
V^* = V^*_{ij} = R_{ij,\sigma}V = (v_1, \dots, v^*_i, \dots., v^*_j, \dots ,
v_N),
$$
where the rotation $R_{ij,\sigma}$ on the $(i,j)$ pair with vector
$\sigma$ is defined by
\begin{equation}\label{vprimvprim*}
  \quad\quad   v^*_i = {w_{ij} \over 2} +  {u^*_{ij} \over 2}, \quad
  v^*_j= {w_{ij} \over 2} - {u^*_{ij} \over 2},
\end{equation}
with
$$
w_{ij} = v_i+v_j, \quad u^*_{ij} = |u_{ij}| \, \sigma, \quad u_{ij} =
v_i-v_j. 
$$
\end{itemize}

\smallskip Scaling the time by a factor $1/N$ and repeating the above
construction lead to the definition of a Markov process $(\VV^N_t)$ on
$(\R^d)^N$. It is associated to a Feller semigroup $(T^N_t)$ with
generator $G^N$. Moreover the \emph{master equation} on the law
$f^N_t$ is given in dual form by
\begin{equation}\label{eq:BoltzBddKolmo}
  \frac{{\rm d}}{{\rm d}t} \langle f^N_t,\varphi \rangle = \langle f^N_t, G^N \varphi \rangle 
\end{equation}
with 
\begin{equation}\label{defBoltzBddGN}
  (G^N\varphi) (V) = {1 \over N} \sum_{1\le i < j\le N}^N 
  \int_{\mathbb{S}^{d-1}} b(\cos\theta_{ij}) \, \left[\varphi^*_{ij} -
    \varphi\right] \, {\rm d}\sigma
\end{equation}
where $\varphi^*_{ij}= \varphi(V^*_{ij})$ and $\varphi = \varphi(V)
\in C_b(\R^{Nd})$.  Finally, the flow $f^N_{\mbox{{\tiny in}}} \mapsto f^N_t$ defines a
semigroup $S^N_t$ for the $N$-particle distributions which is nothing
but the dual semigroup of $T^N_t$.

Note that the collision process is invariant under permutation of the
velocities, and satisfies the microscopic conservations of momentum
and energy at any collision time
$$
\forall \, \alpha = 1, \dots, d, \quad \sum_{k} v^*_{k\alpha} = \sum_{k}
v_{k\alpha}, \qquad |V^*|^2 = |V|^2 := \sum_{k=1}^N |v_k|^2.
$$

We write $V = (v_i)_{1 \le i \le N} = (v_1, \dots, v_N)
\in E^N$ and $v = (v_\alpha)_{1\le \alpha \le d} \in \R^d$, so that $
V = (v_{i\alpha}) \in \R^{Nd}$ with $v_{i\alpha} \in \R$.

As a consequence, for any symmetric initial law $f_{\mbox{{\tiny in}}}^{\otimes N}
\in \PPS(\R^{Nd})$ the law density $f_t^N$ remains a symmetric
probability and conserves momentum and energy
$$
\left\{
\begin{array}{l}\displaystyle
  \forall \, \alpha = 1, \dots, d, \quad \int_{\R^{dN}} \left( \sum_{k=1}^N v_{k\alpha} \right) \, f_t^N ({\rm d}v) =
  \int_{\R^{dN}} \left( \sum_{k=1}^N v_{k\alpha} \right)  \, f_{\mbox{{\tiny in}}}^{\otimes N} ({\rm d}v), \vspace{0.3cm} \\ 
  \displaystyle
  \forall \, \theta : \R_+ \to \R_+, \quad  \int_{\R^{dN}}\theta( |V|^2 ) \, f_t^N ({\rm d}v) = \int_{\R^{dN}} \theta( |V|^2 )  \, f_{\mbox{{\tiny in}}}^{\otimes N} ({\rm d}v).
\end{array}
\right.
$$

The formal limit of this $N$-particle system is the nonlinear
homogeneous Boltzmann equation on $\PP(\mathbb{R}^d)$ defined by
\begin{equation}\label{eq:NLBoltzEq}
\frac{\partial}{\partial t} f_t = Q(f_t,f_t)
\end{equation}
where the quadratic Boltzmann collision operator $Q$ is defined by
\begin{equation}\label{eq:NLBoltzKernel}
\langle Q(f,f), \varphi \rangle := \int_{\R^{2d}\times \Sp^{d-1}}  
b(\theta) \, (\phi (w^*_2) - \phi (w_2)) \, {\rm d}\sigma \, f({\rm d}w_1) \,
f({\rm d}w_2)
\end{equation}
for  $\varphi \in C_b(\R^d)$ and $f \in \PP(\R^d)$, with 
\begin{equation}\label{eq:w1*w2*}
w_1 ^* = {w_1+ w_2 \over 2} + {|w_2 - w_1|\over 2}\, \sigma, \qquad 
w_2 ^* = {w_1+ w_2 \over 2} - {|w_2 - w_1|\over 2}\, \sigma 
\end{equation}
and $\cos\theta = \sigma \cdot (v-w)/|v-w|$. The equation
\eqref{eq:NLBoltzEq}-\eqref{eq:NLBoltzKernel} is the space homogeneous
Boltzmann equation for elastic collisions associated to the Maxwell
molecules cross section with Grad's cutoff. We refer to the textbooks
\cite{Ce88} and \cite{VillaniHandBook} and the numerous references
therein for both the physical background and the mathematical theory
of the Boltzmann equation.  
This equation generates a nonlinear semigroup $S^{N\!L}_t$ on
$\PP(\R^d)$ defined by $S^{N\!L}_t f_{\mbox{{\tiny in}}} := f_t$ for any $f_{\mbox{{\tiny in}}} \in
\PP(\R^d)$, which satisfies conservation of momentum and energy:
$$
\forall \, t \ge 0, \quad \int_{\R^{d}} v \, f_t ({\rm d}v) = \int_{\R^{d}}v \,
f_{\mbox{{\tiny in}}} ({\rm d}v), \quad \int_{\R^{d}} |v|^2 \, f_t ({\rm d}v) = \int_{\R^{d}}|v|^2 \,
f_{\mbox{{\tiny in}}} ({\rm d}v).
$$

\subsection{Statement of the result}
\label{sec:resultEBbounded} On the one hand, it is well known that for
the collision kernel that we have chosen 
the $N$-particle Markov process $(\VV^N_t)$ described above is well
defined for any initial velocity $\VV^N_0$, and in particular, for any
given initial law $f^N_{\mbox{{\tiny in}}} \in \PPS((\R^d)^N)$ there exists a unique
solution $f^N_t \in \PPS((\R^d)^N)$ to equations
(\ref{eq:BoltzBddKolmo})-(\ref{defBoltzBddGN}) so that the
$N$-particle semigroup $S^N_t$ is well defined, see
\cite{Kac1957,Kac1979,T1,Meleard1996}.  On the other hand, it is also
well known that for any $f_{\mbox{{\tiny in}}} \in \PP_q(\R^d)$, $q \ge 0$ the
nonlinear Boltzmann equation
(\ref{eq:NLBoltzEq})-(\ref{eq:NLBoltzKernel}) has a unique solution
$f_t \in \PP_q(\R^d)$. This solution conserves momentum and energy
as soon as $q \ge 2$, see for instance \cite{T1,ToscaniV,FWellPos,VillaniHandBook}.

\medskip Our mean field limit result then states as follows.

\begin{theo}[The Boltzmann equation for Maxwell molecules with Grad's
  cut-off]\label{theo:BddBoltz}
  Consider an initial distribution $f_{\mbox{{\tiny \emph{in}}}} \in \PP_q(\R^d)$, $q \ge 2$,
  the hierarchy of $N$-particle distributions $f^N _t = S^N
  _t(f_{\mbox{{\tiny \emph{in}}}} ^{\otimes N})$ following \eqref{eq:BoltzBddKolmo}, and the
  solution $f_t = S^{N\!L} _t(f_{\mbox{{\tiny \emph{in}}}})$ following \eqref{eq:NLBoltzEq}.

  Then there is a contant $C>0$ and, for any $T>0$, there are 
  constants $C_{T}, \tilde C _T >0$ such that for any
  \[ 
  \varphi = \varphi_1 \otimes  \dots \otimes \,
  \varphi_\ell \in \FF^{\otimes \ell}, \quad \FF := C_b(\R^d) \cap
  \mbox{{\em Lip}}(\R^d), 
  \quad \| \varphi_j \|_\FF \le 1, 
  \] 
  we have for $N \ge 2 \ell$: 
\begin{equation} \label{eq:cvgBddBE}
  \sup_{[0,T]}\left| \left \langle \left( S^N_t(f_{\mbox{{\tiny \emph{in}}}} ^N) - \left(
          S^{N\!L} _t (f_{\mbox{{\tiny \emph{in}}}}) \right)^{\otimes N} \right), \varphi
    \right\rangle \right| \le C \, \frac{\ell^2}{N} + C_{T} \, {\ell^2
    \over N}+ \tilde C_{T} \, 
  \ell \, \Omega^{W_2} _N (f_{\mbox{{\tiny \emph{in}}}})
\end{equation}
where $\Omega^{W_2} _N$ was defined in \eqref{def:OmegaG3N} and $W_2$
is the quadratic MKW distance defined in \eqref{eq:wasserstein}.

As a consequence of \eqref{eq:cvgBddBE} and Lemma~\ref{lem:Rachev},
this implies propagation of chaos with rate $\eps(N) \le
C(\ell,T,f_{\mbox{{\tiny \emph{in}}}}) \, N^{-{1 \over d+4}}$ for any
initial data $f_{\mbox{{\tiny \emph{in}}}} \in P_{d+5}(\R^d)$, where $
C(\ell,T,f_{\mbox{{\tiny \emph{in}}}})$ is an explicitly computable
constant.
\end{theo}

For the Boltzmann equation with bounded kernel, propagation of chaos
has been established by McKean in \cite{McKean1967}, where he adapted
the method introduced by Kac in \cite{Kac1957} based on the Wild sum
representation of the solutions to the Boltzmann equation. Gr\"unbaum
in \cite{Grunbaum} gave an alternative proof based on the same
``duality viewpoint'' as developed in the present paper.  Sznitman in
\cite{S1} also gave a proof of propagation of chaos based on a
nonlinear martingale approach.  In all these works, propagation of
chaos is proved but without any rate of convergence (as the number of
particles goes to infinity).  Graham and M\'el\'eard in
\cite{GrahamM1,GrahamM2,Meleard1996} were then able to prove the
propagation of the chaos with the sharp rate $C(\ell,T)/N$. Their
proof is based on the construction of a stochastic tree associated to
the process $\VV^N_t$ which is specific to the Boltzmann equation with
bounded kernel. More recently Kolokoltsov in \cite{Kolokoltsov} proved
a fluctuation estimate for similar processes using a ``duality view
point'' like the one developed by Gr\"unbaum and used also in our
work. His fluctuation estimate is similar to our Lemma~\ref{estim:T2}
(and of the same order), but pays less attention to the remaing terms
of our estimate.  Fournier and Godinho~\cite{FournierGodinho} prove
the propagation of chaos for a one-dimensional caricature of the
Boltzmann equation using a coupling method in the spirit of
\cite{T1,S6}. Their chaoticity estimate is of the same rate as ours.

\subsection{Proof of Theorem~\ref{theo:BddBoltz} }
\label{sec:proofBddBoltz}
The assumptions  {\bf (A1)-(A2)-(A3)-(A4)-(A5)} needed  to apply
Theorem~\ref{theo:abstract} will be verified step by step. In this
proof we fix $\FF_1=\FF_2=C_0(\R^d)$ and $\FF_3 = 
\mbox{Lip}(\R^d)$ and define $\PP_{\GG_1}(E) = \PP_{\GG_2}(E) :=
\PP(\R^d)$ endowed with the total variation norm $\|\cdot \|_{TV}$,
$\PP_{\GG_3}(E) := \PP_2(\R^d)$ endowed with the quadratic MKW
distance $W_2$. Notice that $(\GG_1,\FF_1)$ and $(\GG_2,\FF_2)$
satisfy a duality inequality of type 1, and $(P_{\GG_3}, \FF_3)$
satisfy a duality of type 2 (see Definition~\ref{def:duality}).

\smallskip
\noindent {\bf Proof of (A1). } The symmetry assumption is satisfied
because of the well-known properties of the Boltzmann-Kac $N$-particle
system, and we refer to the previous works
\cite{McKean1967,Grunbaum,GrahamM1,GrahamM2,Meleard1996} for details.


\smallskip
\noindent {\bf Proof of (A3). }
We claim that there exists $C_1\in \R_+$ such that for all $\Phi \in
C^{1,1}(\PP_{\GG_1}(E),\R)$
 \begin{equation}
    \label{eq:A3BddBoltzmann}
    \left\|  G^N (\Phi \circ \mu^N_V ) -  
      \left\langle Q(\mu^N_V,\mu^N_V), D\Phi[\mu^N_V] \right\rangle
    \right\|_{L^\infty(E^N)} \le {C_1 \over N}   \| \Phi \|_{C^{1,1}(\PP_{\GG_1}(E),\R)},
  \end{equation} 
  which is nothing but {\bf (A3)} with $k=\eta=1$ and $\eps(N)=C_1 \, N^{-1}$.

\smallskip
Take $\Phi \in C^{1,1}(\PP_{\GG_1}(E),\R)$, set $\phi = D\Phi[\mu^N_V]$ and
  compute
\begin{eqnarray*}
  G^N (\Phi \circ \mu^N_V ) 
  &=& \frac{1}{N}\sum_{1\le i<j\le N}
  \int_{\mathbb{S}^{d-1}} b(\theta_{ij})  \left[ \Phi (
    \mu^N_{V^*_{ij}}) - \Phi ( \mu^N_V)\right] \, {\rm d}\sigma \\
  &=& \frac{1}{N}\sum_{1\le i<j\le N} 
  \int_{\mathbb{S}^{d-1}} b(\theta_{ij}) \, 
  \langle \mu^N_{V^*_{ij}} -\mu^N_V, \phi \rangle\, {\rm d}\sigma  
  \qquad\quad (= I_1(V)) \\
  &+& \frac{1}{N}\sum_{1\le i<j\le N} 
  \int_{\mathbb{S}^{d-1}} {\mathcal O} \left( 
    \|\Phi \|_{C^{1,1}} \, 
    \left\|\mu^N_{V^*_{ij}} -\mu^N_V \right\|_{TV}^{2} \right)\, {\rm d}\sigma
  \quad (= I_2(V)).
\end{eqnarray*}
On the one hand, we have
\begin{eqnarray*} 
  I_1 &=& {1 \over 2N^2}\sum_{i,j= 1}^N 
  \int_{\mathbb{S}^{d-1}} b(\theta_{ij}) \, 
  \left[\phi(v^*_i) + \phi(v^*_j) -  \phi(v_i) - \phi(v_j)\right]\, {\rm d}\sigma \\
  &=&{1 \over 2}  \int_{\R^d} \int_{\R^d} \int_{\mathbb{S}^{d-1}}  
  b(\theta) \, 
  \left[\phi(v^*) + \phi(w^*) -  \phi(v) - \phi(w) \right] \, 
  \mu^N_V ({\rm d}v) \, \mu^N_V ({\rm d}w) \, {\rm d}\sigma \\
  &=&  \left\langle Q(\mu^N_V,\mu^N_V), \phi \right\rangle .
\end{eqnarray*}
On the other hand, we have 
\begin{eqnarray*} 
  I_2 (V) 
  &=& {1 \over 2N}\sum_{i,j= 1}^N 
  \int_{\mathbb{S}^{d-1}} {\mathcal O} 
  \left( \|\Phi \|_{C^{1,1}} \, 
    \left( {4 \over N} \right)^{2} \right)\, {\rm d}\sigma \\
  &\le& 8 \,  \, {\|\Phi \|_{C^{1,1}} \over N}  \, 
  \left( \sum_{i,j = 1}^N {1 \over N^2} \right)
  \le 8 \, \|b \| {\|\Phi \|_{C^{1,1}}  \over N}.
\end{eqnarray*}
Collecting these two terms we have proved that (\ref{eq:A3BddBoltzmann}) holds.

\smallskip
\noindent {\bf Proof of (A4). } Here we prove that for any $f, \, h \in
P(\R^d)$ and for any $T>0$
\begin{equation}
  \label{eq:BddBoltzstabTVproof}
  \sup_{t \in [0,T]} \Big\| S^{N\!L} _t(g) - S^{N\!L}_t(f)  
  -  {\mathcal LS}_t^{\infty}[f] (g - f) \Big\|_{TV} \le 
e^{4 \, \|\gamma \|_\infty \, T} \, \| g - f \|_{TV}^2, 
\end{equation}
where ${\mathcal LS}^{\infty}_t[f]$ is the linearization of
$S^{N\!L}_t$ at $f$. As a consequence, this implies that {\bf (A4)}
holds with $k=\eta=1$ and the previous definitions of $\PP_{\GG_1}(E)$ and
$\PP_{\GG_2}(E)$.  We denote by $f_t, g_t, h_t$ the solutions to the
following equations:
$$
\left\{
\begin{array}{l} \displaystyle
\partial_t f_t = Q(f_t,f_t), \quad f_{|t=0} = f_{\mbox{{\tiny in}}}, \vspace{0.3cm} \\ \displaystyle
\partial_t g_t = Q(g_t,g_t), \quad g_{|t=0} = g_{\mbox{{\tiny in}}}, \vspace{0.3cm} \\ \displaystyle
\partial_t h_t = 2 \tilde Q(f_t,h_t) := Q(f_t,h_t) + Q(h_t,f_t), 
\quad h_{|t=0} = h_{\mbox{{\tiny in}}}, 
\end{array}
\right.
$$
where $\tilde Q$ denotes the symmetrized form of the bilinear
collision operator. The third equation corresponds to the first-order
variation of the semigroup: ${\mathcal LS}^{\infty}_t[f](h_{\mbox{{\tiny in}}}) =h_t$
solution to this equation. 

Standard Gronwall arguments show the  existence and uniqueness of
such solutions, which moreover satisfy, uniformly on $[0,T]$
$$
\| h_t \|_{TV} \le e^{2 \,  T} \, \| h_{\mbox{{\tiny in}}} \|_{TV}, \quad
\| g_t -f_t\|_{TV} \le e^{2 \,  T} \, \| g_{\mbox{{\tiny in}}} - f_{\mbox{{\tiny in}}} \|_{TV}.
$$
Next, writing $r_t := g_t - f_t - h_t$, we find that this expression
satisfies the equation
$$
\partial_t r_t = \tilde Q ( f_t+g_t, r_t) + \tilde Q(g_t - f_t,h_t),
\qquad r_{\mbox{{\tiny in}}} = 0.
$$
Introducing $y_t := \| r_t \|_{TV}$, we have 
\begin{eqnarray*}
y'_t
&\le& {1 \over 2} \,  \|\tilde Q ( f_t+g_t, r_t) \|_{TV} + \| \tilde Q(g_t - f_t,h_t) \|_{TV} \\
&\le& \|\gamma \|_\infty \, \| f_t+g_t \|_{TV} \, \| r_t \|_{TV} + C \,  \|g_t - f_t \|_{TV} \, \|h_t \|_{TV} \\
&\le& C \, y_t + C \, e^{4t}\, \| h - f \|_{TV}^2, 
\end{eqnarray*} 
from which we deduce 
\[
\forall \, t \in [0,T], \quad y_t \le e^{4T}\, \| h - f \|_{TV}^2\,.
\]
This concludes the proof of \eqref{eq:BddBoltzstabTVproof}.

\smallskip
\noindent {\bf Proof of (A2).} Assumption {\bf (A2)-(i)} is clearly a
consequence of {\bf (A4)}. For {\bf (A2)-(ii)} we write
\begin{eqnarray*}
  \|Q(f,f) -
  Q(g,g) \|_{TV} &=& \sup_{\|\varphi \|_{L^\infty}\le 1}\int_{E}(Q(f,f)
  - Q(g,g)) \, \varphi \, {\rm d}v
  \\
  &=& \sup_{\|\varphi \|_{L^\infty}\le 1}\int_{E \times E} (f \, f_* - g
  \, g_*) \int_{\Sp^{d-1}} b \, (\varphi' - \varphi) \, {\rm d}\sigma
  \, {\rm d}v \, {\rm d}v_*
  \\
  &\le& 4 \, \|b \|_{L^1} \, \|f - g \|_{TV}, 
\end{eqnarray*}
so that the function $f \mapsto Q(f,f)$ is Lipshitz from $\PP_{\GG_1}(E)$
to $M^1(\R^d)$.

\smallskip
\noindent {\bf Proof of (A5).} 
It is known since the seminal work of Tanaka \cite{T1} that the
nonlinear Boltzmann flow associated to Maxwellian molecules is a
contraction for the quadratic MLW distance $W_2$: for any $f_{\mbox{{\tiny in}}}, \,
g_{\mbox{{\tiny in}}} \in \PP_1(\R^d)$ the solutions $f_t, \,g_t $ to the Boltzmann
equation (\ref{eq:NLBoltzEq}) satisfy
$$
\sup_{[0,T]} W_2(f_t,fg_t) \le W_2(f_{\mbox{{\tiny in}}},g_{\mbox{{\tiny in}}}). 
$$ 
That  immediately implies {\bf (A5)} in the space $P_{\GG_3}(E)$
defined above. \qed


\section{Vlasov and McKean-Vlasov equations}
\label{sec:McK}
\setcounter{equation}{0}
\setcounter{theo}{0}


\subsection{The model}
\label{sec:modelMcK}

\smallskip\smallskip\smallskip
\noindent
In this section we assume that $E = \R^m$ (where $m=d$ or $m=2d$ with
$d$ the physical space dimension, see later) and we consider an
$N$-particle system which undergoes McKean-Vlasov type stochastic
dynamics, i.e. a drift deterministic force field combined with
diffusion. We refer to the lecture notes \cite{S6,Meleard1996} and the
references therein for more details on the model, and among many
references, we highlight the recent paper \cite{BolleyGM} for recent
results and references (using the so-called ``coupling'' method). The
method we shall present here does not rely on any of these
references. The results in this section are mostly not new but compare
to the latest results of mean field limit on this equation as far as
we know. Indeed we shall make strong smoothness assumptions on the
coefficients of the evolution equation in order to avoid technical
difficulties and our goal is to a{\rm d}vocate for our new method and
show its power and ability to deal with very different models.

We assume that the $N$ particles $\ZZ^N_t = (\ZZ_{1,t}, \dots,
\ZZ_{N,t})$ satisfies the stochastic differential equation
\begin{equation}\label{eqSDDE}
{\rm d} \ZZ_{i,t} = \sigma_i(\ZZ_{i,t}) \, {\rm d} \mathscr B_{i,t} +
\mathscr T \ZZ_{i,t} \, {\rm d} t + F^{N} _i (\ZZ^N_t) \,
{\rm d}t \qquad 1 \le i \le N, 
\end{equation}
where the $\sigma(z_ i)$ are the diffusion $m \times m$-matrices, the
$\mathscr B_{i,t}$ are independent standard Wiener processes valued in $\R^m$,
$\mathscr T$ is an $m\times m$-matrix 
and the $F^{N} _i : \R^m \to \R^m$ are the force fields acting on each
particle. Because of indistinguishability we assume
$$
F^{N} _i (Z) := F^N\left(z_i, \mu^{N-1}_{\hat Z^N _i}\right)
$$
with $\hat Z^N _i:= (z_1, \dots,z_{i-1},z_{i+1}, \dots , z_N)$ and $F^N:
\R^m \times \PP(\R^m) \to \R^m$.  (Note that here and below the
latin letters ``$i$, $j$, \dots''  label the particles, whereas the
greek letters ``$\alpha$, $\beta$, \dots'' label the {\em
  coordinates}). 

We assume that $F^N$ is uniformly bounded and Lipschitz in both
variables (when endowing $\PP(\R^m)$ with a distance inherited from a
negative Sobolev norm). More precisely, we assume that for any $k >
m/2$ there exists $C_{F,k} > 0$ such that for any $z, \tilde z \in \R^m$, $f,
\tilde f \in \PP(\R^m)$
\begin{equation}\label{bdd:McVFN} 
  \forall \, N \in \N^*, \quad \left| F^N(z,f) - F^N(\tilde z,\tilde f) \right| 
  \le C_{F,k} \, \Big[ |z - \tilde z| + \| f - \tilde f \|_{H^{-k}} \Big]. 
\end{equation} 
It is also natural for the limit to exist to assume that there exists
a function $F : \R^m \times \PP(\R^m) \to \R^m$ such that $F^N\to F$, in the
sense that there is a constant $C_{F,lim} >0$ such that
\begin{equation}\label{bdd:McVFNF} 
  \forall \, N \in \N^*, \ \forall \, z
  \in \R^m, \ \forall \,  f \in \PP(\R^m), \quad 
  \left| F^N(z,f) - F(z,f) \right| \le \frac{C_{F,lim}}N.  
\end{equation}

A simple example which satisfies these assumptions is
\begin{equation}\label{eq:McKtypical}
  F^N_i \left(Z,\mu^{N-1} _{\hat Z_i}\right) = \frac{N}{N-1} \, F^N\left(z_i, \hat
    Z_i\right), \quad F^N\left(z_i, \hat
    Z_i\right) := {1 \over N}\, \sum_{j \not= i} \mathscr U (z_i - z_j)
\end{equation}
for a smooth vector field $\mathscr U: \R^m \to \R^m$,  so that 
$
F(z,f) = (\mathscr U * f)(z)$. 

Under the smoothness assumptions~\eqref{bdd:McVFN} on the $N$-particle
force fields, for any $N \ge 1$ there exists a Markov process
$(\ZZ^N_t)_{t \ge 0}$ which solves the system of stochastic
differential equations (\ref{eqSDDE}), see \cite{S6,Meleard1996}.

The time-dependent law $f^N_t$ of the process $\ZZ^N_t$ satisfies the
following linear \emph{master equation} corresponding to
(\ref{eqSDDE}), given in dual form by
\begin{equation}\label{masterMcKVN}
  \forall \, \varphi \in \DD(\R^m), \quad 
  \partial_t \left\langle f^N_t,\varphi \right\rangle 
  = \left\langle f^N_t,G^N \, \varphi \right\rangle
\end{equation}
where $G^N$ is defined by
\begin{multline*}
\forall \, Z \in \R^{mN}, \quad (G^N \varphi)(Z) = \sum_{i=1}^N A(z_i)
: \nabla^2_i \varphi + \sum_{i=1} ^N (\mathscr Tz_i) \cdot \nabla_i \varphi \\
+
\sum_{i=1}^N F^N\left(z_i, \mu^{N-1}_{\hat Z_i}\right) \cdot \nabla_i
\varphi\,.
\end{multline*}
The nonnegative diffusion matrix $A$, the
gradient $\nabla_i$ and the Hessian matrix $\nabla^2_i$ associated to
the variable $z_i = (z_{i,1}, \dots, z_{i,m}) \in \R^m$ corresponding
to the $i$-th particle are given by
$$
A = {1 \over 2} \, \sigma \, \sigma^* = \left(A_{\alpha,\beta}\right)_{1 \le
  \alpha,\beta \le m}, \quad A_{\alpha,\beta} = \sum_{\gamma=1}^d
\sigma_{\alpha,\gamma} \,
\sigma_{\beta,\gamma},
$$
and
$$
\nabla_i \varphi = \left( \partial_{z_{i,\alpha}} \varphi \right)_{1 \le
  \alpha \le m}, \quad 
\nabla^2_i \varphi = \left(\partial^2_{z_{i,\alpha} z_{i,\beta}}
  \varphi\right)_{1 \le \alpha,\beta \le m}.
$$

We also introduce the nonlinear mean field McKean-Vlasov equation on
$\PP(\R^m)$:
\begin{equation}\label{eq:nlMcKV} {\partial \over \partial t} f = Q(f_t),
  \quad f_{|t=0} = f_{\mbox{{\tiny in}}} \quad\hbox{in}\quad \PP(\R^m),
\end{equation}
with
$$
Q(f) = \sum_{\alpha,\beta =1}^m \partial^2_{\alpha,\beta}
\left(A_{\alpha,\beta}\, f\right) - \sum_{\alpha=1} ^m \partial_\alpha
[(\mathscr Tz)_\alpha f] - \sum_{\alpha=1}^m \partial_\alpha
\left(F_\alpha(z,f) \, f\right).
$$

There is an important literature on this class of nonlinear partial
differential equations. See fore example \cite{MR1970276} and
\cite{CMcCV1} where more details and more references can be found. 

In the sequel, we make the following strong structure, smoothness and
boundedness assumptions on the coefficients: 
\begin{equation}\label{eq:mkvreg1}
 (A \equiv 0, \ \kappa :=0)  \,\, \qquad\hbox{or}\qquad \,\, (A \ge \kappa \,
\hbox{Id}, \quad \kappa >0, \quad A \in W^{k,\infty}(\R^m)),
\end{equation}
as well as
\begin{equation}\label{eq:mkvreg2}
  \forall \, z \in \R^m, \ \forall \, f \in \PP(\R^m), \quad F (z,f)
  = \int_{\R^m \times \R^m} \mathscr U (z-\tilde z) \, f({\rm d}\tilde z)
\end{equation}
where $\mathscr U \in H^{2k}_6(\R^m)$ for some $k \in \N$, $k > m/2 +
3$.

When the diffusion matrix $A=0$ is zero, $m=2d$, $z=(x,v) \in
\R^{2d}$,
our assumptions cover the case of the mean-field Vlasov
equation. Indeed the classical Vlasov equation reads
\[
\partial_t f + v \cdot \nabla_x f + (\nabla_x \psi \ast \rho[f] )\cdot
\nabla_v f = 0, \quad f=f(t,x,v), \quad x, v \in \R^d, 
\]
with 
\[
\rho[f](t,x) = \int_{\R^d} f(t,x,v) \, {\rm d}v,
\]
and it falls into our structural assumptions with $z = (x,v) \in \R^d
\times \R^d$, $\mathscr U (z) = \mathscr U (x)= (0,\nabla_x
\psi( x))$ with $\nabla_x \psi$ is $H^{d+6+0}$ and 
\[
\mathscr T = \left( \begin{matrix} 0_{xx} & \mbox{Id}_{xv} \\ 0_{vx} &
    0_{vv} \end{matrix} \right).
\]
Then
$F=(F_x,F_v)$ defined by \eqref{eq:mkvreg2} is given by
\[ 
F_x(x,v) = 0, \quad 
F_v(x,v) = \nabla_x \psi \ast \rho[f]
\]
for the limiting system, and, with $X \in (\R^d)^N$ and $V \in
(\R^d)^N$, we have for the $N$-particle system $\mathscr T$ defined as
above and $F^N=(F^N_X,F^N_V)$ given by
\[
F^N _X = 0, \quad 
(F^N _V)_i = \frac1N \sum_{j \not= i} ^N \nabla_x \psi (X_i - X_j), \
i=1, \dots, N.
\]
Observe in particular that it does not allow for the Coulomb or Newton
interactions in this Vlasov setting due to the smoothness assumption on $\psi$.

\subsection{Statement of the result}
\label{sec:resultMcK}

Our main result in the section is a quantitative propagation of chaos
result for the class of equations described above. We state two
separate results respectively for the McKean-Vlasov case (possibly
non-zero diffusion matrix) and the Vlasov case (zero diffusion
matrix).

\begin{theo}[The McKean-Vlasov equation]\label{theo:McK} 
  Consider an initial distribution $f_{\mbox{{\tiny \emph{in}}}} \in
  \PP_q(\R^m)$, $q \ge 2$, the hierarchy of $N$-particle distributions
  $f^N _t = S^N _t(f_{\mbox{{\tiny \emph{in}}}} ^{\otimes N})$ following
  \eqref{masterMcKVN} and the nonlinear evolution $f_t = S^{N\!L}
  _t(f_{\mbox{{\tiny \emph{in}}}})$ following \eqref{eq:nlMcKV}. Assume that
  \eqref{eq:mkvreg1} and \eqref{eq:mkvreg2} hold.

  Then there is $k \in \N$ and a constant $C>0$ and, for any $T>0$,
  there are constants $C_{T}, \tilde C _T >0$ such that for any
  \[ 
  \varphi = \varphi_1 \otimes \dots \otimes \, \varphi_\ell \in
  \FF^{\otimes \ell}, \quad \FF := H^k_6(\R^m) \cap \mbox{{\em
      Lip}}(\R^m), \quad \| \varphi_j \|_\FF \le 1,
  \] 
  we have for $N \ge 2 \ell$:
  \begin{equation} \label{eq:cvgBddMcK} \sup_{[0,T]}\left| \left
        \langle \left( S^N_t(f_{\mbox{{\tiny \emph{in}}}} ^N) - \left(
            S^{N\!L} _t (f_{\mbox{{\tiny \emph{in}}}})
          \right)^{\otimes N} \right), \varphi \right\rangle \right|
    \le C \, {\ell^2 \over N} + C_{T} \, {\ell^2 \over N} + \tilde
    C_{T} \, \ell\, \Omega^{W_2} _N (f_{\mbox{{\tiny \emph{in}}}}).
\end{equation}

As a consequence of \eqref{eq:cvgBddMcK} and Lemma~\ref{lem:Rachev},
this implies the propagation of chaos with rate $\eps(N) \le
C(\ell,T,f_{\mbox{{\tiny \emph{in}}}}) \, N^{-{1 \over m+4}}$ for any
initial data $f_{\mbox{{\tiny \emph{in}}}} \in \PP_{m+5}(\R^m)$.
\end{theo}

Now we consider the case of the Vlasov equation. As will be clear from
the proof, when $A=0$ and $\mathscr U(0)=0$, the error $\eps(N) = 0$
vanishes in assumption {\bf (A3)}. This leads to the following
improved result.

\begin{theo}[The Vlasov equation]\label{theo:Vlasov} Suppose, in addition to the
  assumptions for Theorem~\ref{theo:McK}, that 
 $A \equiv 0$ and $\mathscr U (0) = 0$.  Then there is a constant
  $C >0$ and, for any $T>0$,  a constant $\tilde C_T>0$ such
  that for any $ \varphi \in \mbox{{\em Lip}}(\R^{\ell m}) $ and
  any $N \ge \ell$:
  \begin{multline} \label{eq:cvgMKV} \sup_{[0,T]}\left| \left \langle
        \left( S^N_t(f_{\mbox{{\tiny \emph{in}}}} ^{\otimes N}) -
          \left( S^{N\!L} _t (f_{\mbox{{\tiny \emph{in}}}})
          \right)^{\otimes N} \right), \varphi \right\rangle \right|
    \le C \, \|\nabla \varphi \|_{L^\infty(\R^{\ell m})} \, {\ell
      \over N} + \tilde C_{T} \, \Omega_N ^{W_1} (f_{\mbox{{\tiny
          \emph{in}}}})
\end{multline}
(observe the replacement of $W_2$ by $W_1$ in the last term) which in
turn implies
\[
\sup_{[0,T]} {1 \over N}W_1 \left( ( S^N_t(f_{\mbox{{\tiny \emph{in}}}} ^{\otimes N}),
  \left( S^{N\!L} _t (f_{\mbox{{\tiny \emph{in}}}}) \right)^{\otimes N} \right) \le {C \over
  N} + \frac{\tilde C_{T} \, \Omega_N ^{W_1} (f_{\mbox{{\tiny \emph{in}}}})}{N}.
\]
\end{theo}

\begin{rem}
  Note that the coupling method introduced in \cite{S6} leads to a
  rate of chaoticity of order $\OO(1/\sqrt{N})$ for the normalized
  Wasserstein distance $W_2$ between the law of $\ZZ^N_t$ and the
  tensor product $f_t^{\otimes N}$. This   is better than our
  estimate, which is limited by the estimate
  in~Lemma~\ref{lem:Rachev}. However, the coupling method is usually
  limited to the   quadratic interaction given by \eqref{eq:McKtypical}.
\end{rem}

\subsection{Proof of Theorem~\ref{theo:McK}}
\label{sec:proofMcK}

As in the proof of Theorem~\ref{theo:BddBoltz} we prove that
Theorem~\ref{theo:McK} is a consequence of
Theorem~\ref{theo:abstract}, by verifying that assumptions {\bf
  (A1)-(A2)-(A3)-(A4)-(A5)} hold. However, in the present model we
cannot, as in Section~\ref{sec:BddBoltzmann}, use the total variation
norm for the key consistency estimate {\bf (A3)} and differential
stability estimate {\bf (A4)}. The reason is that $G^N\pi^N\Phi$
involves derivatives of $Z\mapsto \Phi(\phi_Z^N)$, and hence of $Z
\mapsto \mu^N_Z$ which is not differentiable from $\R^{mN}$ to
$\PP(\R^m)$ when $\PP(\R^m)$ is endowed with the total variation norm.
We therefore make the following choice of functional spaces: $E :=
\R^m$ with
\begin{equation*}
\left\{ 
\begin{array}{lll}
\GG_1 := H^{-s_1}_{-2} (\R^m), & \FF_1 =
H^{s_1}_2(\R^m), &  s_1 > {m \over 2} + 2 \vspace{0.2cm} \\
\GG_2 := H^{-s_2}_{-6} (\R^m), & 
\FF_2 = H^{s_2}_6(\R^m), & s_2 :=  s_1 + 2,
\end{array}
\right.
\end{equation*}
and the weight $m_{\GG_1}(z) = m_{\GG_2}(z) = 1$, 
and $\FF_3 =
\mbox{Lip}(\R^m)$ and $\PP_{\GG_3}(E) := \PP_2(\R^m)$ endowed with the
quadratic MKW distance $W_2$.

\smallskip\noindent {\bf Proof of assumption (A1). }  The symmetry
assumption is a consequence of the fact that first \eqref{masterMcKVN} is
well posed for any $f^N_{\mbox{{\tiny in}}} \in \PPS(\R^{mN})$ so that $f^N_t$ is a
probability measure for any $t \ge 0$, and second the generator $G^N$
commutes with the permutations. 

\smallskip\noindent {\bf Proof of assumption (A3). } We claim that for
any $s_1 >m/2+1$ there exists a constant $C_{s_1}$ such that for all
$\Phi \in C^{2,1}(\PP_{\GG_1}(E),\R)$
\begin{equation}
  \label{eq:H1VMcK}
  \left\|  G^N (\Phi \circ \mu^N_Z ) -  
    \left\langle Q(\mu^N_Z), D\Phi[\mu^N_Z] \right\rangle
  \right\|_{L^\infty(E^N)} \le {C_{s_1} \over N}   \| \Phi \|_{C^{2,1}(\PP_{\GG_1}(E),\R)},
\end{equation} 
  which is nothing but {\bf (A3)} with $k=2$, $ \eta=1$ and $\eps(N) = C_{s_1} \, N^{-1}$.
\smallskip

\begin{proof}[Proof of \eqref{eq:H1VMcK}] 
  First, the map 
\[
\R^{m N} \to H^{-s_1}(\R^m), \quad Z \mapsto \mu^N_Z
\]
is $C^2$ with
$$
\partial_{z_{i,\alpha}} \mu^N_Z = {1 \over N} \, \partial_{{\alpha}} 
\delta_{z_i}, \qquad
\partial^2_{z_{i,\alpha},z_{i,\beta}} \mu^N_Z = {1 \over N^2} \,
\partial^2_{\alpha\beta} \delta_{z_i}.
$$
Take $\Phi \in C_b^{2,1}(\PP_{\GG_1}(E))$. Then the map 
\[
\R^{mN} \to \R,\quad Z \mapsto \Phi(\mu^N_Z)
\]
is $C_b^2$. Indeed, denoting $\phi = \phi_Z (\cdot) =
D\Phi\!\left[\mu^N_Z\right] \in (H^{-s_1}_{-2}(\R^m))' = H^{s_1}_2(\R^m) $,
we can write:
\begin{eqnarray*}
 \partial_{z_{i,\alpha}} \Phi\left(\mu^N_Z\right)
  &=& \left\langle D\Phi\!\left[\mu^N_Z\right],  {1 \over N} \, 
  \partial_\alpha \delta_{z_i} \right\rangle  
  = {1 \over N}\, \partial_{\alpha} \phi_Z (v_i)\\
  \partial^2_{z_{i,\alpha}, z_{i,\beta}} \Phi\left(\mu^N_Z\right) 
  &=& \left\langle D\Phi\!\left[\mu^N_Z\right], {1 \over N} \, 
  \partial^2_{z_{i,\alpha}, z_{i,\beta}} \delta_{z_i} \right\rangle \\
  && \quad \qquad \qquad
  + D^2\Phi\!\left[\mu^N_Z\right] \left( {1 \over N} \, 
  \partial_{z_{i,\alpha}} \delta_{z_i},  {1 \over N} \, 
  \partial_{z_{i,\beta}} \delta_{z_i} \right) \\
  &=& {1 \over N} \, \partial^2_{\alpha,\beta} \phi_Z (z_i) 
  +  {1 \over N^2} \, D^2\Phi\!\left[\mu^N_Z\right] 
  \left( \partial_{z_{i,\alpha}} \delta_{z_i},  \partial_{z_{i,\beta}} 
  \delta_{z_i} \right)
\end{eqnarray*}
and both $\partial_{z_{i,\alpha}} \delta_{z_i}$ and
$\partial^2_{z_{i,\alpha}, z_{i,\beta}} \delta_{z_i}$ belong to
$H^{-s_1}_{-2}(\R^m)$ thanks to the condition $s_1 > m/2+2$.

As a consequence, we compute
\begin{eqnarray*}
  &&\left(G^N \pi^N \Phi\right) (Z) 
  = G^N \,  \Phi(\mu^N_Z) \\
  &&\qquad = \sum_{i=1}^N A(z_i) : \nabla^2_i
  \left( \Phi(\mu^N_Z) \right) \\ 
 && \qquad \qquad \qquad + \sum_{i=1} ^N (\mathscr Tz_i) \cdot \nabla_i
  \left(\Phi(\mu^N_Z)\right)  
  +  \sum_{i=1}^N F^N\left(z_i, \mu^{N-1}_{\hat Z_i}\right) \cdot 
  \nabla_i \left(\Phi(\mu^N_Z)\right) \\
  &&\qquad =: I_1(Z) + I_2(Z)
\end{eqnarray*}
with 
\begin{multline*}
I_1(Z) := \frac{1}{N} \, \sum_{i=1}^N \sum_{\alpha,\beta=1}^m 
  A_{\alpha,\beta}(z_i) \, 
  \partial^2_{\alpha,\beta} \phi_Z(z_i) \\
  + \sum_{i=1} ^N \sum_{\alpha, \beta =1} ^m (\mathscr T_{\alpha \beta}
  z_i) \, \partial_\beta \phi_Z(z_i) + \frac1N \, \sum_{i=1}^N
  \sum_{\alpha=1}^m F_\alpha \left(z_i, \mu^N _Z\right)
  \, \partial_{\alpha} \phi_Z (z_i)
\end{multline*}
and 
\begin{multline*}
I_2(Z) := \frac{1}{N^2} \, \sum_{i=1}^N \sum_{\alpha,\beta=1}^m
A_{\alpha,\beta}(z_i) \, 
D^2\Phi\!\left[\mu^N_Z\right] \left( \partial_{z_{i,\alpha}} \delta_{z_i},  
  \partial_{z_{i,\beta}} \delta_{z_i} \right) \\
+ \frac1N \, \sum_{i=1}^N \sum_{\alpha=1}^m \left[F^N_\alpha\left(z_i,
    \mu^{N-1}_{\hat Z_i}\right) - F_\alpha\left(z_i,
    \mu^{N}_Z\right) \right] \,D\Phi\!\left[\mu^N_Z\right]
\left( \partial_{z_{i,\alpha}} \delta_{z_i}\right).
\end{multline*}

On the one hand, using that 
\[
\left\| \mu^{N-1}_{\hat Z_i} - \mu^{N}_Z\right\|_{H^{-s_1}_{-2}(\R^m)}
\le \frac{2}N \, \sup_{z_i} \left\| \delta_{z_i}
\right\|_{H^{-s_1}(\R^m)} \le \frac{C}N
\]
as well as (\ref{bdd:McVFN}) and (\ref{bdd:McVFNF}), we deduce that
\begin{eqnarray*}
  |I_2(Z)| &\le& N \, \frac{m^2}{N^2} \, \|A \|_\infty  \, 
  \|D^2 \Phi \|_\infty \, \|\partial_1 \delta \|^2_{H^{-s_1}_{-2}(\R^m)}  \\
  && + N \, m \, \left( \frac{C_{F}}{N} \right) \,  
  {1 \over N} \, \|D \Phi \|_\infty \, \|\partial_1 \delta \|_{H^{-s_1}_{-2}(\R^m)} 
  \le \frac{C_\Phi}N.
\end{eqnarray*}

On the other hand, we recognize 
\begin{eqnarray*}
  I_1 (Z) 
  &=& \left\langle \mu^{N}_Z \, ,  \,   
    \sum_{\alpha,\beta=1}^m A_{\alpha,\beta} \, 
    \partial^2_{\alpha,\beta} \phi_Z \right\rangle \\
  && \qquad \qquad 
  +  \left\langle \mu^{N}_Z \, ,  \,   
    \sum_{\alpha=1}^m (\mathscr T \cdot)_{\alpha} \, 
    \partial_{\alpha} \phi_Z \right\rangle + \left\langle \mu^{N}_Z \, ,  \,  
    \sum_{\alpha=1}^m  F_\alpha\left(\cdot, \mu^{N}_Z\right) \, 
    \partial_{\alpha} \phi_Z   \right\rangle \\
  &=& \left\langle Q(\mu^{N}_Z), \phi_Z \right\rangle 
  = \left\langle Q(\mu^{N}_Z), D\Phi(\mu^{N}_Z) \right\rangle 
  = \left(\pi^N G^\infty \Phi\right)(Z),
\end{eqnarray*}
thanks to the calculation of the limit dual generator made in
Subsection~\ref{sec:calculus-gen}. 
\end{proof}

\smallskip\noindent
{\bf Proof of assumption (A4). }
We need here to perform a {\em second}-order expansion of the limit
semigroup. 

We consider
  \begin{itemize}
  \item for any two given initial data $f_{\mbox{{\tiny in}}}, g_{\mbox{{\tiny in}}} \in \PP(\R^m)$
    the corresponding solutions $f_t$ and $g_t$ to the nonlinear
    McKean-Vlasov (or Vlasov) equation (\ref{eq:nlMcKV}),
  \item for any given initial data $h_{\mbox{{\tiny in}}} \in \PP(\R^m)$ the solution
    $h_t$ to the following equation, which is the linearization around
    $f_t$:
  \begin{equation}\label{eq:nlMcKVLin}
  \partial_t h = \nabla^2  : (A \, h) - \nabla \cdot ((\mathscr Tz) \, h)
  - \nabla \cdot \left[ h \, (\mathscr U*f) + f \, (\mathscr U *h) \right],
  \quad h_{t=0} = h_{\mbox{{\tiny in}}},
  \end{equation}
  \item $r_t$ the solution to the following second variation 
  equation around $f_t$ 
  \begin{equation}\label{eq:nlMcKVQuad}
  \left\{ 
    \begin{array}{l} \displaystyle
      \partial_t r = \nabla^2 : (A \, r) - \nabla \cdot ((\mathscr Tz) \, r)
      - \nabla \cdot \left[ r \, (\mathscr U *f) +  f \, (\mathscr U *r)\right] \vspace{0.2cm} \\
      \displaystyle 
      \qquad \qquad \qquad \qquad \qquad \qquad - \frac{1}{2} \nabla
      \cdot 
      \left[ \tilde h \, (\mathscr U *h) \right] - \frac{1}{2} \nabla
      \cdot 
      \left[ h \, (\mathscr U * \tilde h) \right], \vspace{0.12cm} \\
      \displaystyle       r_{|t=0} = r_{\mbox{{\tiny in}}} = 0
\end{array}
\right.
\end{equation}
for two solutions $h, \tilde h$ of the first variation equation. 
  \end{itemize}

Then we shall prove the following a priori estimates.
\begin{lem}\label{lem:VMcK3} 
  For any $s_1 \in \N$, $s_1> m/2 +1$, $\ell \in \{1,2,3\}$ and for any
  $T > 0$, there exists $C_T$ such that
  \begin{eqnarray}
  \label{lem:VMcK3-1}
  && \sup_{[0,T]} \|g_t - f_t \|_{H^{-s_1}_{-\ell}(\R^m)} 
  \le C_T \, \|g_{\mbox{{\tiny \emph{in}}}} - f_{\mbox{{\tiny \emph{in}}}} \|_{H^{-s_1}_{-\ell}(\R^m)},
   \\ \label{lem:VMcK3-1bis}
  && \sup_{[0,T]} \| h_t \|_{H^{-s_1}_{-\ell}(\R^m)}
  \le C_T \, \| h_{\mbox{{\tiny \emph{in}}}} \|_{H^{-s_1}_{-\ell}(\R^m)},
   \\ \label{lem:VMcK3-1ter}
   && \sup_{[0,T]} \|r_t \|_{H^{-(s_1+1)}_{-4}(\R^m)}
   \le C_T \, \|h_{\mbox{{\tiny \emph{in}}}} \|_{H^{-s_1}_{-2}(\R^m)}
   \, 
   \| \tilde h_{\mbox{{\tiny \emph{in}}}} \|_{H^{-s_1}_{-2}(\R^m)}, 
  \end{eqnarray} 
and when $\tilde h_{\mbox{{\tiny \emph{in}}}} = h_{\mbox{{\tiny
      \emph{in}}}} =
g_{\mbox{{\tiny \emph{in}}}} - f_{\mbox{{\tiny \emph{in}}}}$ we have 
\begin{eqnarray}
  \label{lem:VMcK3-2}
  && \sup_{[0,T]} \|g_t - f_t - h_t \|_{H^{-(s_1+1)}_{-4}(\R^m)} 
  \le C_T \, \|g_{\mbox{{\tiny \emph{in}}}} - f_{\mbox{{\tiny \emph{in}}}} \|^2_{H^{-s_1}_{-2}(\R^m)}, \\
  \label{lem:VMcK3-3}
  && \sup_{[0,T]} \|g_t - f_t - h_t - r_t \|_{H^{-(s_1+2)}_{-6} (\R^m)} 
  \le C_T \, \|g_{\mbox{{\tiny \emph{in}}}} - f_{\mbox{{\tiny \emph{in}}}} \|^3_{H^{-s_1}_{-2}(\R^m)}.
  \end{eqnarray} 

  This shows that the nonlinear semigroup $S^{N\!L}_t$ associated to the
  nonlinear McKean-Vlasov equation (\ref{eq:nlMcKV}) is
  $C_b^{2,1}(\PP_{\GG_1}(E),\PP_{\GG_2}(E))$.
 \end{lem}

\begin{proof}[Proof of Lemma~\ref{lem:VMcK3}] 
The proof is carried out in several steps. 

\noindent 
{\sl Step 1. } 
We will several times consider  the equation
\begin{equation}\label{eq:McVlinearizedeq}
  \partial_t \zeta_t = \nabla^2 : (A \, \zeta_t)  - \nabla \cdot
  ((\mathscr Tz)
  \zeta_t) -  \nabla \cdot ( u_1 \, \zeta_t + u_2 \, (\mathscr U \ast \zeta_t) )
\end{equation}
with given initial data $\zeta_{\mbox{{\tiny in}}}$ and with an
$\R^m$-valued function $u_1$ and an $\R$-valued measure $u_2$ to be
specified (chosen in order to ``match'' equations (\ref{eq:nlMcKV}),
(\ref{eq:nlMcKVLin}) and (\ref{eq:nlMcKVQuad})).  We claim that for
any $k \in \N$, $k > m/2 + 1$, $\ell \in \{1,2,3\}$ and any $T >0$,
\begin{equation}\label{estim:McVLinearH-k} 
  \forall \, t \in
  [0,T], \quad \left\| \zeta_t   \right\|_{H^{-k}_{-\ell}(\R^m)} \le
  \left\|\zeta_{\mbox{{\tiny in}}} \right\|_{H^{-k}_{-\ell}(\R^m)} \, e^{C_k(\mathscr U ,u_1,u_2) \, T} 
\end{equation}
with
\[
C_k(\mathscr U ,u_1,u_2) := C(k) \, \sup_{t \in [0,T]} \Bigg[ \left\|
  u_1 \right\|_{W^{k,\infty}(\R^m)} + \|\mathscr U \|_{H^{k}_\ell
  (\R^m)} \, \| u_2 \|_{TV(\R^m)} \Bigg].
\]

We argue by duality and we consider a smooth solution $\theta$ to the
following linear equation (which is the dual equation of
\eqref{eq:McVlinearizedeq})
\begin{equation}\label{eq:dualSG}
  \left\{ 
    \begin{array}{l} \displaystyle
      \partial_t \theta  =  L^*_1 \theta + L^*_2 \theta, \vspace{0.2cm} \\ \displaystyle
      L^*_1 \theta := A : \nabla^2 \theta + (\mathscr T z) \cdot \nabla
      \theta, \vspace{0.2cm} \\ \displaystyle
      L^*_2 \theta :=  u_1  \cdot \nabla \theta + \check{\mathscr U}  \ast \left(
        u_2  \, \nabla \theta \right),
    \end{array}
\right.  
\end{equation}
with $\check{\mathscr U}(x) := \mathscr U (-x)$. 

For a given multi-index $\nu \in \N^m$ with $|\nu| = k' \le k$, we
compute
\begin{multline*}
  {{\rm d} \over {\rm d}t} \int_{\R^m} \left|\partial^{\nu}
    \theta\right|^{2\ell} \, \langle z \rangle^{2\ell} \, {\rm d}z \\ =
  \int_{\R^m} (\partial^{\nu} L^*_1 \theta ) \, \partial^{\nu} \theta
  \, \langle z \rangle^{2\ell} \, {\rm d}z + \int_{\R^m}
  (\partial^{\nu} L^*_2 \theta ) \, \partial^{\nu} \theta \, \langle z
  \rangle^{2\ell} \, {\rm d}z =: \LL_1 + \LL_2.
\end{multline*}
%
By integrations by parts, we get  
\begin{equation*}
  \LL_1 \le - \kappa \, \int_{\R^m} \left| \nabla \partial^\nu \theta
  \right|^2 \,  \langle z \rangle^{2\ell} \, {\rm d}z + C \, \left( \| \mathscr T\|_\infty + \| A
    \|_{W^{k'+2,\infty}(\R^m)} \right) \, \left\| \theta \right\|^2_{H^{k'}_\ell (\R^m)} 
\end{equation*}
and (using Sobolev embedding inequalities for the last term)
\begin{equation*}
  \LL_2 \le C \left( \| u_1\|_{W^{k',\infty}(\R^m)} + \| u_2
    \|_{TV(\R^m)} \, \| \mathscr U  \|_{H^{k'}_\ell(\R^m)} \right)
  \| \theta_t \|^2_{H_\ell^{k}(\R^m)}
\end{equation*}
which shows that 
\[ 
\forall \, t \in [0,T], \quad \left\| \theta_t \right\|_{H^k_\ell(\R^m)} \le
\left\|\theta_{\mbox{{\tiny in}}}\right\|_{H^k_\ell(\R^m)} \, e^{C_k(\mathscr U ,u_1,u_2) \, T}.
\]

Denoting by $U_t$ the linear semigroup associated to
\eqref{eq:McVlinearizedeq}, the associated dual semigroup $U^*_t$ is
generated by \eqref{eq:dualSG}.  As a consequence, for any $\theta_{\mbox{{\tiny in}}} \in
H^k(\R^m)$, we have
\begin{multline*} 
  \langle \zeta_t,\theta_{\mbox{{\tiny in}}} \rangle = \left\langle
    \zeta_{\mbox{{\tiny in}}}, U^*_t \theta_{\mbox{{\tiny in}}}
  \right\rangle \le \| \zeta_{\mbox{{\tiny in}}}
  \|_{H^{-k}_{-\ell}(\R^m)} \, \| U^*_t \theta_{\mbox{{\tiny in}}}
  \|_{H^{k}_\ell(\R^m)}
  \\
  \le e^{C_k(\mathscr U ,u_1,u_2) \, T} \, \| \zeta_{\mbox{{\tiny
        in}}} \|_{H^{-k}_{-\ell}(\R^m)} \, \|
  \theta_{\mbox{{\tiny in}}} \|_{H^{k}_\ell(\R^m)},
  \end{multline*} 
which concludes the proof of the claim \eqref{estim:McVLinearH-k}. 

\smallskip\noindent {\sl Step 2.  Proof of (\ref{lem:VMcK3-1}). }  The
equation satisfied by the difference $\mathsf d_t = g_t - f_t$ is
\begin{equation} \label{eq:defUt}
  \left\{ 
    \begin{array}{l} 
      \partial_t \mathsf d_t =  \nabla^2 : (A \, \mathsf d_t) - \nabla
      \cdot ( (\mathscr T z) \mathsf d_t)
      - \nabla \cdot \left( \mathsf d_t \, (\mathscr U  \ast f) +
    g \, (\mathscr U  \ast \mathsf d_t) \right), \vspace{0.2cm} \\ 
  \mathsf d_{|t=0} = \mathsf d_{\mbox{{\tiny in}}} = g_{\mbox{{\tiny in}}} - f_{\mbox{{\tiny in}}},
\end{array}
\right.
\end{equation}
which fits in the form \eqref{eq:McVlinearizedeq} with $u_1 := \mathscr U  *f$ and $u_2 = g$.
Now, since 
\[
\left\| \nabla^k (\mathscr U \ast f) \right\|_{L^\infty(\R^m)} =
\left\| ( \nabla^k \mathscr U ) \ast f \right\|_{L^\infty(\R^m)} \le
\left\| \nabla^k \mathscr U \right\|_{L^\infty(\R^m)}
\]
we conclude that 
\[
C_k(\mathscr U , \mathscr U \ast f ,g) \le C \, \|  \mathscr U
  \|_{H^k_3 \cap W^{k,\infty} (\R^m)}
\]
and that \eqref{lem:VMcK3-1} holds.  Proceeding in the same way for
the function $h$ we end up with
\begin{equation}\label{eq:VMch&delta}
  \sup_{[0,T]} \| h_t \|_{H^{-s_1}_{-\ell}(\R^m)} \le C_T \, \| h_{\mbox{{\tiny in}}} \|_{H^{-s_1}_{-\ell}(\R^m)},
\end{equation} 
for any $s_1 \in \N$, $s_1 > m/2 + 1$. 

\smallskip\noindent {\sl Step 3. Inequalities for products and
  convolutions in Sobolev spaces.} We define the weighted
$L^\infty$-based Sobolev spaces as usual: 
\begin{equation*}
  \| f \|_{L^\infty _{-\ell}(\R^m)} := \left\| \langle \cdot
  \rangle^{-\ell} f \right\|_{L^\infty(\R^m)},
\end{equation*}
\begin{equation*}
\| f \|_{W^{k,\infty}_{-\ell}(\R^m)} := \sum_{0 \le k' \le k} \left\| \langle \cdot
  \rangle^{-\ell} \partial^{k'} f \right\|_{L^\infty(\R^m)}.
\end{equation*}

We have the three following inequalities on functions $S, \psi$ in the
appropriate spaces:
\begin{equation*}
  \|S\ast \mathscr U  \|_{L^\infty_{-\ell} (\R^m)} \le \|\mathscr U  \|_{H^k_\ell(\R^m)} \, \|S \|_{H^{-k}_{-\ell}(\R^m)}
\end{equation*}
and more generally 
\begin{equation*}
  \|S \ast \mathscr U \|_{W^{k,\infty}_{-\ell} (\R^m)} \le \|\mathscr U  \|_{H^{2k}_\ell(\R^m)} \, \|S \|_{H^{-k}_{-\ell} (\R^m)}
\end{equation*}
and finally 
\begin{equation*}
\|S \, \psi \|_{H^{-k}_{-\ell}(\R^m)} \le C_{k,\ell} \, \|S\|_{H^{-k}_{-\ell}(\R^m)} \, \|\psi \|_{W^{k,\infty}(\R^m)}
\end{equation*}
for any $k,\ell \in \N$, and some constant $C_{k,\ell} >0$. The proofs
are elementary and we omit them for the sake of conciseness.


\smallskip\noindent {\sl Step 4.  Proof of (\ref{lem:VMcK3-2}). }
Let  $\omega_t := g_t - f_t - h_t = \mathsf d_t - h_t$, which satisfies the equation 
\begin{equation}\label{eq:defphit}
\partial_t \omega =  L \, \omega + \Sigma, \quad \omega_{|t=0} =
\omega_{\mbox{{\tiny in}}} = 0,
\end{equation}
with 
\begin{equation*}
\left\{ 
\begin{array}{l}
L \, \omega := \nabla^2 : (A \, \omega ) - \nabla \cdot ( (\mathscr T
z) \omega_t ) - \nabla \cdot \left( \omega
  \, (\mathscr U\ast f) + f \, (\mathscr U \ast\omega) \right),
\vspace{0.3cm} \\
\Sigma_t = \nabla \cdot \left( \mathsf d_t \, (\mathscr U \ast \mathsf d_t)
\right).
\end{array}
\right.
\end{equation*}

Denoting by $\Theta_{s,t}w$ the unique solution of the linear,
non-autonomous equation
\[
  \partial_t w_t = L w_t,\qquad w_s=w,
\]
%
the Duhamel
formula for equation \eqref{eq:defphit} yields
$$
\omega_t = \int_0^t \Theta_{s,t} \, \Sigma_s \, {\rm d}s.
$$
Therefore we obtain, using \eqref{estim:McVLinearH-k} and the
estimates established in the Step~3, that for any $t \in [0,T]$
\begin{eqnarray*} 
  \| \omega_t \|_{H^{-k}_{-4}(\R^m)}
  &\le& 
  C_T \, \int_0^t \left\| \nabla \left( \mathsf d_s \, (\mathscr U
      \ast \mathsf d_s)  
    \right) \right\|_{H^{-k}_{-4}(\R^m)} \, {\rm d}s 
  \\
  & \le & 
  C_{T,k} \, \int_0^t \| \nabla \mathsf d_s \|_{H^{-k}_{-2}(\R^m)} \, 
  \| \mathscr U \ast \mathsf d_s  \|_{W^{k,\infty}_{-2}(\R^m)} \, {\rm d}s \\
  && \qquad \qquad \quad 
  +  C_{T,k} \, \int_0^t  \|  \mathsf d_s \|_{H^{-k}_{-2}(\R^m)} \| \mathscr U \ast
  (\nabla \mathsf d_s) \|_{W^{k,\infty}_{-2}(\R^m)} \, {\rm d}s
  \\
  & \le & C_{T,k} \, \|\mathscr U \|_{H_2^{2k}(\R^m)} \, \int_0^t \left\| \mathsf d_s
  \right\|_{H_{-2}^{-(k-1)}(\R^m)} \,  \left\| \mathsf d_s \right\|_{H^{-k}_{-2}(\R^m)} \, {\rm d}s,  
\end{eqnarray*} 
which together with \eqref{lem:VMcK3-1} for the control of the norms
of $\mathsf d_t$ implies (\ref{lem:VMcK3-2}).

\smallskip\noindent {\sl Step 5.  Proof of \eqref{lem:VMcK3-1ter} and
  (\ref{lem:VMcK3-3}).}  The second variation $r$ satisfies the
equation
\[
\partial_t r =  L \, r + R_t, \quad r_{|t=0} = 0,
\]
with $L$ as above and 
$$
R_t := \frac12 \nabla \cdot \left( \tilde h_t \, (\mathscr U \ast h_t)
\right) + \frac12 \nabla \cdot \left( h_t \, (\mathscr U \ast \tilde h_t) \right).
$$
We proceed as in Step 4, taking advantage of the bound
\eqref{eq:VMch&delta}, and we obtain
\begin{equation}\label{eq:controlrMKV}
\sup_{[0,T]} \| r_t \|_{H^{-k}_{-4}(\R^m)} \le C_T \, \| h_{\mbox{{\tiny in}}}
\|_{H^{-(k-1)}_{-2}(\R^m)} \, \| \tilde h_{\mbox{{\tiny in}}} \|_{H^{-(k-1)}_{-2}(\R^m)},
\end{equation}
which is nothing but \eqref{lem:VMcK3-1ter}. 

Finally we introduce $\psi_t := g_t - f_t - h_t -r_t = \mathsf d_t -
h_t - r_t = \omega_t - r_t$, with the initial data $\tilde h_{\mbox{{\tiny in}}} =
h_{\mbox{{\tiny in}}} = g_{\mbox{{\tiny in}}} - f_{\mbox{{\tiny in}}}$. It satisfies the equation
\begin{equation}\label{eq:defphitA}
\left\{ 
\begin{array}{l} 
\partial_t \psi=  L \, \psi + \Psi_t , \quad \psi_{|t=0} = 0,
\vspace{0.3cm} \\
\Psi_t :=  \nabla \left( \omega_t \, (\mathscr U \ast \mathsf d_t) 
  + h_t \, (\mathscr U \ast \omega_t) \right)
\end{array}
\right.
\end{equation}
Therefore, we deduce
\begin{multline*}
\forall \, t \in [0,T], \quad \| \psi_t \|_{H^{-k}_{-6}(\R^m)} 
\le \left\| \int_0^t \Theta_{s,t} \, \Psi_s \, {\rm d}s \right \|_{H^{-k}_{-6} (\R^m)}  \\
\le C_T  \, 
\int_0^t \left( \|  h_s \|_{H^{-(k-1)}_{-2}(\R^m)} + \| \mathsf d_s
  \|_{H^{-(k-1)}_{-2}(\R^m)} \right) \, \|  \omega_s \|_{H^{-(k-1)}_{-4} (\R^m)} \, {\rm d}s,
\end{multline*}
which together with \eqref{lem:VMcK3-1}-\eqref{lem:VMcK3-1bis} and (\ref{lem:VMcK3-2}) implies
(\ref{lem:VMcK3-3}).
\end{proof}

\smallskip\noindent {\bf Proof of (A2).}  The first property {\bf
  (A2)-(i)} is a consequence of \eqref{lem:VMcK3-1} in
Lemma~\ref{lem:VMcK3}.  
we have \begin{equation}\label{estim:McKVLipQHk}
  \left\{ \begin{array}{l} \displaystyle
      \|Q(f_1) \|_{H^{-k}_{-2}(\R^m)} \le C_{\mathscr U,1} \vspace{0.2cm} \\
      \displaystyle \|Q(f_2) \|_{H^{-2}_{-2}(\R^m)} \le C_{\mathscr
        U,1} \vspace{0.2cm} \\ \displaystyle \|Q(f_2) - Q(f_1)
      \|_{H^{-k}_{-2}(\R^m)} \le C_{\mathscr U,2} \, \|f_2 - f_1
      \|_{H^{-k}_{-2}(\R^m)} ^{1/5}.  \end{array}
  \right.  \end{equation} We write $$ Q(f) = Q_1(f) + Q_2(f) $$
with $$ Q_1(f) = \nabla^2 :(A \, f) - \nabla \cdot ( (\mathscr Tz) f),
\quad Q_2(f) = - \nabla \cdot ( (\mathscr U*f) \, f).  $$ The linear
term $Q_1$ satisfies the first part of \eqref{estim:McKVLipQHk} (first
equation) by direct inspection combined with the use of Sobolev
embeddings.  In order to see that $Q_1$ also satisfies the second part
of \eqref{estim:McKVLipQHk} (on the difference), we
write \begin{eqnarray*} 
  \|Q_1(f_2) - Q_1(f_1) \|_{H^{-k}_{-2}(\R^m)}
  &\le& \| A \|_{W^{k-2,\infty}(\R^m)} \, \| f_2 - f_1
  \|_{H^{-(k-2)}_{-2}(\R^m)}
  \\
  &&+ \| \mathscr T \|_{\infty} \, \| f_2 - f_1
  \|_{H^{-(k-1)}_{-1}(\R^m)}, 
\end{eqnarray*} 
and we conclude by using interpolation and Sobolev embeddings
(noticing that $k - 2 - 1/2 > m/2$ allows for the Sobolev embedding in
the first term).

Concerning the quadratic term $Q_2$, using
estimates proved in Step 3 of the proof of {\bf (A4)}, on the one hand we get
\begin{multline*}
  \|Q_2(f) \|_{H^{-k}_{-2}(\R^m)} \le   \|Q_2(f) \|_{H^{-k}(\R^m)} \le  \| \mathscr U \ast f   \|_{W^{k-1,\infty}(\R^m)} \,
  \|f \|_{H^{-(k-1)}(\R^m)}
  \\
  \le C_k \, \| \mathscr U \|_{H^{2(k-1)}(\R^m)} \, \|f \|_{H^{-(k-1)}(\R^m)}^2 \le C_k \,
  \| \mathscr U \|_{H^{2(k-1)}(\R^m)} 
\end{multline*}
where we have used $\PP(\R^m) \subset H^{-(k-1)/2}(\R^m)$ with
continuous embedding. On the other hand, we have  with $\mathsf d := f_2 - f_1$
\begin{multline*}
 \|Q_2(f_2) - Q_2(f_1) \|_{H^{-k}_{-2}(\R^m)} 
 \le 
  \|(\mathscr U \ast \mathsf d) \, \langle \cdot \rangle^{-2}
  \|_{W^{k-1,\infty}(\R^m)} \, \|f_2 \|_{H^{-(k-1)}(\R^m)} \\ 
  + \|\mathscr U \ast f_1 \|_{W^{k-1,\infty}(\R^m)} \, \| \mathsf d \|_{H^{-(k-1)}_{-2}(\R^m)}.
\end{multline*} 
In order to estimate the first term in the above inequality, we remark that
\begin{eqnarray*}
  \langle z \rangle^{-2} \, |(\partial^\alpha \mathscr U \ast \mathsf d)(z)| 
  &\le&
  \langle z \rangle^{-2} \, \| \partial^\alpha \mathscr U (z-\cdot) \, \langle \cdot \rangle^2
  \|_{H^{k}(\R^m)} \, 
  \| \mathsf d \, \langle \cdot \rangle^{-2} \|_{H^{-k}(\R^m)}
  \\
  &\le& C \, \|  \partial^\alpha \mathscr U   \|_{H^{k}_2(\R^m)} \, \| \mathsf
  d \|_{H^{-k}_{-2}(\R^m)}
\end{eqnarray*}
uniformly for any $z \in \R^m$. All together, we have for the quadratic term 
$$
 \|Q_2(f_2) - Q_2(f_1) \|_{H^{-k}_{-2}(\R^m)} \le C'_{\mathscr U,2} \,  \| f_2-f_1   \|_{H^{-(k-1)}_{-2}(\R^m)}, 
$$
and we  conclude the proof of {\bf (A2)-(ii)} by using interpolation and Sobolev embeddings again.
%

\smallskip\noindent {\bf Proof of (A5).}  We use the well known
following estimate (see \cite{S6}): for any $q \ge 1$, $f_{\mbox{{\tiny in}}}, \,
g_{\mbox{{\tiny in}}} \in \PP_q(\R^d)$ and $T > 0$ there exists $C_T$ such that
$$
\sup_{t \ge 0}W_q( S^{N\!L}_t (f_{\mbox{{\tiny in}}}), S^{N\!L}_t (g_{\mbox{{\tiny in}}})) \le C_T \, W_q (f_{\mbox{{\tiny in}}},g_{\mbox{{\tiny in}}}),
$$
that we use with $q=2$. Alternatively, estimate \eqref{lem:VMcK3-1}
precisely says that assumption {\bf (A5)} holds in $\PP_{\GG_1}(E)$.


\section{Inelastic collisions with thermal bath}
\label{sec:thermostat}
\setcounter{equation}{0}
\setcounter{theo}{0}


\subsection{The model}
\label{sec:modelThermostat}

In this section we assume that $E= \R^d$, $d \ge 1$, and we are
interested in the following \emph{Boltzmann equation for diffusively
  excited granular media} on the distribution $f(t,v) \ge 0$, $v \in
\R^d$ of particles:
\begin{equation}\label{eq:nlMix} 
{\partial f_t\over \partial t} = Q(f_t), 
\quad f_{|t=0} = f_{\mbox{{\tiny in}}} \quad\hbox{in}\quad \PP(\R^d),
\end{equation}
with
$$
Q(f) =  Q_\alpha(f,f) + \nu \, \Delta \, f
$$
for some $\nu >0$, and where the quadratic Boltzmann collision kernel
$Q_\alpha$ is defined by the following dual formulation
\begin{equation}\label{defQalpha}
\langle Q_\alpha(f,f), \varphi \rangle := \int_{\R^{2d} \times \mathbb{S}^{d-1}} 
b(\cos \theta) \, \left( \phi (w^*_2) -  \phi (w_2) \right)  \, {\rm d}\sigma \, f({\rm d}w_1) \, f({\rm d}w_2)
\end{equation}
for any $\varphi \in C_0(\R^d)$, $f \in \PP(\R^d)$, and with $\cos
\theta = \sigma \cdot (w_2-w_1)/|w_2-w_1|$ and  similarly as in equations~\eqref{vprimvprim*} and \eqref{eq:w1*w2*}
$$
w^*_2 = {w_1+ w_2 \over 2}+ {u^*\over 2},
$$
but with 
$$
 u^* = \left( {1- \alpha
    \over 2} \right) \, (w_1-w_2) + \left( {1+\alpha \over 2} \right) \,
|w_2 - w_1| \, \sigma,
$$
for some $\alpha \in (0,1)$. (Note that the case $\alpha = 1$ and $\nu
=0$ would correspond to the elastic Boltzmann kernel considered in
Section~\ref{sec:BddBoltzmann}.) This corresponds to a situation where
particles lose energy when they collide.  We refer to \cite{BCG00,CCG}
for a physical motivation to these equations. The mathematical theory
is treated in e.g. \cite{BCG00,BisiCT,BolleyC}, where, for example,
it is proven that this equation generates a nonlinear semigroup
$S^{N\!L}_t f_{\mbox{{\tiny in}}} := f_t$ for any $f_{\mbox{{\tiny
      in}}} \in \PP_q(\R^d)$, $q \ge 2$. Notice that unlike the
classical Boltzmann equation the kinetic energy is not conserved. For
the sake of simplicity we make the normalization assumptions
\[
\|b \|_{L^1(\mathbb S^{d-1})} = \int_{\mathbb{S}^{d-1}} b(\sigma_1) \,
{\rm d}\sigma = 1, \qquad \nu = 1.
\]
One of these quantities, say the first,  can be set to one just by a
rescaling of time but the two cannot be changed
independently. However, the result would be the same for any value of $\nu$.
Note that due to the normalization $\|b \|_{L^1(\mathbb{S}^{d-1})} =
1$ and the fact that $f \in \PP(\R^d)$, the bilinear operator
$Q_\alpha$ splits into a quadratic part and a linear part
 $$
 Q(f) = Q^+_\alpha(f,f) - f +  \Delta f,
 $$
 where $Q^+$ is defined through the positive part of the expression
 \eqref{defQalpha}.
 
 \smallskip We now want to introduce a $N$-particle system associated
 to the above Boltzmann equation for diffusively excited granular
 media by mimicking the Kac's construction.  We consider the
 velocities process $(\VV^N_t)$ with values in $\R^{dN}$, of mixed
 jump and diffusion nature, defined through the
 stochastic differential equations
\begin{multline*}
\VV^N_{t} = \VV^N_{0} + \\ \int_0^t \int_{\Sp^{d-1}} \sum_{i,j=1} ^N
\Gamma_{i,j,\sigma}(\VV^N_{s^-}) \, {\bf 1}_{z < b ( \sigma \cdot \hat
  u_{i,j}(\VV^N_{s^-}) ) } \, \NN^{N}(ds,{\rm d}\sigma,i,j,{\rm d}z) + \sqrt{2} \, \mathscr B^N_t.
\end{multline*}
Here $\mathscr B^N_t$ is a $\R^{dN}$ valued standard Brownian motions,
$\NN^{N}({\rm d}s,{\rm d}\sigma,i,j,{\rm d}z)$ is a
Poisson measure on $[0,\infty) \times \Sp^{d-1} \times \{1, \dots , N
\}^2 \times \R_+$ with intensity
$$
{\rm d}s \, {\rm d}\sigma \, {1 \over N} \sum_{i',j'=1}^N {\bf 1}_{i'
  \not= j'} \delta_{(i',j')}(i,j) \, {\rm d}u
$$
independent of $\mathscr B^N_t$, and the two functions $
\Gamma_{i,j,\sigma} : \R^{dN} \to \R^{dN}$ and $\hat u_{i,j} : \R^{dN}
\to \Sp^{d-1}$ are a.e. defined through the following expressions: for
any $V = (v_1, \dots, v_N) \in \R^{dN}$ we set
$$
\hat u_{i,j} (V) := \frac{u_{ij}}{|u_{ij}|}, \quad u_{ij} := v_i-v_j 
$$
and 
$$
 \Gamma_{i,j,\sigma} (V) := V^*_{ij} - V, 
$$
where 
\[
V^*_{ij} = (v_1, \dots, v_{i-1}, v^*_i, v_{i+1}, \dots, v_{j-1},
v^*_j, v_{j+1}, \dots, v_N)
\]
and, as in equation~\eqref{vprimvprim*},
\begin{equation}\label{eq:v*iv*j}
  v^*_i = {w_{ij} \over 2} +  {u^*_{ij} \over 2}, \quad 
v^*_j= {w_{ij} \over 2} - {u^*_{ij} \over 2}, 
\end{equation}
but here with
$$
w_{ij} = v_i+v_j, \quad u^*_{ij} = \left( {1- \alpha \over 2} \right)
\, u_{ij} + \left( {1+\alpha \over 2} \right) \, |u_{ij}| \, \sigma.
$$ 

The associated forward Kolmogorov equation
on the probability law $f^N_t$ of $(\VV^N_t)$ in $\R^{dN}$ reads
\begin{equation}\label{eq:MixKolmo}
\partial_t \langle f^N_t,\varphi \rangle = \langle f^N_t, G^N \varphi \rangle
\end{equation} 
with generator $G^N = G^N_1 + G^N_2$, where $G^N_1$ is associated to
an \emph{inelastic} Boltzmann collision process whose collision kernel
only depends on the deviation angle as in Section~\ref{sec:BddBoltzmann}
\begin{equation}\label{def:BbddGN1}
  (G^N_1\varphi) (V) = {1 \over N} \, 
  \sum_{i,j= 1}^N  \int_{\mathbb{S}^{d-1}} b(\cos \theta_{ij}) \, 
  \left[\varphi(V^*_{ij}) - \varphi(V)\right]\, {\rm d}\sigma,
\end{equation}
with $\cos \theta_{ij} = \sigma \cdot (v_j-v_i)/|v_j-v_i|$ and
$V^*_{ij}$ defined in \eqref{eq:v*iv*j}, and $G^N_2$ is the
generator associated to the Brownian motion
\begin{equation}\label{def:MixG2}
(G^N_2\varphi) (V) =   \sum_{i=1}^N \Delta_{i} \varphi,
\end{equation}
where $v_i:=(v_{i,1}, \dots, v_{i,d})$ and $\Delta_i$ 
denotes the Laplacian  in $\R^d$
associated to the $i$-th particle:
\[
\Delta_i := \sum_{\alpha=1} ^d \partial^2_{v_{i,\alpha},v_{i,\alpha}}.
\]
It is classical to prove that $\VV^N_t$ is a Feller process and we
refer to the textbooks \cite{EthierKurtz,Pazy} where the theory is set
up with full details (one can also refer to
\cite{S6,Meleard1996,FournierGodinho} where similar processes are
considered).

\subsection{Statement of the result}
\label{sec:resultThermostat}

The main result in this section is a quantitative estimate of
propagation of chaos for the mixed collision and diffusion model
introduced above.
\begin{theo}\label{theo:Mix} 
  Consider an initial distribution $f_{\mbox{{\tiny in}}} \in
  \PP_q(\R^d)$, $q \ge 2$, the hierarchy of $N$-particle distributions
  $f^N _t = S^N _t(f_{\mbox{{\tiny \emph{in}}}} ^{\otimes N})$ following the
  evolution \eqref{eq:MixKolmo}, and the nonlinear semigroup $f_t =
  S^{N\!L} _t(f_{\mbox{{\tiny \emph{in}}}})$ following the evolution
  \eqref{eq:nlMix}.
 
  Then there is a constant $C>0$ and, for any $T>0$, there are
  constants $C_{T}, \tilde C_T \in (0,\infty)$ only depending $T \in
  (0,\infty)$ such that for any
  \[ 
  \varphi = \varphi_1 \otimes  \dots \otimes \,
  \varphi_\ell \in \FF^{\otimes \ell}, \quad \FF := W^{9,1}(\R^d) \cap
W^{1,\infty} (\R^d), \quad \| \varphi_j \|_\FF \le 1,
  \] 
  we have for $N \ge 2 \ell$:
  \begin{equation} \label{eq:cvgmix} \sup_{[0,T]}\left| \left \langle
        \left( S^N_t(f_{\mbox{{\tiny \emph{in}}}} ^N) - \left( S^{N\!L} _t (f_{\mbox{{\tiny \emph{in}}}})
          \right)^{\otimes N} \right), \varphi \right\rangle \right|
    \le C \, {\ell^2 \over N} + C_{T} \, {\ell^2 \over N}+ \tilde
    C_{T} \, \ell \, \Omega^{W_2} _N (f_{\mbox{{\tiny \emph{in}}}}) .
\end{equation}
As a consequence of \eqref{eq:cvgmix} and Lemma~\ref{lem:Rachev}, this
shows the quantitative propagation of chaos with rate $\eps(N) \le
C(\ell,T,f_{\mbox{{\tiny \emph{in}}}}) \, N^{-{1 \over d+4}}$ for any
initial data $f_{\mbox{{\tiny \emph{in}}}} \in \PP_{d+5}(\R^d)$.
\end{theo}

We are not aware of any result of propatation of chaos in this
setting. A conceivable alternative approach would be to use the
general nonlinear martingale approach, but that would most likely not
provide any quantitative rate of propagation of chaos. The techniques
developed recently in \cite{FournierGodinho} for the elastic Kac
equation without cut-off is yet another alternative technique that
could be tried on this model, but we have not made any attempts in
this direction, and, would it work, it is not clear as to what kind of
convergence rate one could hope to achieve.

\subsection{Proof of Theorem~\ref{theo:Mix}}
\label{sec:proofThermo}
We shall prove that Theorem~\ref{theo:Mix} is a consequence of
Theorem~\ref{theo:abstract} by proving that the assumptions {\bf
  (A1)-(A2)-(A3)(A4)-(A5)}. We consider the phase space $E=\R^d$ and
the following choice of functional spaces
$$
\left\{
\begin{array}{l} 
\GG_1 := \HH^{-s_1}(\R^d), \,\, s_1 := 3, \vspace{0.2cm} \\
\GG_2 := \HH^{-s_2}(\R^d), \,\, s_2 := 3 s_1=9, \vspace{0.2cm} \\
\FF_1 = W^{s_1,1}(\R^d), \vspace{0.2cm} \\
\FF_2 = W^{s_2,1}(\R^d),
\end{array}
\right.
$$
where the Fourier based space $\HH^{-s}(\R^d)$ and the norms
$|\cdot|_s$ are defined in example~\ref{expleFourier}, and the
corresponding spaces $\PP_{\GG_1}(E)$ and $\PP_{\GG_2}(E)$ (without
weight). We finally define $\FF_3 = \mbox{Lip}(\R^d)$ and
$\PP_{\GG_3}(E) := \PP_2(\R^d)$ endowed with the quadratic MKW
distance $W_2$.

\medskip\noindent {\bf Proof of (A1). } The well-posedness of equation
\eqref{eq:MixKolmo}-\eqref{def:BbddGN1} is a variation on the
well-posedness result for equation \eqref{eq:nlMix} as obtained in
\cite{BisiCT,BolleyC}.  We also refer to
\cite{EthierKurtz,S6,Meleard1996,FournierGodinho} for a proof of the
fact that t $\ZZ^N_t$ is a Feller process.

\medskip\noindent {\bf Proof of (A2).} First we prove {\bf (A2)-(i)}, and
more precisely we prove that
\[
S^{N\!L}_t \in C^{0,1}(\PP_{\GG_1}(E),\PP_{\GG_1}(E)),
\]
which is a consequence of the following result:
\begin{lem} \label{lem:contraction} For any $f_{\mbox{{\tiny
        \emph{in}}}},g_{\mbox{{\tiny \emph{in}}}} \in \PP(\R^d)$ and
  any final time $T >0$, the associated solutions $f_t$ and $g_t$ to
  the diffusive inelastic Boltzmann equation \eqref{eq:nlMix} satisfy
  for any $s\ge0$
\begin{equation}\label{estim:C01TrueMax}
  \sup_{t \in [0,T]} \left|f_t-g_t\right|_s \le e^{2T}\, 
  \left|f_{\mbox{{\tiny \emph{in}}}}-g_{\mbox{{\tiny \emph{in}}}}\right|_s.
\end{equation}
\end{lem}

\noindent
{\it Proof of Lemma~\ref{lem:contraction}.} 
We recall Bobylev's identity for Maxwellian inelastic collision kernel
(see for instance \cite{BisiCT})
$$
\FF\left(Q^+ _\alpha(f,g)\right) (\xi) = \hat Q^+ _\alpha(F,G) (\xi) =: {1
  \over 2} \int_{\Sp^{d-1}} b\left(\sigma \cdot \hat\xi\right) \, [F^+
\, G^- + F^- \, G^+ ]\, {\rm d}\sigma,
$$
with $F = \hat f$, $G = \hat g$, $F^\pm= F(\xi^\pm)$, $G^\pm=
G(\xi^\pm)$
and
$$
\xi^+ = {3-\alpha \over 4} \, \xi +   {1+\alpha \over 4}  \, |\xi| \, \sigma,
 \quad
\xi^- = {1+\alpha \over 4} (\xi -  |\xi| \, \sigma).
$$
Denoting $D = \hat g - \hat f$, $S = \hat g + \hat f$, the following
equation holds
\begin{equation}\label{eq:BoltzMaxD}
  \partial_t D = \int_{S^2} b
  \left(\sigma \cdot \hat \xi \right) \, 
  \left[ \frac{D^+ \, S^-}{2} + \frac{D^- \, S^+}{2} \right] \, {\rm d}\sigma - D - |\xi|^2 \, D.
\end{equation}
Using that $\| S \|_\infty \le 2$ and then $|\xi^\pm|\le |\xi|$, we deduce in distributional sense
\begin{eqnarray*}
\frac{{\rm d}}{{\rm d}t} {|D| \over\langle \xi \rangle^s}  
&\le &
\left( \sup_{\xi \in \R^d} {|D| \over \langle \xi \rangle^s} \right) \, \left( \sup_{\xi \in \R^d}
\int_{\Sp^{d-1}} b(\sigma \cdot \hat \xi )  \, \left\{ {\langle \xi^+ \rangle^s \over \langle \xi \rangle^s}
+ {\langle \xi^- \rangle^s \over \langle \xi \rangle^s} \right\} \, {\rm d}\sigma \right)
\\
&\le&
2 \, \sup_{\xi \in \R^d}  {|D| \over\langle \xi \rangle^s},
\end{eqnarray*}
from which we conclude that \eqref{estim:C01TrueMax} holds. 
\qed  

\medskip
Next we prove {\bf (A2)-(ii)}, as a consequence of the following result:
\begin{lem}
  \label{lem:estimQtoscani1}
  For any $f,g \in \PP(\R^d)$ and $s \ge 0$, we have
  \begin{equation}\label{eq:Qtoscani1}
    \left| Q_\alpha(f,f) \right|_s \le 2 
      \end{equation}
   and 
     \begin{equation}\label{eq:Qbiltoscani2}
    \left| Q_\alpha(f+g,f-g) \right|_s \le  3 \,  | f-g |_s .
  \end{equation}

  Moreover for any $s > 2$ there exists $\delta \in (0,1)$ such that
 \begin{equation}\label{eq:Deltatoscani}
    \left| \Delta f - \Delta g \right|_s \le  2 \,  | f-g |^\delta_s .
  \end{equation}
\end{lem}

\noindent
{\it Proof of Lemma~\ref{lem:contraction}.} 
We prove the second inequalities \eqref{eq:Qbiltoscani2}. We write in
Fourier:
\begin{multline*}
\mathcal{F}\left( Q_\alpha(f+g,f-g) \right) = \hat Q_\alpha(D,S) \\ =
\frac12 \, \int_{\mathbb{S}^{d-1}} b(\sigma \cdot \hat \xi) \, \left(
  S(\xi^+) \, D(\xi^-) + S(\xi^-) \, D(\xi^+) - 2 \, D(\xi) \right) \,
{\rm d}\sigma
\end{multline*}
where $\hat Q_\alpha$ is the Fourier transform of the symmetric
version of the collision operator $Q_\alpha$, which yields
\[
\frac{\left| \hat Q_\alpha (D,S) \right|}{\langle \xi \rangle^s} \le
\TT_1 + \TT_2 + \TT_3,
\]
with 
\begin{eqnarray*}
\TT_1 
&:=& \left|  \frac{1}{2 \, {\langle \xi \rangle^s}}  \,
\int_{\mathbb{S}^{d-1}} b(\sigma \cdot \hat \xi) \,  S(\xi^+) \,
  D(\xi^-) \, 
{\rm d}\sigma\right|
\\
&\le&  \int_{\mathbb{S}^{d-1}} b(\sigma \cdot \hat \xi) \, 
{\left|S(\xi^+)\right| \over 2}\, \frac{\left| D(\xi^-) \right|}{\langle \xi^- \rangle^s} \,
\frac{\langle \xi^- \rangle^s}{\langle \xi \rangle^s} \, {\rm d}\sigma 
\le  | D |_s.
\end{eqnarray*}

Similar estimates hold for the two other terms $\TT_2$ and $\TT_3$.
The proof of the first inequality \eqref{eq:Qtoscani1} is similar (and
simpler): we use the Fourier representation of $Q_\alpha(f,f)$ and the
bound $\| \hat f \|_{L^\infty} \le 1$.  We finally prove the last
inequality. We compute
$$
\left| \Delta f - \Delta g \right|_s = \sup_{\xi \in \R^d} |\xi|^2 {|F
  - G|\over \langle \xi \rangle^s} \le \sup_{\xi \in \R^d} \left( |F -
G|^{1-\delta} \, \left( {|F - G|\over \langle \xi
    \rangle^s}\right)^{\delta} \right) 
 $$
 with $\delta := (s-2)/s$.
\qed

\medskip\noindent {\bf Proof of (A3).}  We claim that for any $s_1 \ge
3$ there exists $C_1\in \R_+$ such that for all $\Phi \in
C^{2,1}(\PP_{\GG_1}(E),\R)$
 \begin{equation}
    \label{eq:H1VMix}
    \left\|  G^N (\Phi \circ \mu^N_Z ) -  \left\langle
        Q(\mu^N_Z,\mu^N_Z), 
        D\Phi[\mu^N_Z] \right\rangle
    \right\|_{L^\infty(E^N)} \le {C_1 \over N}   \| \Phi \|_{C^{2,1}(\PP_{\GG_1}(E),\R)},
  \end{equation} 
  which is {\bf (A3)} with $k=2$, $ \eta=1$ and $\eps(N) = C_{1} \,
  N^{-1}$.

  \smallskip We begin with a technical lemma which shows that the norm
  $| \cdot |_s$ is well-adapted for obtaining differentiability of the
  empirical measures.  It is worth emphasizing that the choice of $s_1
  = 3$ (in fact we only need $s_1 > 2$ by modifying slightly the
  arguments) comes from the requirement  that the function $V \mapsto \Phi
  (\mu^N_Z)$ be $C^2$.

\begin{lem}\label{lem:muvHHs} 
  The map  $\R^{dN}\to \PP_{\GG_1}(E)$, $V
  \mapsto \mu^N_Z$ is $C^{2,1}$ and
  \[
  \partial_{(v_i)_\alpha} (\mu^N_Z) = \frac1N \, \partial_{\alpha}
  \delta_{v_i}, \quad 
  \partial^2_{(v_i)_\alpha,(v_i)_\beta} (\mu^N_Z) = N^{-1}
  \, \partial^2_{\alpha\beta} \delta_{v_i}.
\]
\end{lem}

\noindent {\sl Proof of Lemma~\ref{lem:muvHHs}.} 
For $v,w \in \R^d$, we have
\begin{eqnarray*}
| \delta_v - \delta_w |_s 
&=& \sup_{\xi \in \R^d} { \left| e^{- i \, v \cdot \xi} - e^{- i \, w\cdot \xi}  \right|\over \langle \xi \rangle^s }
\le
|v-w| \,  \sup_{\xi \in \R^d} { \left\| \nabla _v e^{- i \, v \cdot \xi}    \right\|_{L^\infty(\R^d_v)} \over \langle \xi \rangle^s }
\\
&\le&|v-w| \,  \sup_{\xi \in \R^d} { |\xi|  \over \langle \xi \rangle^s }   \le |v-w|
\end{eqnarray*}
which shows that $v \mapsto \delta_v$ is $C^{0,1}$. For the sake of simplicity we present the proof of differentiability
when $d=1$, the case $d>1$ being similar.
For $v \in \R$ and $h \in \R^*$, we have
\[
  \left|\delta_{v+h} - \delta_v - h \, \delta'_v\right|_s
  = \sup_{\xi \in \R} { \left| (e^{-i \, \xi \, h} - 1 
      + i \, \xi \, h) \, e^{-i \, v \, \xi} \right| \over \langle \xi \rangle^s}
  \le \sup_{\xi \in \R} { \left| \xi \, h\right|^2  
    \over \langle \xi \rangle^s} \le |h|^2,
\]
from which we deduce that $v \mapsto \delta_v$ is $C^{1,1}$. 
Similarly we can go to second order:
\begin{multline*}
  \left|\delta_{v+h} - \delta_v - h \, \delta'_v + \frac{h^2}2 \,
    \delta''_v \right|_s \\
  = \sup_{\xi \in \R} { \left| \left(e^{-i \, \xi \, h} - 1 
      + i \, \xi \, h - \xi^2 \, h^2\right) \, e^{-i \, v \, \xi} 
  \right| \over \langle \xi \rangle^s}
  \le \sup_{\xi \in \R} { \left| \xi \, h\right|^3  
    \over \langle \xi \rangle^s} \le |h|^3,
\end{multline*}
and we easily conclude that $v \mapsto \delta_v$ is $C^{2,1}$. When
the dimension $d$ is greater than $1$, one can perform the
same argument for the partial derivatives of the Dirac mass.\qed

\smallskip
We come back to the proof of \eqref {eq:H1VMix}. Take $\Phi \in
C^{2,1}(\PP_{\GG_1}(E), \R)$ and compute separately the contributions of
$G^N_i$, $i=1, 2$. Proceeding as in the proof of {\bf (A3)} in
Theorem~\ref{theo:BddBoltz} we have
\[
G^N_1 \left(\Phi \circ \mu^N_Z \right) = \left\langle
  Q_\alpha\left(\mu^N_Z,\mu^N_Z\right), D \Phi(\mu^N_Z) \right\rangle
+ I_2(V)
 \]
 with 
 \begin{multline*}
 |I_2(V)|
 \le {1 \over 2N} \, \sum_{i,j= 1}^N  \int_{\mathbb{S}^{d-1}} b(\cos(\theta_{ij}) )
   \|\Phi \|_{C^{2,1}} \,  \left| \mu^N_{V^*_{ij}} -\mu^N_Z \right|_{s_1} ^{2} \,
   {\rm d}\sigma \\ \le {8 \over N} \,  \|\Phi \|_{C^{2,1}(\PP_{\GG_1}(E), \R)} ,
 \end{multline*}
 since for any $i \not = j$
  \begin{eqnarray*}
\left| \mu^N_{V^*_{ij}} -\mu^N_Z \right|_{s_1}
 &=&  
 {1 \over N} \, \left| \delta_{v_i'} +  \delta_{v_j'} -  \delta_{v_i} -  \delta_{v_j} \right|_{s_1}
 \\
 &\le&  
 {1 \over N} \, \left( \left| \delta_{v_i'}  \right|_{s_1} + \left|
     \delta_{v_j'} \right|_{s_1} + 
\left|  \delta_{v_i}  \right|_{s_1} + \left| \delta_{v_j}
\right|_{s_1} \right) = {4 \over N}. 
  \end{eqnarray*}

  On the other hand, as in the proof of assumption {\bf (A3)} in
  Section~\ref{sec:McK}, the map $\R^{dN} \to \R$, $V \mapsto
  \Phi(\mu^N_Z)$ is $C^{2,1}$ thanks to Lemma~\ref{lem:muvHHs} and 
  denoting $\phi_Z = D\Phi\!\left[\mu^N_Z\right] \in
  (\HH^{s_1}(\R^d))' $, we compute
\begin{eqnarray*}
 G^N_2 (\Phi(\mu^N_Z))
   &= &\sum_{i=1}^N \Delta_i \Phi(\mu^N_Z) 
   \\
  &=&  \sum_{i=1}^N  \left\{  \frac{1}{N} \,  (\Delta  \phi_Z) (v_i) 
  + \frac{1}{N^2} \sum_{\alpha=1}^d D^2\Phi\!\left[\mu^N_Z\right] 
  \left( \partial_{\alpha} \delta_{v_i},  
    \partial_{\alpha} \delta_{v_i} \right)  \right\}
\\
&=& \langle \Delta \mu^N_Z,   \phi_Z \rangle + \OO \left( \frac{\|\Phi \|_{C^{2,1}(\PP_{\GG_1}(E), \R)}}N \right).
  \end{eqnarray*}
We conclude the proof by combining the previous estimates. \qed


\medskip\noindent {\bf Proof of (A4). } For $f_{\mbox{{\tiny in}}}, g_{\mbox{{\tiny in}}} \in
\PP(\R^d)$, we define the associated solutions $f_t$ and $g_t$ to the
nonlinear Boltzmann equation; we define $h_t := \LL^{N\!
  L}_t[f_{\mbox{{\tiny in}}}](g_{\mbox{{\tiny in}}}-f_{\mbox{{\tiny in}}})$ the solution of the linearized
Boltzmann equation around $f_t$; and we define $r_t$ the solution to
the ``second variation'' equation around $f_t$. More precisely, we
define
\[
\left\{
\begin{array}{l}
\partial_t f_t = Q_\alpha(f_t,f_t) + \Delta \, f_t , \quad f_{|t=0} = f_{\mbox{{\tiny in}}} \vspace{0.3cm} \\
\partial_t g_t = Q_\alpha(g_t,g_t) + \Delta \, g_t , \quad g_{|t=0} = g_{\mbox{{\tiny in}}} \vspace{0.3cm} \\
\partial_t h_t = Q_\alpha(f_t,h_t) + Q_\alpha(h_t,f_t)  + \Delta \, h_t , \quad h_{|t=0} = h_{\mbox{{\tiny in}}}, 
\vspace{0.3cm} \\
\partial_t r_t =  Q_\alpha(f_t,r_t) + Q_\alpha(r_t,f_t) + \Delta \,
r_t + \frac12 
Q_\alpha(h_t,\tilde h_t) + \frac12 Q_\alpha(\tilde h_t, h_t), \quad r_{|t=0} = 0
\end{array}
\right.
\]
where in the last equation (second-order variation) $h_t$ and $\tilde
h_t$ are two solutions to the third equation (first-order variation). 

We then define when $h_{\mbox{{\tiny in}}} = \tilde h_{\mbox{{\tiny
      in}}} = g_{\mbox{{\tiny in}}} - f_{\mbox{{\tiny in}}}$ the
following error terms
\[
\left\{
\begin{array}{l}
  \mathsf d_t := g_t - f_t  \vspace{0.3cm} \\
  \omega_t := g_t - f_t - h_t = S^{N \! L}_t(g_{\mbox{{\tiny in}}}) - S^{N \! L}_t(f_{\mbox{{\tiny in}}})-
  \LL^{N \! L}_t[f_{\mbox{{\tiny in}}}] (g_{\mbox{{\tiny in}}} - f_{\mbox{{\tiny in}}}) \vspace{0.3cm} \\
  \psi_t := g_t - f_t - h_t - r_t.
\end{array}
\right.
\]

\begin{lem}\label{lem:a4maxwellfourier} Fix $s \ge 0$ and $T \in
  (0,\infty)$. There exists  $C_T$ such that for any  
$f_{\mbox{{\tiny \emph{in}}}}, g_{\mbox{{\tiny \emph{in}}}} \in \PP(\R^d)$, the following estimates hold
\begin{eqnarray}
\label{ineq;tMMC01}
&&  \forall \, t \in [0,T], \quad \left| h_t \right|_s \le C_T \, \left| h_{\mbox{{\tiny \emph{in}}}} \right|_s, 
\\
\label{ineq;tMMC02}
&&  \forall \, t \in [0,T], \quad  \left| r_t \right|_{2s}  \le C_T \,
\left| h_{\mbox{{\tiny \emph{in}}}} \right|_s \, | \tilde h_{\mbox{{\tiny \emph{in}}}} |_s, 
\end{eqnarray}
and when $h_{\mbox{{\tiny \emph{in}}}} = \tilde h_{\mbox{{\tiny \emph{in}}}} = g_{\mbox{{\tiny \emph{in}}}} - f_{\mbox{{\tiny \emph{in}}}}$ we have furthermore
\begin{eqnarray}
\label{ineq;tMMC11}
&&  \forall \, t \in [0,T], \quad  \left| \omega_t \right|_{2s} \le C_T \, \left| f_{\mbox{{\tiny \emph{in}}}} - g_{\mbox{{\tiny \emph{in}}}} \right|^2_s, 
\\
\label{ineq;tMMC21}
&&  \forall \, t \in [0,T], \quad  \left| \psi_t \right|_{3s} \le C_T \, \left| f_{\mbox{{\tiny \emph{in}}}} - g_{\mbox{{\tiny \emph{in}}}} \right|^3_s.
\end{eqnarray}
This proves that $S^{N\!L}_t \in C^{2,1}(\PP_{\GG_1}(E),P_{\GG_2}(E))$.
\end{lem}

\begin{proof}[Proof of Lemma~\ref{lem:a4maxwellfourier}]
  We skip the proof of \eqref{ineq;tMMC01} since it is similar to the
  proof of \eqref{estim:C01TrueMax}.  We then deal with each term
  successively. We work in Fourier variable and we introduce the
  notations $F = \hat f$, $D = \hat{\mathsf d}$, $H = \hat h$, $\tilde
  H = (\tilde h)^{\hat{}}$, $\Omega = \hat \omega$, $R = \hat r$ and $\Psi
  = \hat \psi$.

\smallskip\noindent
{\sl Step 1. } The evolution equation satisfied by $\Omega$ is
\begin{equation}\label{eq:BIevold}
  \partial_t \Omega= \hat Q_\alpha(\Omega,F) + \hat Q_\alpha(F,\Omega) - |\xi|^2 \, \Omega
  - \hat Q_\alpha(D,D).
\end{equation}
We deduce in distributional sense
$$
 \frac{{\rm d}}{{\rm d}t} {|\Omega (\xi)| \over \langle\xi\rangle^{2s}} 
\le \TT_1 + \TT_2,
$$
where
\begin{eqnarray*}
  \TT_1
  &:=& \sup_{\xi \in \R^d} \int_{\mathbb{S}^{d-1}} 
  {b\left(\sigma \cdot \hat\xi\right) \over \langle\xi\rangle^{2s}} \, 
  \Bigg( \left| \frac{\Omega (\xi^+) \, F (\xi^-)}{2} \right|  
    + \left|\frac{\Omega (\xi^-) \, F
        (\xi^+)}{2} \right| \\ 
    && \hspace{6cm} - F(\xi) \Omega(0) - F(0) \Omega(\xi) \Bigg) \, {\rm d}\sigma
  \\
  &\le& \sup_{\xi \in \R^d} \int_{\mathbb{S}^{d-1}} 
  b\left(\sigma \cdot \hat\xi\right) \, 
  \Bigg( {\left|\Omega (\xi^+)\right| \over \langle\xi^+\rangle^{2s}} \, 
   { \langle\xi^+\rangle^{2s} \over \langle\xi\rangle^{2s}}
    + {\left|\Omega (\xi^-)\right| \over \langle\xi^-\rangle^{2s}} \, 
    {\langle\xi^-\rangle^{2s} \over\langle\xi\rangle^{2s}} \\ 
    && \hspace{6cm} +
    \frac{|\Omega(\xi)|}{\langle \xi \rangle^{2s}} + |\Omega(0)| \,
    \frac{|F(\xi)|}{\langle \xi \rangle^{2s}} \Bigg) \, {\rm d}\sigma
  \\
  &\le& C \, \sup_{\xi \in \R^d} 
    {\left|\Omega (\xi)\right| \over \langle\xi\rangle^{2s}} + |\Omega(0)| \, \sup_{\xi \in \R^d} 
    {\left|F (\xi)\right| \over \langle\xi\rangle^{2s}},
\end{eqnarray*}
for some constant $C>0$, and
\begin{eqnarray*}
  \TT_2 &:=& {1\over 2} \,  \sup_{\xi \in \R^d} \int_{\mathbb{S}^{d-1}} 
  {b\left(\sigma \cdot \hat\xi\right) \over \langle\xi\rangle^{2s}} \, 
  \left| D (\xi^+) \, D (\xi^-) +  D (\xi^-) \, 
    D (\xi^+) \right|  \,  {\rm d}\sigma
  \\ 
  &\le& {1 \over 2} \, \sup_{\xi \in \R^d} \int_{\mathbb{S}^{d-1}} 
  b\left(\sigma \cdot \hat\xi\right)
  \, \left(  {| D (\xi^+) | \over \langle\xi^+\rangle^s} \, {| D (\xi^-) | \over \langle\xi^-\rangle^s}
   + {| D (\xi^+) |^2 \over\langle\xi^+\rangle^s} \,  {| D (\xi^-) |^2 \over \langle\xi^-\rangle^s}  \right) \, {\rm d}\sigma
  \\
  &\le& \left| \mathsf d_t \right|_s ^2 
  \le  C _T \, \left| f_{\mbox{{\tiny in}}} - g_{\mbox{{\tiny in}}}
  \right|_s ^2, 
\end{eqnarray*}
using the estimates \eqref{estim:C01TrueMax}. 
We then conclude thanks to a Gronwall lemma. 

\medskip\noindent
{\sl Step 2. } The evolution equation satisfied by $R$ is
\begin{equation}\label{eq:RIevold}
  \partial_t R = \hat Q_\alpha(F,R) + \hat Q_\alpha(R,F) - |\xi|^2 \,
  R + \frac12 \hat
  Q_\alpha(H,\tilde H)  + \frac12 \hat Q_\alpha(\tilde H, H), \quad R_{|t=0} = 0.
\end{equation}
Equation \eqref{eq:RIevold} being similar to equation
\eqref{eq:BIevold}, with the same computations as in Step 1 we deduce
that \eqref{ineq;tMMC02} holds.

\smallskip\noindent
{\sl Step 3. } Choosing now $h_{\mbox{{\tiny in}}} = \tilde
h_{\mbox{{\tiny in}}} = \mathsf d_{\mbox{{\tiny
      in}}}$, the equation satisfied by $\Psi$ is 
\[  
\partial_t \Psi = \hat Q_\alpha(F,\Psi) + \hat Q_\alpha(\Psi,F) - |\xi|^2 \, \Psi - \hat Q_\alpha(\Omega,H) 
  - \hat Q_\alpha(D,\Omega), \quad \Psi_{|t=0}  = 0.
\]
Observe that by conservation of mass $\Psi_t(0)=0$ for all times, and
with these choices of initial data $\hat Q_\alpha(\Omega,H) = \hat
Q_\alpha ^+(\Omega,H)$ and $\hat Q_\alpha(D,\Omega) = \hat
Q_\alpha ^+(D,\Omega)$.  Then we perform similar computations as in Step
1, and we deduce in distributional sense
$$
 \frac{{\rm d}}{{\rm d}t} {|\Psi (\xi)| \over \langle\xi\rangle^{3s}} 
\le \TT_1 + \TT_2 + \TT_3,
$$
where
\[
\TT_1 := \sup_{\xi \in \R^3} {| \hat Q_\alpha(F,\Psi) + \hat
  Q_\alpha(\Psi,F)|\over \langle\xi\rangle^{3s}} \le C \, \sup_{\xi \in
  \R^3} {\left|\Psi (\xi)\right| \over \langle\xi\rangle^{3s}}, 
\]
\[
  \TT_2
 := \sup_{\xi \in \R^3} {| \hat Q_\alpha^+(\Omega,H) |\over  \langle\xi\rangle^{3s}} 
 \le 2 \, \left( \sup_{\xi \in \R^3} 
    {\left| \Omega (\xi)\right| \over \langle\xi\rangle^{2s}} \right) \,
  \left( \sup_{\xi \in \R^3} 
    {\left| H(\xi)\right| \over \langle\xi\rangle^{s}} \right),
\]
\[
  \TT_3
 := \sup_{\xi \in \R^3} {| \hat Q^+_\alpha(D,\Omega) |\over  \langle\xi\rangle^{3s}} 
 \le 2 \, \left( \sup_{\xi \in \R^3} 
    {\left| D (\xi)\right| \over \langle\xi\rangle^{s}} \right) \,
  \left( \sup_{\xi \in \R^3} 
    {\left| \Omega (\xi)\right| \over \langle\xi\rangle^{2s}} \right).
\]
Finally we then conclude the proof of \eqref{ineq;tMMC21} using the
already established estimates \eqref{estim:C01TrueMax},
\eqref{ineq;tMMC01}, \eqref{ineq;tMMC11}, and the Gronwall lemma.
\end{proof}

\medskip\noindent {\bf Proof of (A5). } We use the following result proved in
\cite{BolleyC} (see also \cite{BisiCT} for a similar result)
$$
\sup_{t \ge 0}W_2( S^{N\!L}_t f_{\mbox{{\tiny in}}}, S^{N\!L}_t
g_{\mbox{{\tiny in}}}) \le W_2 (f_{\mbox{{\tiny in}}},g_{\mbox{{\tiny
      in}}}),
$$
which concludes the proof. 


\bibliographystyle{acm}
\bibliography{./meanfield}


\signsm \signcm \signbw

\end{document}